\documentclass[11pt]{amsart}

\usepackage[utf8]{inputenc}
\usepackage{amsmath}
\usepackage{amssymb}
\usepackage{amsthm}
\usepackage{color}

\usepackage{hyperref}

\newtheorem{theorem}{Theorem}[section]
\newtheorem{corolario}[theorem]{Corollary}
\newtheorem{lema}[theorem]{Lemma}
\newtheorem{prop}[theorem]{Proposition}
\newtheorem{defn}[theorem]{Definition}
\theoremstyle{definition}
\newtheorem{remark}[theorem]{Remark}
\newtheorem{notation}[theorem]{Notation}
\newtheorem{ejemplo}[theorem]{Example}

\title{Computing a Saito basis from a standard basis}
\author{Felipe Cano}
\address{Felipe Cano. Departamento de Álgebra, Análisis Matemático,
Geometr\'\i a y Topolog\'\i a. Universidad de Valladolid.
Paseo de Bel\'en 7,
47011 -- Valladolid, SPAIN}
\email{fcano@uva.es}
\author{Nuria Corral}
\address{Nuria Corral. Departamento de Matemáticas, Estadística y Computación. Universidad de Cantabria. Avda. de los Castros s/n, 39005 -- Santander, SPAIN}
\email{nuria.corral@unican.es}
\author{David Senovilla-Sanz}
\address{David Senovilla-Sanz. Departamento de Matemáticas, Estadística y Computación. Universidad de Cantabria. Avda. de los Castros s/n, 39005 -- Santander, SPAIN}
\email{david.senovilla@unican.es}
\subjclass[2020]{14H20, 14H15, 32S65, 32S15, 32S05}
\keywords{Saito basis, analytic invariants, equisingularity, semimodule, cusp, standard basis, differential values}
\thanks{The authors are supported by the Spanish research project PID2022-139631NB-I00 funded by the Agencia Estatal de Investigación - Ministerio de Ciencia e Innovación. The third author is also supported by a predoctoral contract “Concepción Arenal” of the  Universidad de Cantabria}

\begin{document}
\maketitle
\begin{abstract}
In this paper we describe how to compute a Saito basis of a cusp, a plane curve with only one Puiseux pair. Moreover, the 1-forms of the Saito basis that we compute are characterized in terms of their divisorial orders associated to the ``cuspidal" divisor of the minimal reduction of singularities of the cusp.
We also introduce  a new family of analytic invariants for plane curves computed in terms of Saito bases.

%
\end{abstract}
\tableofcontents
\section{Introduction}\label{sec:intro}
Let $C$ be a plane curve in $(\mathbb{C}^2,\mathbf{0})$ and consider the
$\mathcal{O}_{\mathbb C^2,\mathbf 0}$-module $\Omega_{\mathbb C^2,\mathbf 0}^1[C]$ of holomorphic 1-forms    which have the curve $C$ invariant. In \cite{saito}, K. Saito proved that  $\Omega_{\mathbb C^2,\mathbf 0}^1[C]$ is a free module of rank 2 and a basis of $\Omega_{\mathbb C^2,\mathbf 0}^1[C]$ is called a {\em Saito basis} for the curve $C$. The objective of this paper is to compute a Saito basis for a curve $C$ with only one Puiseux pair. Moreover, the 1-forms of the Saito basis that we compute are characterized in terms of their critical values. Let us precise the statements.

Consider an irreducible plane curve $C$ in $(\mathbb{C}^2,\mathbf{0})$ and let $\phi(t)$ be a primitive parametrization of $C$. Given $h \in \mathcal{O}_{\mathbb C^2,\mathbf{0}}$, we denote $\nu_C(h)=\text{ord}_t(h (\phi(t)))$.  Recall that the semigroup $\Gamma$ of $C$, defined as
$$
\Gamma=\{\nu_C(h):h\in\mathcal{O}_{\mathbb C^2,\mathbf{0}}\},
$$
is equivalent to the equisingularity data of the curve $C$. Given a 1-form $\omega \in \Omega^1_{\mathbb C^2,\mathbf 0}$, we denote $\nu_C(\omega)=\text{ord}_t(\alpha(t))+1$ with $\phi^*(\omega)=\alpha(t)dt$.  The set of differential values $\Lambda$ of the curve $C$, defined as
$$
\Lambda=\{ \nu_C(\omega) \ : \ \omega \in  \Omega^1_{\mathbb C^2,\mathbf 0}\},
$$
is a $\Gamma$-semimodule and $\Lambda$ is an analytic invariant of the curve $C$.
There exists a basis of $\Lambda$, that is, a strictly increasing sequence $\mathcal{B}=(\lambda_{-1},\lambda_0,\lambda_1, \ldots, \lambda_s)$ of elements of $\Lambda$, with $s$ minimal, such that
$$
\Lambda=\bigcup_{i=-1}^s (\lambda_i + \Gamma).
$$
A set of $1$-forms $(\omega_{-1},\omega_0,\omega_1,\ldots, \omega_s)$ such that $\nu_C(\omega_i)=\lambda_i$, for $-1 \leq i\leq s$, is called a {\em minimal standard basis} of the curve $C$.

Assume now that $C$ is a cusp, that is, an irreducible curve with a single Puiseux pair $(n,m)$. In this situation, the semigroup of the curve $C$ is equal to  $\Gamma=n \mathbb{Z}_{\geq 0} + m \mathbb{Z}_{\geq 0}$, and we say that $\Gamma$ is cuspidal and that $\Lambda$ is a cuspidal semimodule.  Let us introduce some structural values associated to $\Lambda$.

The basis of the semimodule $\Lambda$ allows to define a chain
$\Lambda_{-1}\subset \Lambda_0 \subset \Lambda_1 \subset \cdots \subset \Lambda_s=\Lambda$
with $\Lambda_{i}=\cup_{k=1}^{i} (\lambda_k + \Gamma)$, for $i=-1,0,1,\ldots,s$, such that $\lambda_i \not \in \Lambda_{i-1}$ for $i=0,1, \ldots,s$.
The {\em axes} $u_i^n$, $u^m_i$, $u_i$ and $\tilde{u}_i$ of $\Gamma$ are defined as
$$
u_i^n=\min \{\lambda_{i-1} + n \ell \in \Lambda_{i-2} \ : \ \ell \geq 1\}, \qquad u_i^{m}=\min \{\lambda_{i-1} + m\ell \in \Lambda_{i-2} \ : \ \ell \geq 1\},
$$
$$
u_i= \min \{u_i^n,u_i^m\}, \qquad \tilde{u}_i=\max\{u_i^n,u_i^m\},
$$
with $1 \leq i \leq s+1$, and the {\em critical values} $t_i^n, t_i^m,t_i,\tilde{t}_i$ are given by $t_{-1}=\lambda_{-1}=n$, $t_{0}=\lambda_0=m$ and
$$
t_i^n=t_{i-1} + u_i^n-\lambda_{i-1}, \qquad t_i^m=t_{i-1} + u_i^m-\lambda_{i-1},
$$
$$
t_i=\min\{t_i^n,t_i^m\}, \qquad \tilde{t}_i = \max\{t_i^n,t_i^m\},
$$
for $1 \leq i \leq s+1$.
Note that the semimodule $\Lambda$ of a cusp is increasing which means that $\lambda_i > u_i$ for $1 \leq i \leq s$ (see \cite{delorme}).
In a previous work  \cite{Can-C-SS-2023}, we have proved that these values allow to characterized the elements of an {\em extended standard basis} of the curve $C$, that is, a set of 1-forms $\mathcal{E}=(\omega_{-1},\omega_0,\omega_1,\ldots, \omega_s,\omega_{s+1})$ such that $\omega_{-1},\omega_0,\omega_1,\ldots, \omega_s$ is a minimal standard basis of $C$ and $\omega_{s+1}$ is a 1-form with $C$ as invariant curve ($\nu_C(\omega_{s+1})=\infty$) and  divisorial order with respect to the cuspidal divisor equal to $t_{s+1}$.

Let us recall the notion of divisorial order. Consider any sequence $\pi: M \to (\mathbb{C}^2,\mathbf{0})$ of punctual blowing-ups and let $D \subset \pi^{-1}(\mathbf{0})$ be an irreducible component of the exceptional divisor. Given a point $Q \in D$, we can take coordinates $(u,v)$ at $Q$ such that $D=(u=0)$ and, for any  $h \in \mathcal{O}_{\mathbb C^2,\mathbf{0}}$, we can write $f \circ \pi = u^\beta \tilde{h}$,  locally at $Q$, such that $u$ does not divide $\tilde{h}$. We define the {\em divisorial order} $\nu_D(h)$ {\em with respect to the divisor} $D$ by $\nu_D(h)=\beta$. Given a 1-form $\omega \in \Omega^1_{\mathbb C^2,\mathbf 0}$, that can be written as $\omega=a(x,y)dx+b(x,y)dy$, where $(x,y)$ are coordinates in $(\mathbb{C}^2,\mathbf{0})$, we define the {\em divisorial order} $\nu_D(\omega)$ of $\omega$  {\em with respect to the divisor} $D$ by $\nu_D(\omega)=\min\{\nu_D(xa),\nu_D(yb)\}$.

If $C$ is a cusp, we can consider the minimal reduction of singularities $\pi_C: M \to (\mathbb{C}^2,\mathbf{0})$ and we denote $D_C$ the ``cuspidal divisor'', that is, the only irreducible component of $\pi^{-1}_C(\mathbf{0})$ such that the strict transform of $C$ intersects $D_C$. In this situation, if $(x,y)$ are adapted coordinates for $C$, the divisorial order with respect to the divisor $D_C$ is a monomial order since it can be computed as $\nu_{D_C}(h)=\min \{n i + mj \ : \ h_{ij} \neq 0\}$ where $h(x,y)=\sum_{i,j \geq 0} h_{ij} x^i y^j$.

Now we can state the main result of this article.

\begin{theorem}\label{th:main}
  Let $C$ be a curve in $(\mathbb{C}^2,\mathbf{0})$ with only one Puiseux pair and  $\mathcal B=(\lambda_{-1},\lambda_{0},\linebreak \lambda_1, \ldots,\lambda_s)$ be the basis of the semimodule $\Lambda$ of differential values for $C$. There exist
  two $1$-forms $\omega_{s+1},\tilde\omega_{s+1}$ having $C$ as an invariant curve and such that
  $$\nu_{D_C}(\omega_{s+1})=t_{s+1}, \qquad \nu_{D_C}(\tilde\omega_{s+1})=\tilde t_{s+1},$$
  where  $t_{s+1}$ and $\tilde t_{s+1}$ are the last critical values of $\Lambda$.

    Moreover, for any pair of $1$-forms as above, the set $\{\omega_{s+1},\tilde \omega_{s+1}\}$ is a Saito basis for $C$.
\end{theorem}
The proof of the existence of the 1-forms in the theorem above is done in a constructive way.

The notion of divisorial order with respect to a divisor  allows to introduce a new analytic invariant of any plane curve (which is not necessarily a cusp). Given a divisor $D$ as above, we define the {\em the pair $(\mathfrak{s}_{D}(C),\widetilde{\mathfrak s}_{D}(C))$ of Saito multiplicities at $D$} by
\begin{eqnarray*}
\mathfrak s_{D}(C)&=&\min\{\nu_{D}(\omega); \; \omega \text{ belongs to a Saito basis of } C\}.\\
\widetilde{\mathfrak s}_{D}(C)&=&\max\{\nu_{D}(\omega);\; \omega \text{ belongs to a Saito basis of } C\}.
\end{eqnarray*}
Then the pair $(\mathfrak s_{D}(C),\widetilde{\mathfrak s}_{D}(C))$ is an analytic invariant of any plane curve $C$ in $(\mathbb{C}^2,\mathbf{0})$. In \cite[p. 8-9]{Genzmer2}, Y. Genzmer introduces an  analytic invariant of a curve $C$ directly related with  the pair  $(\mathfrak s_{D_1}(C),\widetilde{\mathfrak s}_{D_1}(C))$ associated to the divisor $D_1$ which appears after one blow-up. More precisely, he proves that the pair of multiplicities $(\nu_0(\omega),\nu_0(\tilde\omega))$, where $\{\omega,\tilde{\omega}\}$ is a Saito basis for $C$, with  $\nu_0(\omega) \leq \nu_0(\tilde\omega)$ and such that $\nu_0(\omega)+\nu_0(\tilde\omega)$ is maximal, is an analytic invariant for the curve $C$.
Note that the multiplicity at the origin of a 1-form $\omega$ can also be computed as  $\nu_0(\omega)=\nu_{D_1}(\omega)-1$.

Moreover, we prove that the pair $(\mathfrak s_{D_C}(C),\widetilde{\mathfrak s}_{D_C}(C))$ is determined in terms of the critical values of the semimodule. More precisely
\begin{theorem}\label{th:main2}
Let $C$ be a cusp in $(\mathbb{C}^2,\mathbf{0})$. Then we have that
$$
	(\mathfrak s_{D_C}(C),\widetilde{\mathfrak s}_{D_C}(C))=(t_{s+1},\tilde t_{s+1}),
$$
where $t_{s+1}$ and $\tilde t_{s+1}$ are the last critical values of the semimodule of differential values of the curve $C$.
\end{theorem}
However, we  give an example of two cusps $C_1$ and $C_2$ with the same semimodule of differential values but with $(\mathfrak s_{D_1}(C_1),\widetilde{\mathfrak s}_{D_1}(C_1)) \neq  (\mathfrak s_{D_1}(C_2),\widetilde{\mathfrak s}_{D_1}(C_2))$.

\medskip

The structure of the semimodule  $\Lambda$  of differential values of a plane curve plays a key role in the proofs of the results of this article. In \cite{zariski}, O. Zariski pointed the importance of the semimodule $\Lambda$ in the analytic classification of plane curves. The work of Zariski was based on the computation of parametrizations of irreducible plane curves as simple as possible. Following these ideas, C. Delorme described the structure of the semimodule of differential values for a curve with only one Puiseux pair (see \cite{delorme}). The complete analytic classification of irreducible plane curves was given by A. Hefez and M. E. Hernandes in 2011 (\cite{hefez2}, see also \cite{Hef-H-handbook-2021}). Recently, M. E. Hernandes and M. E. R. Hernandes give the analytic  classification of plane curves in the general case (\cite{Her-H-2024}).

Moreover, some analytic invariants of a plane curve can be computed from the semimodule $\Lambda$; for instance, the Milnor $\mu(C)$ and Tjurina number $\tau(C)$ of an irreducible plane curve $C$ can be computed as $\mu(C)-\tau(C)=\sharp (\Lambda \setminus \Gamma)$ (see \cite{berger,Zar-66}). In \cite{Gen-H-2020}, Y. Genzmer and M. E. Hernandes compute the difference $\mu(C) - \tau(C)$ for an irreducible plane curve $C$ when $C$ admits a Saito basis of a special kind called good Saito basis. In some recent works, Y. Genzmer describes invariants associated to the Saito module to compute the the generic dimension of the moduli space of a curve $C$ (see {\cite{Genzmer2,Gen-JSymComp-2022,Gen-2024}).

\medskip
The article is organized as follows. In Section~\ref{sec:semi} we describe the structure of a cuspidal increasing semimodule extending some results obtained in \cite{Can-C-SS-2023}. In Section~\ref{sec:bases} we generalize the decomposition given by Delorme in \cite{delorme} (see also \cite{Can-C-SS-2023}) of the elements of a minimal standard basis of a cusp.

The proof of the main result of the paper  is given in Section~\ref{sec:saito}.  We introduce the notion of a {\em standard system} $(\mathcal{E},\mathcal{F})$ for a cusp which is the data of an extended standard basis $\mathcal{E}=(\omega_{-1},\omega_0,\omega_1,\ldots, \omega_s,\omega_{s+1})$ and a family $\mathcal{F}=(\tilde{\omega}_1,\tilde{\omega}_2,\cdots,\tilde{\omega}_s,\tilde{\omega}_{s+1})$ of 1-forms  with divisorial order $\nu_{D_C}(\tilde{\omega}_k)=\tilde{t}_j$ and such that $C$ is an invariant curve of each $\tilde{\omega}_j$, $1 \leq j \leq s+1$. Given an {\em standard system} $(\mathcal{E},\mathcal{F})$, we have that
the set $\mathcal T=\{\omega_{s+1},\tilde \omega_1,\tilde \omega_2,\ldots,\tilde \omega_{s+1}\}$ is a generator system of the Saito $\mathcal O_{\mathbb C^2,\mathbf 0}$-module $\Omega^1_{\mathbb C^2,\mathbf 0}[C]$ (see Proposition~\ref{prop:generadoresSaito}). We show the existence of 1-forms $\omega_{s+1},\tilde{\omega}_{s+1}$  having $C$ as invariant curve  and with divisorial values $\nu_{D_C}(\omega_{s+1})=t_{s+1}$ and $\nu_{D_C}(\tilde \omega_{s+1})=\tilde t_{s+1}$. Then, given two 1-forms $\omega_{s+1}$ and $\tilde{\omega}_{s+1}$ as above, we prove that it is possible to construct a   {\em special standard system} that contains $\omega_{s+1}$ and $\tilde{\omega}_{s+1}$, that is, a standard system such that the 1-forms $\tilde{\omega}_i$, $1 \leq i\leq s$, can be written in terms of $\omega_{s+1}$ and $\tilde{\omega}_{s+1}$ (see Proposition~\ref{prop:existencia-special-system}). Finally, we prove that the generator system $\mathcal{T}$ of $\Omega^1_{\mathbb C^2,\mathbf 0}[C]$ can be reduced to obtain a basis with the properties given in Theorem~\ref{th:main}.

The last section of the paper is devoted to introduce the analytic invariants given by the pair of Saito multiplicities at a divisor and we prove Theorem~\ref{th:main2}. Finally, we describe the examples which show that the Saito pairs of multiplicities with respect to the first divisor are not determined by the semimodule of differential values of the curve.

\section{Structure of Cuspidal  Semimodules}\label{sec:semi}

In this section we enlarge the description given in \cite{Can-C-SS-2023} of the structure of a cuspidal increasing semimodule.

Take $\Gamma \subset \mathbb{Z}_{\geq 0}$ an additive numerical semigroup, that is, $\Gamma$ is a monoid such is generated by $\langle \overline \beta_0,\overline \beta_1, \ldots,\overline \beta_g\rangle$ with $gcd(\overline \beta_0,\ldots,\overline\beta_g)=1$.  A set $\Lambda\subset \mathbb{Z}_{\geq 0}$ is a $\Gamma$-\emph{semimodule}, if $\gamma+\lambda\in \Lambda$  for all $\gamma\in \Gamma$ and $\lambda\in \Lambda$. The {\em  basis} of a $\Gamma$-semimodule $\Lambda$ is the only increasing sequence
$$
\mathcal B=(\lambda_{-1},\lambda_0,\lambda_1,\ldots,\lambda_s)
$$
satisfying that $\Lambda=\cup_{i=-1}^s (\lambda_i+\Gamma)$ and that $\lambda_j\notin \Lambda_{j-1}$, for any $j=0,1,\ldots,s$, where $\Lambda_i=\cup_{k=-1}^{i} (\lambda_k+\Gamma)$, for $i=-1,0,1,\ldots,s$.
 The basis induces the \emph{decomposition chain} of $\Lambda$:
$$
\Gamma+\lambda_{-1}=\Lambda_{-1}\subset \Lambda_0\subset \Lambda_1\subset \ldots\subset \Lambda_s=\Lambda,
$$
where each $\Lambda_i$ is a semimodule with basis $\mathcal B_i=(\lambda_{-1},\lambda_0,\ldots,\lambda_i)$, for $i=-1,0,\ldots,s$. When $\lambda_{-1}=0$, we say that $\Lambda$ is a \emph{normalized} semimodule. The number $s$ above is called the \emph{length} of the semimodule $\Lambda$.

Denote by $n=\min(\Gamma\setminus\{0\})$, given the basis $\mathcal B=(\lambda_{-1},\lambda_0,\ldots,\lambda_s)$, we have that $\lambda_i\neq \lambda_j$ mod $n$. Hence, the length $s$ is bounded by $n-2$.

We say that a numerical semigroup $\Gamma$ is \emph{cuspidal} if it is generated by two coprime integer numbers $n$, with $2\leq n<m$. A $\Gamma$-semimodule $\Lambda$ is \emph{cuspidal} when $\Gamma$ is cuspidal. From now on, we fix a cuspidal semigroup $\Gamma$, and we denote $n<m$ its generators.

Any cuspidal semimodule $\Lambda$ has an element $c_\Lambda\in\Lambda$ which is the minimum one satisfying the property that for every integer $k\geq c_\Lambda$, we have that $k\in \Lambda$. The element $c_\Lambda$ is called the \emph{conductor} of $\Lambda$. In the case $\Lambda=\Gamma$, the conductor takes value $c_\Gamma=(n-1)(m-1)$. Furthermore, since $\Lambda_i\subset \Lambda_{i+1}$, then $c_{\Lambda_{i}}\geq c_{\Lambda_{i+1}}$. Note that if $\Lambda$ is normalized, we have that $\Lambda_{-1}=\Gamma$ and in this case we have that
$
c_\Lambda\leq c_\Gamma
$.

\subsection{Axes, limits and critical values}

Let us introduce here some structural values for a cuspidal semimodule $\Lambda$ of basis
$\mathcal B=(\lambda_{-1},\lambda_0,\ldots,\lambda_s)$.

For $1\leq i\leq s+1$, we define the \emph{axes} $u_i^n,u_i^m,u_i$ and $\tilde{u}_i$ of $\Lambda$, as follows:
\begin{itemize}
\item $u_i^n=\min\{\lambda_{i-1}+n\ell\in \Lambda_{i-2};\, \ell\geq 1\}$. We write $u_i^n=\lambda_{i-1}+n\ell_i^n$.
\item $u_i^m=\min\{\lambda_{i-1}+m\ell\in \Lambda_{i-2};\, \ell\geq 1\}$. Similarly, we  put $u_i^m=\lambda_{i-1}+m\ell_i^m$.
\item $u_i=\min\{u_i^n,u_i^m\}$ and $\tilde{u}_i=\max\{u_i^n,u_i^m\}$.
\end{itemize}
The numbers $\ell_i^n$ and $\ell_i^m$ are called \emph{limits} of $\Lambda$.
\begin{remark}
If we consider the semimodule $\Lambda'=\Lambda-\lambda$, the new basis and the axes are shifted by $\lambda$ and we obtain the same limits as for $\Lambda$. This is particularly interesting when $\lambda=\lambda_{-1}$ and hence $\Lambda'$ is a normalized semimodule.
\end{remark}
\begin{remark}
\label{rk:cotasparalimites}
Let us note that $1\leq \ell_i^m<n$ and that $1\leq \ell_i^n<m$.
To see this we can suppose that $\Lambda$ is normalized and thus $c_{\Lambda_j}\leq c_\Gamma=(n-1)(m-1)$ for any $j=-1,0,1,\ldots,s$. Assume that $\ell_i^m\geq n$, we have
$$
\lambda_{i-1}+m(\ell_i^m-1)\geq (n-1)m\geq c_\Gamma\geq c_{\Lambda_{i-2}}.
$$
Then $\lambda_{i-1}+m(\ell_i^m-1)\in \Lambda_{i-2}$ in contradiction with the minimality of $\ell_i^m$. A similar argument proves that $\ell_i^n<m$.
\end{remark}

\begin{remark}
Notice that $u_i^n\neq u_i^m$ for each index $1\leq i\leq s+1$. Indeed, if $u_i^n=u_i^m$, then $n\ell_i^n=m\ell_i^m$; given that $n$ and $m$ are coprime, then $mk=\ell_i^n$, for a positive integer $k$ and hence $\ell_i^n\geq m$ which is a contradiction.
\end{remark}

\begin{lema}
\label{lema:propiedadescrecimiento}
	Let $\Lambda$ be a cuspidal semimodule of length $s$. Take $1\leq i\leq s+1$.  If $\lambda_{i-1}+na+mb\in \Lambda_{i-2}$, where $a,b\in\mathbb{Z}_{\geq0}$, then either $a\geq \ell_i^n$ or $b\geq \ell_i^m$.
\end{lema}
\begin{proof}
	(See also \cite{Can-C-SS-2023}, Lemma 6.9). By definition, we have that:
	$$
	\lambda_{i-1}+na+mb=\lambda_r+nc+md,\quad r<i-1,
	$$
	where $c,d$ are non negative integers.
	We proceed by induction on $\alpha=ac+bd\geq0$. If $\alpha=0$, then $ac=bd=0$. This implies that $ab=0$, otherwise, $ab\neq 0$ and hence $c=d=0$, that is
	$$
	\lambda_{i-1}+na+mb=\lambda_r,
	$$
	which is a contradiction because $\lambda_{r}<\lambda_{i-1}$. Now if $a=0$ we end with the minimality of $\ell_i^m$ and, similarly, if $b=0$, we end by the minimality of $\ell_i^n$.
	
	Assume that $\alpha>0$. Then $ac\neq 0$ or $bd\neq 0$. If $ac\neq 0$, let us put $a'=a-1\geq 0$ and $c'=c-1\geq 0$. We have that:
	$$
	\lambda_{i-1}+na'+mb=\lambda_r+nc'+md.
	$$
	By induction we are done. We apply a similar argument if $bd\neq 0$.
\end{proof}
For $-1\leq i\leq s+1$, we define inductively the {\em critical values} $t_i^n,t_i^m,t_i$, and $\tilde{t}_i$, for  $-1\leq i\leq s+1$ by putting $t_{-1}=\lambda_{-1}$ and $t_0=\lambda_0$ and
\begin{equation}
\begin{array}{cc}
 	t_i^n=t_{i-1}+n\ell_i^n,& \quad t_i^m=t_{i-1}+m\ell_i^m,\\
  t_i=\min\{t_i^n,t_i^m\},& \quad \tilde{t}_i=\max\{t_i^n,t_i^m\},
 \end{array}\quad 1\leq i\leq s+1.
 \end{equation}

\begin{remark}
	\label{rk:relacionejesyvalorescriticos}
  Noting that $n\ell_i^n=u_i^n-\lambda_{i-1}$ and $m\ell_i^m=u_i^m-\lambda_{i-1}$, we see that:
  \begin{equation*}
  \begin{aligned}
  t_i^n& =t_{i-1}+(u_i^n-\lambda_{i-1}),\quad
  	& t_i^m & =t_{i-1}+(u_i^m-\lambda_{i-1}),\\
  	t_i & =t_{i-1}+(u_i-\lambda_{i-1}), \quad 
  	&\tilde t_i& =t_{i-1}+(\tilde u_i-\lambda_{i-1}).\\
  \end{aligned}
  \end{equation*}
\end{remark}

\begin{defn}
	We say that the cuspidal semimodule $\Lambda$ is \emph{increasing} if we have that $\lambda_i>u_i$, for any $ 1\leq i\leq s.$
\end{defn}
Notice that if $\Lambda$ is increasing, then each $\Lambda_i$ is also increasing, for $1\leq i\leq s$. The notion of increasing semimodule was introduced in \cite{patricio}.
 \begin{lema}
 	\label{lema:lambdasytes}
 	Let $\Lambda$ be an increasing cuspidal semimodule. For any index $1\leq i\leq s$, we have that $\lambda_i-\lambda_j>t_i-t_j$,  for  $-1\leq j<i$.
 \end{lema}
\begin{proof}
	(See also \cite{Can-C-SS-2023},  Lemma 7.10). By a telescopic argument, it is enough to prove the following statements:
	\begin{itemize}
		\item $\lambda_r-\lambda_{r-1}>t_r-t_{r-1}$, for $1\leq r\leq s$.
		\item $\lambda_0-\lambda_{-1}\geq t_0-t_{-1}$.
	\end{itemize}
	The second statement is straightforward, because $t_{-1}=\lambda_{-1}$ and $t_0=\lambda_0$. Let us prove that $\lambda_r-\lambda_{r-1}>t_r-t_{r-1}$, for $1\leq r\leq s$.
	
	The inequality $\lambda_r-\lambda_{r-1}>t_r-t_{r-1}$  is equivalent to:
	$$
	t_r=t_{r-1}+u_r-\lambda_{r-1}>t_r+u_r-\lambda_r
	$$
	and this is equivalent to say that $\lambda_r>u_r$. We are done because $\Lambda$ is increasing.
\end{proof}
\begin{corolario}
	\label{cor: ejesyvalorescriticos}
	 Let $\Lambda$ be an increasing cuspidal semimodule. For any $1\leq i\leq s$, we have that
	$$
	u_{i+1}^n>t_{i+1}^n\quad \text{ and }\quad u_{i+1}^m>t_{i+1}^m.
	$$
\end{corolario}
\begin{proof}
Recalling that $t^n_{i+1}=u_{i+1}^n-(\lambda_i-t_i)$, it is enough to prove that $\lambda_i-t_i>0$. In view of Lemma \ref{lema:lambdasytes} and putting $j=-1$, we have that $\lambda_i-t_i>\lambda_{-1}-t_{-1}=0$.
\end{proof}

\begin{remark}
	The property of being increasing is true for cuspidal semimodules corresponding to the differential values of analytic branches of cusps (see \cite[Lemma 12]{delorme} and \cite[Theorem 7.13]{Can-C-SS-2023}). For this reason, here, we mainly consider increasing cuspidal semimodules.
\end{remark}

\subsection{Circular Intervals}

The circular intervals we describe here are useful for understanding the distribution of the elements of an increasing cuspidal semimodule. The notion of circular interval was introduced in \cite{Can-C-SS-2023}. Let us note that we have fixed two natural numbers $2\leq n<m$, with $\gcd(n,m)=1$ along this manuscript. We are going to consider the unit circle $\mathbb S^1\subset\mathbb C$ as a clock with $n$-hours as we explain below.

Let $\varepsilon: \mathbb R\rightarrow \mathbb S^1$ be the map given by
$$
\varepsilon(t)=\exp\left(-\frac{2\pi t\sqrt{-1}}{n}\right).
$$
We define the {\em $n$-clock\/  $\mathbb S^1_n$} to be
$\mathbb S^1_n=\varepsilon(\mathbb Z)$. Note that there is a bijection
$$
c:\mathbb Z/n\mathbb Z\rightarrow \mathbb S^1_n
$$
given by $c(k+n\mathbb Z)=\varepsilon(k)$. More than that, the bijection $c$ is an isomorphism of abelian groups, where $\mathbb S^1_n\subset \mathbb C$ has the induced multiplicative structure coming from the complex numbers $\mathbb C$. Note that
$$
c((k+n\mathbb Z)+ (k'+n\mathbb Z))=\varepsilon(k)\varepsilon(k').
$$
In particular $c(k+1+n\mathbb Z)= \varepsilon(k)\varepsilon(1)$.
\begin{notation} In order to visualize in a better way the arithmetic of the abelian multiplicative group $\mathbb S^1_n$, we introduce the following notations:
$$
\varepsilon(k)=k_\varepsilon,\quad \varepsilon(k)\varepsilon(k')=k_\varepsilon+k'_\varepsilon.
$$
Note that there is no confusion possible with the addition in $\mathbb C$. For instance, we have $(-1)_\varepsilon=(n-1)_\varepsilon$, $(k+1)_\varepsilon=k_\varepsilon+1_\varepsilon$ and $(k-1)_\varepsilon=k_\varepsilon-1_\varepsilon=
k_\varepsilon+(n-1)_\varepsilon$.
\end{notation}

Let us consider two points $P,Q\in \mathbb S^1_n$.  There are $\alpha\in \mathbb Z$ and  an integer number $\beta$ with $0\leq \beta\leq n-1 $ such that $P=\varepsilon(\alpha)$ and $Q=\varepsilon(\alpha+\beta)$. This number $\beta$, with $0\leq \beta\leq n-1$, does not depend on the chosen $\alpha$ such that $P=\varepsilon(\alpha)$ and we call it the {\em separation $S(P,Q)$ from $P$ to $Q$}. That is, if $P=\varepsilon(\alpha)$, we have that $Q=\varepsilon(\alpha+S(P,Q))$. We have that $S(P,P)=0$  and that
$$
S(P,Q)+S(Q,P)=n,\; \text{ if } Q\ne P.
$$
We define the {\em circular interval $<P,Q>$} to be
$$
<P,Q>=\{\varepsilon(\alpha+k);\; k=0,1,\ldots,S(P,Q)\}\subset \mathbb S^1_n.
$$
Note that if $P\ne Q$, we have that
$$
<P,Q>\cup<Q,P>=\mathbb S^1_n,\quad <P,Q>\cap<Q,P>=\{P,Q\}.
$$

\begin{remark} Take three points $P,Q, R\in\mathbb S^1_n$ with $P\ne Q$ and such that  $$R\in <P,Q>.$$
	 Then, we have that
	$
	S(P,Q)=S(P,R)+S(R,Q)\leq n-1
	$.
\end{remark}

Consider a list $B=(z_{-1},z_0,z_1,\ldots,z_s)$ of two by two distinct points $z_j\in\mathbb S^1_n$, with $s\geq 0$. For any index $0\leq i\leq s$ we define the {\em $i$-left bound $b_i^{\ell}(B)$ and the $i$-right bound ${b_i^{r}}(B)$ of $B$} to be integer numbers such that
$$
-1\leq b_i^{\ell}(B),b_i^r(B)\leq i-1
$$
and, moreover, the following holds:
\begin{enumerate}
\item If $k=b_i^{\ell}(B)$, then $S(z_{k},z_i)\leq S(z_q,z_i), \text{ for any } -1\leq q\leq i-1.$
\item If $\tilde k=b_i^{r}(B)$, then $S(z_i,z_{\tilde k})\leq S(z_i,z_q), \text{ for any } -1\leq q\leq i-1.$
\end{enumerate}
\begin{remark} Denote $k=b_i^{\ell}(B)$ and $\tilde k=b_i^{r}(B)$. The bounds are the integer numbers $k,\tilde k$ with $-1\leq k,\tilde k\leq i-1$  defined by the two following properties:
\begin{enumerate}
\item $z_i\in <z_k,z_{\tilde k}>$.
\item If $z_j\in <z_k,z_{\tilde k}>$ with $-1\leq j\leq i$, then $j\in \{i,k,\tilde k\}$.
\end{enumerate}
\end{remark}

\subsection{Circular Intervals in a Cuspidal Semimodule}
Let us recall that $\Gamma$ is the semigroup generated by $n,m$, with $2\leq n<m$ and $n,m$ are without common factors.

 Let  $\rho: \mathbb{Z}\rightarrow \mathbb{Z}/n\mathbb{Z}$ be the quotient map, that we also denote by $\rho(k)=\overline{k}$. Since $(n,m)=1$, the class $\bar m$ is a unit in $\mathbb Z/n\mathbb Z$, thus, we have a ring isomorphism
 $$
 \xi:\mathbb Z/n\mathbb Z\rightarrow \mathbb Z/n\mathbb Z,\quad  \bar m \xi(\bar k)=\xi(\bar m\bar k)=\bar k.
 $$
 Along this paper, we consider the map
 $\zeta:\mathbb Z\rightarrow\mathbb S^1_n$ defined by
 $
 \zeta(k)=(c\circ\xi\circ\rho)(k)
 $.
 Let us note that $\zeta(k+an)=\zeta(k)$ and that $\zeta(mk)=\varepsilon(k)=k_\varepsilon$.

Consider the intervals $I_q=\{nq,nq+1,\ldots,nq+n-1\}\subset \mathbb Z, $ $q\in\mathbb Z$. For a set $S\subset \mathbb{Z}$, we define the \emph{q-level set $R_q(S)$} by
$$
R_q(S)=\zeta(S\cap I_q)\subset \mathbb S^1_n.
$$
\begin{remark}If $S\subset \mathbb Z$ satisfies the property that
$ n+S\cap I_{q-1}\subset S\cap I_q$,
 we have that $R_{q-1}(S)\subset R_{q}(S)$. This is the case of cuspidal semimodules.
\end{remark}

Let us consider a cuspidal semimodule $\Lambda$ of length $s\geq 0$ with basis
$$
\mathcal B=(\lambda_{-1},\lambda_0,\lambda_1,\ldots,\lambda_s).
$$
We see the basis $\mathcal B$ in the clock $\mathbb S^1_n$ as
$
B=\zeta(\mathcal B)=(z_{-1},z_0,z_1,\ldots,z_s)$, where we have that $z_j=\zeta(\lambda_j)$, for $j=-1,0,1,\ldots,s$.

 Note that $z_i\ne z_j$, if $i\ne j$; indeed, saying that $z_i=z_j$ means that $\lambda_i-\lambda_j\in n\mathbb Z$, that is not possible in view of the definition of basis.

 Take an index $1\leq i\leq s+1$.
 We define the \emph{tops} $q_i^n$ and $q_i^m$ of $\Lambda$ by the property that $u_i^n\in I_{q_i^n}$ and $u_i^m\in I_{q_i^m}$. We also define the {\em tops } $q_i$ and $\tilde q_i$ to be such that $u_i\in I_{q_i}$ and $\tilde u_i\in I_{\tilde q_i}$. Recall that
$$
\{u_i^n,u_i^m\}=\{u_i,\tilde u_i\}.
$$
As a consequence, we have that $\{q_i^n,q_i^m\}=\{q_i,\tilde q_i\}$. Note that $q_i\leq \tilde q_i$, since $u_i<\tilde u_i$.

We also need to consider the integers $v_i$ that indicate the first levels $R_{v_i}(\Lambda)$ such that $z_i\in R_{v_i}(\Lambda)$. In other words, each $v_i$ is defined by the property that $\lambda_i\in I_{v_i}$, for $i=-1,0,1,\ldots,s$.

The following statements concern the properties of being circular intervals for the levels of $\Lambda$ and some derived properties of the conductor.

\begin{lema}[\cite{Can-C-SS-2023}, Lemma A.3]\label{lem:semi:circ-gamma}
Take $\mu\in \mathbb{Z}$, then $R_q(\mu +\Gamma)$ is a circular interval for all $q\in \mathbb Z$ (up to consider the emptyset as a circular interval).
\end{lema}

\begin{prop}[\cite{Can-C-SS-2023}, Proposition A.5]\label{prop:semi:paper1} Assume that $\Lambda$ is normalized (that is $\lambda_{-1}=0$) and that $R_q(\Lambda_{s-1})$ is a circular interval for any $q\geq v_s$. We have:
\begin{enumerate}
\item $<0_\varepsilon,z_s-1_\varepsilon>\subset R_q(\Lambda_{s-1})$, for $q\geq q^n_{s+1}-1$.
\item $<z_s,(n-1)_\varepsilon> \subset R_q(\Lambda)$, for $q\geq q^m_{s+1}-1$.
\end{enumerate}
In particular, we have that $R_q(\Lambda)=\mathbb S^1_n$, for any $q\geq \tilde q_{s+1}-1$.
\end{prop}

\begin{prop}[\cite{Can-C-SS-2023}, Proposition A.6]\label{prop:semi:paper2} Assume that
$\Lambda$ is normalized and  increasing. Then $R_q(\Lambda)$ is a circular interval for any $q\geq q_{s+1}$.
\end{prop}


\begin{corolario}\label{prop:cota-conductor}
 Assume that $\Lambda$ is normalized and  increasing. Then $\tilde{u}_{s+1}\geq c_{\Lambda}+n$, where $c_{\Lambda}$ is the conductor of $\Lambda$.
\end{corolario}
\begin{proof}
	 (See also \cite{Can-C-SS-2023}, Corollary A.7). First, let us show that $R_q(\Lambda_{s-1})$ is a circular interval for $q\geq v_s$.
	
	If $s=0$, we have that $\Lambda_{s-1}=\Lambda_{-1}=\Gamma$, we apply Lemma \ref{lem:semi:circ-gamma} by taking $\mu=0$. Assume now that $s\geq 1$. By Proposition \ref{prop:semi:paper2}, we know that $R_q(\Lambda_{s-1})$ is a circular interval for any $q\geq q_s$. Moreover, we have that $\lambda_s>u_s$ since $\Lambda$ is an increasing semimodule. This implies that $v_s\geq q_s$, hence we get that $R_q(\Lambda_{s-1})$ is a circular interval for any $q\geq v_s$, as desired.
	
	We end the proof as follows. By Proposition \ref{prop:semi:paper1}, we have that $R_q(\Lambda)=\mathbb S^1_n$, for any $q\geq \tilde q_{s+1}-1$.
 This implies that for any $k\geq n\tilde q_{s+1}-n$, we have that $k\in\Lambda$, and hence $k\geq c_\Lambda$. Finally, by definition of the tops, we have that
	$
	\tilde{u}_{s+1}\geq n\tilde q_{s+1}
	$
	and we are done.
\end{proof}

\begin{remark}\label{rem:semi:paper}
	Notice that Proposition \ref{prop:semi:paper1}, \ref{prop:semi:paper2} and \ref{prop:cota-conductor} are also true for increasing cuspidal semimodules such that $\lambda_{-1}$ is a multiple $nk$ of $n$. Indeed, in this case, we obtain the desired statements by applying the propositions to $\Lambda-nk$.
\end{remark}

\subsection{Distribution of the elements of the basis}

Along this section, we consider a cuspidal semimodule $\Lambda$ of length $s\geq 0$ with basis $\mathcal B$, that we read in the clock $\mathbb S^1_n$ as $B=\zeta(\mathcal B)$ as in the previous section.  We are going to describe a pattern for the distribution of the points $z_i$ in
$$
B=(z_{-1},z_0,z_1,\ldots,z_s)
$$
by computing the bounds $b_i^\ell(B)$ and $b_i^r(B)$ of $B$ in terms of the axes $u_{i+1}^n$ and $u_{i+1}^m$.

\begin{lema} Take $0\leq i\leq s$.
There are unique integer numbers $k_i^n$ and $k_i^m$ such that:
\begin{enumerate}
	\item $-1\leq k_i^n,k_i^m\leq i-1$.
	\item There is $b_{i+1}\geq 0$ such that  $u_{i+1}^n=\lambda_{i}+n\ell_{i+1}^n=\lambda_{k_i^n}+mb_{i+1}$.
	\item There is $a_{i+1}\geq 0$ such that $u_{i+1}^m=\lambda_{i}+m\ell_{i+1}^m=\lambda_{k_i^m}+na_{i+1}$.
\end{enumerate}
\end{lema}

\begin{proof}
	The existence of $k_i^n$ and $k_i^m$ comes from the definition of axes and limits. Let us show their uniqueness. Assume that there is another $k\neq k_i^n$ with $-1\leq k\leq i-1$ and a natural number $b$  such that
	$$
	u_{i+1}^n=\lambda_{i}+n\ell_{i+1}^n=\lambda_{k_i^n}+mb_{i+1}=\lambda_k+mb.
	$$
    Then either $\lambda_k\in (\lambda_{k_i^n}+\Gamma)$ or $\lambda_{k_i^n}
\in (\lambda_{k}+\Gamma)$ in contradiction with definition of basis. The uniqueness of $k_i^m$ is shown in the same way.
\end{proof}
The numbers $b_{i+1}$ and $a_{i+1}$ are the \emph{colimits} of $\Lambda$.

\begin{notation}\label{defn:bounds}
	 We denote $k_i$  and $\tilde k_i$ of $\tilde{u}_{i+1}$ by
	$$
	k_i=\left\{
	\begin{array}{lc}
		k_i^n,&\text{ if }u_{i+1}=u_{i+1}^n,\\
		k_i^m,&\text{ if }u_{i+1}=u_{i+1}^m.
	\end{array}
	\right.
	\quad
	\tilde k_i=\left\{
	\begin{array}{lc}
		k_i^n,&\text{ if }\tilde u_{i+1}=u_{i+1}^n,\\
		k_i^m,&\text{ if }\tilde u_{i+1}=u_{i+1}^m.
	\end{array}
	\right.
	$$
\end{notation}
\begin{remark}
\label{rk:separracionycolimites}
Note that $1\leq b_{i+1}<n$. To see this, it is enough to consider the case of a normalized $\Lambda$. Indeed, if $b_{i+1}\geq n$, we have that
$$
\lambda_i+n(\ell_{i+1}^n-1)=\lambda_{k_{i}^n}+m(b_{i+1})-n\geq (m-1)n\geq c_\Gamma\geq c_{\Lambda_{i-1}}.
$$
Thus $\lambda_i+n(\ell_{i+1}^n-1)\in\Lambda_{i-1}$ in contradiction with the minimality of $\ell_{i+1}^n$. Now, as a consequence, we have that the separation $S(z_{k_i^n},z_i)$ is given by $S(z_{k_i^n},z_i)=b_{i+1}$.  Recalling that $1\leq \ell_{i+1}^m<n$, see Remark \ref{rk:cotasparalimites}, we have that the separation $S(z_i, z_{k_i^m})$ is given by $S(z_i, z_{k_i^m})=\ell_{i+1}^m$.
\end{remark}
\begin{remark}
	Later on, we will show that $k_i^n=b_i^\ell(B)$ and $k_i^m=b_i^r(B)$.
\end{remark}

\begin{ejemplo}\label{ej:semi:t_1}
Take the semimodule $\Lambda=\Gamma\setminus \{0\}$. The basis of $\Lambda$ is $\mathcal B=(n,m)$. Note that $\lambda_{-1}=n$ and $\lambda_0=m$; thus, we have
$\Lambda_{-1}=n+\Gamma$ and $\Lambda_0=\Lambda=\Gamma\setminus \{0\}$.

The limit $\ell_1^n$ is the smallest positive integer such that
$$
m+n\ell_1^n=\lambda_0+n\ell_1^n\in\Lambda_{-1}=n+\Gamma.
$$
After solving the equation $m+n\ell_1^n=n+mb_1$, we obtain that  $\ell_1^n=1=b_1$. Moreover, we have $u_1^n=n+m=t_1^n$.

In the same way, in order to compute $\ell_1^m$, we solve $m+m\ell_1=n+na_1$, obtaining $\ell_1^m=n-1$ and $a_1=m-1$. Therefore, $u_1^m=nm=t_1^m$.

We conclude that $u_1=u_1^n=n+m$, $\tilde{u}_1=u_1^m=nm$, $t_1=t_1^n$ and $\tilde{t}_1=t_1^m$. As expected, we have that $k_0^n=k_0^m=-1$, that are the $0$-bounds of the list
$$
B=(0_\varepsilon, 1_\varepsilon)=(z_{-1},z_0),
$$
(note that $\zeta(m)=1_\varepsilon$).

Any cuspidal
semimodule  $\Lambda$ with basis $(n,m,\ldots)$ has the same first axes, first critical values, first limits, first colimits and $0$-bounds as the ones computed above, since their computation depends only on $\Lambda_{0}=\Gamma\setminus\{0\}$.
\end{ejemplo}

\begin{lema}
\label{lema: lemaees} Consider $0\leq i\leq s$ and take integer numbers $-1\leq k,k'\leq i-1$, with $k\neq k'$. Assume that we have the following equalities:
\begin{equation}\label{eq:lema:semi:1}
\lambda_i+ne=\lambda_k+mb; \qquad \lambda_i+ne'=\lambda_{k'}+mb',
\end{equation}
where $e,e'\in\mathbb Z$ and $0\leq b,b'<n$. Then we have that
$e<e'$ if and only if $b<b'$.
\end{lema}
\begin{proof} Equations \eqref{eq:lema:semi:1} lead us to:
\begin{eqnarray*}
\lambda_k&=&\lambda_{k'}+n(e-e')+m(b'-b),\\
\lambda_{k'}&=&\lambda_{k}+n(e'-e)+m(b-b').
\end{eqnarray*}
Note that $\lambda_k\notin\lambda_{k'}+\Gamma$ and $\lambda_{k'}\notin \lambda_k+\Gamma$, since $\lambda_k$ and $\lambda_{k'}$ are different elements of the basis of $\Lambda$. We conclude that $b<b'$ if and only if $e<e'$.
\end{proof}

\begin{prop} Consider \/ $0\leq i\leq s$ and take integer numbers $-1\leq k,k'\leq i-1$, with $k\neq k'$. We have
\begin{enumerate}
\item Assume that
$\lambda_i+ne=\lambda_k+mb$, $\lambda_i+ne'=\lambda_{k'}+mb'$,
where $e,e'\in\mathbb Z$ and $0\leq b,b'<n$. Then
$e<e'\Leftrightarrow
\lambda_i+ne<\lambda_i+ne'\Leftrightarrow S(z_k,z_i)<S(z_{k'},
z_{i})
$.
In particular, taking $k=k_i^n$, we have
$S(z_{_{k_i^n}}, z_i)<
S(z_{k'}, z_i)
$.
\item Assume that
$
\lambda_i+mf=\lambda_k+na$,  $\lambda_i+mf'=\lambda_{k'}+na'$
where $a,a'\in\mathbb Z$ and $0\leq f,f'<n$. Then
$ f<f'\Leftrightarrow
\lambda_i+mf<\lambda_i+mf'\Leftrightarrow
S(z_i,z_k)<S(z_{i},
z_{k'})
$.
In particular, taking $k=k_i^m$, we have
$S(z_i, z_{_{k_i^m}})<
S(z_i, z_{k'})
$.
\end{enumerate}
\end{prop}
\begin{proof}
Notice that  $S(z_i,z_k)=f$ and  $S(z_i,z_{k'})=f'$, this proves the second statement.
 For the first statement, we apply Lemma \ref{lema: lemaees}, by noting that
$
S(z_k,z_i)=b$ and $S(z_{k'},z_i)=b'$.
\end{proof}

\begin{corolario} We have that $k_i^n=b_i^\ell(B)$ and $k_i^m=b_i^r(B)$, for $0\leq i\leq s$.
\end{corolario}

\begin{remark}
	\label{rk:shiftedsemimodule} Take an integer number $\lambda\in\mathbb Z$. Then
$\mathcal B'=\lambda+\mathcal B$ is the basis of $\Lambda'=\lambda+\Lambda$ and $B'=\zeta(\mathcal B')=B+\lambda_\varepsilon$. Thus, the bounds of $B'$ are the same ones as the bounds of $B$. Anyway, the axes  for $\Lambda '$ are the ones of $\Lambda$ shifted by $\lambda$, this implies also that bounds, limits and colimits coincide for both semimodules.
\end{remark}

For the particular case when the semimodule $\Lambda$ is increasing, we can give a more accurate description of the bounds, as shown in  next proposition:

\begin{prop}\label{prop:semi:limites}
Assume that $\Lambda$ is increasing and  take $1\leq i\leq s$. We have
\begin{enumerate}
	\item If $u_{i}=u_{i}^n$, then $k_{i}^n=i-1$ and $k_{i}^m=k_{i-1}^m$.
	\item If $u_i=u_i^m$, then $k_{i}^n=k_{i-1}^n$ and $k_{i}^m=i-1$.
\end{enumerate}
\end{prop}

\begin{proof} In view of Remark \ref{rk:shiftedsemimodule}, it is enough to consider the normalized case $\lambda_{-1}=0$. Let us do the proof of (1); the proof of (2) is similar and we do not explicit it. Thus, we take the assumption that $u_{i}=u_i^n$.

First, let us suppose that $i=1$. By considering the bounds in the list $(z_{-1},z_0,z_1)$,  we deduce that $k_0^n=k_0^m=-1$ and either $z_1\in <z_{-1},z_0>$, or $z_1\in <z_{0},z_{-1}>$.
Let us show that we actually have that $z_1\in <z_{0},z_{-1}>$, this gives
$k_1^n=0$ and $k_1^m=-1$ as desired.

Since $\Lambda_{-1}=\Gamma$, we have that $R_q(\Lambda_{-1})$ is a circular interval for $q\geq 0$, due to Lemma \ref{lem:semi:circ-gamma}. Recall that $u_1^n=u_1\in I_{q_1}$. Noting that $z_{-1}=0_\varepsilon$ and applying  Proposition \ref{prop:semi:paper1} we have that
$$
<z_{-1},z_0-1_\varepsilon>\subset R_{q}(\Lambda_{-1}),\quad q\geq q_1^n-1=q_1-1.
$$
On the other hand, we have that $z_0\in  R_{q}(\Lambda_{0})$, for any $q\geq q_1$ since $\lambda_0< u_1$ and hence $v_0\leq q_1$. Thus, we have
$
<z_{-1},z_0>\subset R_{q_1}(\Lambda_{0})
$.
Note that $\lambda_1>u_1$, since
$\Lambda$ is increasing;  this implies that $z_1\notin R_{q_1}(\Lambda_0)$ and thus we necessarily have that $z_1\in <z_0,z_{-1}>$.

Now, assume that $i>1$. Our first step is to show that
$
z_i\in <z_{k_{i-1}^n}, z_{k_{i-1}^m}>
$.
By Proposition \ref{prop:semi:paper2}, we have that $R_{q}(\Lambda_{i-2})$ is a circular interval for $q\geq q_{i-1}$. Since $z_{k_{i-1}^n}$ and $z_{k_{i-1}^m}$ belong to $R_{q_{i-1}}(\Lambda_{i-2})$ we have that
 $$
 \text{Either }
 <z_{k_{i-1}^n},z_{k_{i-1}^m} >\subset R_{q_{i-1}}(\Lambda_{i-2}),\text{ or }
 <z_{k_{i-1}^m},z_{k_{i-1}^n} >\subset R_{q_{i-1}}(\Lambda_{i-2}).
 $$
 Noting that and $z_{i-1} \notin R_{q_{i-1}}(\Lambda_{i-2})$ and $z_{i-1}\in <z_{k_{i-1}^n},z_{k_{i-1}^m} > $, we conclude that
 $$
 <z_{k_{i-1}^m},z_{k_{i-1}^n} >\subset R_{q_{i-1}}(\Lambda_{i-2}).
 $$
 Noting also that $z_i\notin R_{q_{i-1}}(\Lambda_{i-2})$, we obtain that
$z_i\in <z_{k_{i-1}^n},z_{k_{i-1}^m} >$,
as desired.

Thus, we have $z_i,z_{i-1}\in <z_{k_{i-1}^n}, z_{k_{i-1}^m}>$ and hence there are two possibilities: either $z_i\in <z_{k_{i-1}^n},z_{i-1}>$, or $z_i\in <z_{i-1}, z_{k_{i-1}^m}>$.
Let us  show that $z_i\in <z_{i-1}, z_{k_{i-1}^m}>$ holds; in this way we are done. By Proposition \ref{prop:semi:paper1} we have that
$$
<0_\varepsilon, z_{i-1}-1_\varepsilon>\subset  R_{q_i-1}(\Lambda_{i-2})\subset R_{q_i-1}(\Lambda_{i-1}).
$$
Since $z_{i-1}\in R_{q_i}(\Lambda_{i-1})$, we have that $
<0_\varepsilon, z_{i-1}>\subset  R_{q_i}(\Lambda_{i-1})
$. Recalling that $k_{i-1}^n=b_{i-1}^\ell(B)$, we necessarily have that $z_{k_{i-1}^n}\in <0_\varepsilon,z_{i-1}>$, noting that $0_\varepsilon=z_{-1}$. Hence, we have
$$
<z_{k_{i-1}^n},z_{i-1}>\subset R_{q_i}(\Lambda_{i-1}).
$$
Since $z_i\notin R_{q_i}(\Lambda_{i-1})$, we obtain that $z_i\in <z_{i-1},z_{k_{i-1}^m}>$ as desired.
\end{proof}
\begin{remark}
\label{rk:semi:limites}
Note that Proposition \ref{prop:semi:limites} implies the following statements:
\begin{enumerate}
	\item If $k_{i}^n=i-1$, then $u_{i}=u_{i}^n$.
\item If $k_{i}^m=i-1$, then $u_{i}=u_{i}^m$.
\end{enumerate}
Indeed, we have that $i-1\in\{k_{i}^n, k_{i}^m\}$; if $k_{i}^n=i-1$, then necessarily $k_i^m\ne i-1$ (note that $i\geq 1$) and we are in the situation of the first statement of Proposition \ref{prop:semi:limites}. Similar argument when $k_{i}^m=i-1$.
\end{remark}

\subsection{Relations between parameters}

Let $\Lambda$ be an increasing cuspidal semimodule with basis
$
\mathcal B=(\lambda_{-1},\lambda_0,\lambda_1,\ldots,\lambda_s)
$ and let us denote
$$
B=\zeta(\mathcal B)=(z_{-1},z_0,z_1,\ldots,z_s).
$$
In this section we describe inductive features of axes, limits and co-limits of $\Lambda$.

\begin{lema}\label{lem:semi:u-tilde}
	Take $1\leq k<i\leq s+1$.  We have:
	\begin{enumerate}
		\item The axes and the critical values $u_i,t_i$ satisfy that $u_i>u_k$ and $t_i>t_k$.
		\item The axes and the critical values $\tilde u_i,\tilde t_i$ satisfy that $\tilde u_i<\tilde u_k$ and $\tilde t_i<\tilde t_k$.
	\end{enumerate}
\end{lema}

\begin{proof}
	It is enough to consider the case $k=i-1$.
	
	Let us prove Property (1). By definition of the axes, we have that $u_i>\lambda_{i-1}$. Since the semimodule is increasing, we have that $\lambda_{i-1}>u_{i-1}$. We get that $u_i>u_{i-1}$. Moreover, by Remark
	\ref{rk:relacionejesyvalorescriticos}, we see that
	$$
	t_i=t_{i-1}+(u_i-\lambda_{i-1})>t_{i-1}.
	$$
	This ends the proof of Property (1).
	
	Let us prove Property (2).  We do it for the case that
	$
	\tilde{u}_i=u_i^n=\lambda_{i-1}+n\ell_i^n
	$, the proof for the case $\tilde{u}_i=u_i^m$ runs in a similar way. By Proposition \ref{prop:semi:limites}, there are two cases: either $k_{i-1}^n=i-2$ or $k_{i-1}^m=i-2$. We shall see that $\tilde u_i<\tilde u_{i-1}$ and that $\tilde t_i<\tilde t_{i-1}$ simultaneously in each of the cases above.
	
{\em Case  $k_{i-1}^n=i-2$.} By Remark \ref{rk:semi:limites} we see that $u_{i-1}=u_{i-1}^n$ and $\tilde u_{i-1}=u_{i-1}^m$. Hence we can write:
	\begin{eqnarray}
		\tilde{u}_i =u_i^ n  = &\lambda_{i-1}+n\ell_{i}^n=\lambda_{i-2}+mb_i& \label{eq:lem:utilde1}\\
		\tilde{u}_{i-1}  =u^m_{i-1}  = &\lambda_{i-2}+m\ell_{i-1}^m=\lambda_{k}+na_{i-1}& \text{ with } k<i-2.\label{eq:lem:utilde2}
	\end{eqnarray}
	In order to see that $\tilde u_i<\tilde u_{i-1}$, we need to show that $b_i<\ell^m_{i-1}$. In order to do this, we are going to exclude the possibility $b_i\geq\ell^m_{i-1}$:
	\begin{itemize}
		\item If $\ell^m_{i-1}=b_i$,  we deduce that $\tilde{u}_{i-1}=\tilde{u}_i$ from equations \eqref{eq:lem:utilde1} and \eqref{eq:lem:utilde2}. Hence, we have
		$$
		\tilde u_i= \lambda_{i-1}+n\ell_{i}^n=\lambda_{k}+na_{i-1}=\tilde u_{i-1}, \quad k<i-2.
		$$
		Then $\rho({\lambda_{i-1}})=\rho({\lambda_k})\in \mathbb Z/n\mathbb Z$, contradicting  the fact that $\mathcal B$ is a basis.
		
		\item  If  $\ell^m_{i-1}<b_i$, by  equation \eqref{eq:lem:utilde1} and  by Proposition \ref{prop:cota-conductor} we have that:
		\begin{align*}
			\tilde u_i-n=\lambda_{i-1}+n(\ell_{i}^n-1)=&\lambda_{i-2}+mb_i -n \geq \lambda_{i-2}+m\ell^m_{i-1}+m-n\\
			=&\tilde{u}_{i-1}+m -n\geq c_{\Lambda_{i-2}}+m> c_{\Lambda_{i-2}}.
		\end{align*}
		We get that
		$
		\lambda_{i-1}+n(\ell_{i}^n-1)> c_{\Lambda_{i-2}}
		$ and thus
		$
		\lambda_{i-1}+n(\ell_i^n-1)\in \Lambda_{i-2}
		$,
		contradicting the minimality of $\ell_i^n$.
	\end{itemize}
	We conclude that $b_i<\ell^m_{i-1}$ and then we have that $\tilde{u}_i<\tilde{u}_{i-1}$, in the case $k_{i-1}^n=i-2$.
	
	Let us see now that $\tilde t_{i}<\tilde t_{i-1}$ in this case $k^n_{i-1}=i-2$. From equations \eqref{eq:lem:utilde1}, \eqref{eq:lem:utilde2}, using the fact that $\tilde{u}_i<\tilde{u}_{i-1}$ and the property of increasing semi\-module, we have that:
	$$
	\lambda_{i-2}+m\ell^m_{i-1}=\tilde{u}_{i-1}>\tilde{u}_i=\lambda_{i-1}+n\ell_i^n>u_{i-1}+n\ell_i^n.
	$$
	Consequently, $m\ell^m_{i-1}>u_{i-1}-\lambda_{i-2}+n\ell^n_i$ and
	\begin{eqnarray*}
		\tilde{t}_{i-1}=t_{i-2}+m\ell^m_{i-1}>t_{i-2}+u_{i-1}-\lambda_{i-2}+n\ell_i^n=t_{i-1}+n\ell_i^n=\tilde{t}_{i}.
	\end{eqnarray*}	
This ends the proof that $\tilde t_i<\tilde t_{i-1}$ in this case.

{\em Case $k_{i-1}^m=i-2$.}
Note that $\tilde u_{i-1}=u^n_{i-1}$ and $k^n_{i-1}=k^n_{i-2}$ in view of Remark \ref{rk:semi:limites} and  Proposition \ref{prop:semi:limites} . Thus, we can write
	\begin{eqnarray}
		\tilde{u}_i=u_i^n&=&\lambda_{i-1}+n\ell^{n}_i=\lambda_{k}+mb_i, \text{ with } k=k_{i-1}^n=k^n_{i-2}<i-2,
		\label{eq:cuatro}
		\\
		\tilde{u}_{i-1}=u_{i-1}^n&=&\lambda_{i-2}+n\ell_{i-1}^n=\lambda_{k}+mb_{i-1},
		\text{ with } k=k_{i-2}^n<i-2.
		\label{eq:cinco}
	\end{eqnarray}
 Let us  proceed in a similar way as before to show that $b_{i-1}>b_i$:
	\begin{itemize}
		\item Assume that $b_{i-1}=b_i$. Then $\rho(\lambda_{i-1})=\rho(\lambda_{i-2})$, absurd.
		\item Assume that $b_{i-1}<b_{i}$. Then, we have that
			\begin{eqnarray*}
			\tilde u_i-n=\lambda_{i-1}+n(\ell^n_i-1)&=&\lambda_k+mb_i-n\\
			&>& \lambda_k+(b_i-1)m=\tilde u_{i-1}+(b_i-b_{i-1}-1)m \\
			&\geq& \tilde u_{i-1}\geq n+c_{\Lambda_{i-2}}.
		\end{eqnarray*}
		Then $\lambda_{i-1}+n(\ell^n_i-1)\in\Lambda_{i-2}$, in contradiction with the minimality of $\ell^n_{i}$.
	\end{itemize}
	We conclude that $b_{i-1}>b_i$ and thus $\tilde u_{i-1}>\tilde u_i$.

Let us see now that $\tilde t_{i}<\tilde t_{i-1}$ in this case $k^m_{i-1}=i-2$.  From equations \eqref{eq:cuatro}, \eqref{eq:cinco}, using the fact that $\tilde{u}_i<\tilde{u}_{i-1}$ and the property of increasing semi\-module, we have that:
$$
\lambda_{i-2}+n\ell^n_{i-1}=\tilde{u}_{i-1}>
\tilde{u}_i=\lambda_{i-1}+n\ell_i^n>u_{i-1}+n\ell_i^n.
$$
Consequently, $n\ell^n_{i-1}>u_{i-1}-\lambda_{i-2}+n\ell^n_i$ and
	\begin{eqnarray*}
		\tilde{t}_{i-1}=t_{i-2}+n\ell^n_{i-1}>t_{i-2}+u_{i-1}-\lambda_{i-2}+n\ell_i^n=t_{i-1}+n\ell_i^n=\tilde{t}_{i}.
	\end{eqnarray*}	
	This ends the proof.
\end{proof}
\begin{corolario}
	\label{cor:cotastes}
	Let $\Lambda$ be a cuspidal increasing semimodule with basis
	$$
	\mathcal B=(\lambda_{-1},\lambda_0,\lambda_1,\ldots,\lambda_s)
	$$
	such that $\lambda_{-1}=n$ and $\lambda_0=m$. We have that $\tilde t_1=nm$ and the following holds
	$$
	t_{i+1}^n<\tilde t_1=nm \quad \text{ and }\quad t_{i+1}^m<\tilde t_1=nm,
	$$
	for any $1\leq i\leq s$.
\end{corolario}
\begin{proof}
	It is enough to recall that  $\tilde t_1=nm$ in view of Example \ref{ej:semi:t_1}.
\end{proof}

We end this section with a proposition that connects the limits and the colimits.

\begin{prop}\label{prop:semi:li+ai}
Consider $1\leq i\leq s$. We have
\begin{enumerate}
	\item If $k^n_i=i-1$, then
	$
	\ell^n_{i+1}+a_{i+1}=a_{i}$ and  $\ell^m_{i+1}+b_{i+1}=\ell^m_i
	$.
	\item  If $k^m_i=i-1$, then
	$
	\ell^n_{i+1}+a_{i+1}=\ell^n_i$ and  $\ell^m_{i+1}+b_{i+1}=b_i
	$.
\end{enumerate}
\end{prop}

\begin{proof}
 Notice that shifting the semimodule any integer number does not change the value of the limits and the colimits. Therefore, we can assume  without loss of generality that $\Lambda$ is normalized and thus $\lambda_{-1}=0$.

 Let us prove statement (1). By hypothesis, we have that $k_{i}^n=i-1$. In view of Remark \ref{rk:semi:limites} and  Proposition \ref{prop:semi:limites}, we also have that $k_{i}^m=k^{m}_{i-1}$. Let us write:
 \begin{eqnarray}
 \label{eq:1:lem:li+a1}
u^n_{i+1}&=&\lambda_i+n\ell^n_{i+1}=\lambda_{i-1}+mb_{i+1},\\
\label{eq:2:lem:li+a1}
u^m_{i+1}&=&\lambda_i+m\ell^m_{i+1}=\lambda_{k^m_{i}}+na_{i+1},\\
\label{eq:3:lem:li+a1}
u_i^m&=&\lambda_{i-1}+m\ell_i^m=\lambda_{k_{i-1}^m}+na_i= \lambda_{k_{i}^m}+na_i.
 \end{eqnarray}
From equations \eqref{eq:1:lem:li+a1} and \eqref{eq:2:lem:li+a1} we obtain that
\begin{equation}
\label{eq:4:lem:li+a1}
n\ell^n_{i+1}+na_{i+1}+\lambda_{k_{i}^m}=mb_{i+1}+m\ell^m_{i+1}+\lambda_{i-1}.
\end{equation}
By equation \eqref{eq:3:lem:li+a1} we can substitute $\lambda_{k_{i}^m}=\lambda_{i-1}+m\ell_i^{m}-na_i$ in equation \eqref{eq:4:lem:li+a1} to obtain
\begin{equation}
\label{eq:colimites5}
n(\ell^n_{i+1}+a_{i+1}-a_i)=m(\ell^m_{i+1}+b_{i+1}-\ell_i^m).
\end{equation}
Since $n$ and $m$ have no common factor, we have that  $n$ divides $\ell^m_{i+1}+b_{i+1}-\ell_i^m$.

Let us see that $\ell^m_{i+1}+b_{i+1}-\ell_i^m=0$ and hence $\ell^m_{i+1}+b_{i+1}=\ell_i^m$ as desired. If $\ell^m_{i+1}+b_{i+1}-\ell_i^m\ne 0$ we are in one of the following three cases
\begin{equation*}
a)\; |\ell^m_{i+1}+b_{i+1}-\ell_i^m|\geq 2n,\quad
b)\; \ell^m_{i+1}+b_{i+1}-\ell_i^m=-n, \quad
c)\; \ell^m_{i+1}+b_{i+1}-\ell_i^m=n.
\end{equation*}
Let us see that each of these cases leads to a contradiction.

Assume first that we are in case a).  Noting that $\ell^m_{i+1},b_{i+1},\ell_i^m\geq 1$, there is at least one of them that is strictly bigger than $n$. Let us consider the three possibilities:
\begin{itemize}
\item
If $\ell^m_{i+1}>n$, we have that $m\ell^m_{i+1}>nm$ and then
$
\lambda_i+m\ell^m_{i+1}>nm
$.
This implies that $\lambda_i+m(\ell^m_{i+1}-1)>(n-1)m\geq c_\Gamma\geq c_{\Lambda_{i-1}}$. Then, we have that $\lambda_i+m(\ell^m_{i+1}-1)\in \Lambda_{i-1}$, contradicting the minimality of $\ell^m_{i+1}$.
\item If $\ell_i^m>n$, we do the same argument as before.
\item If $b_{i+1}>n$, we have that $\lambda_{i}+n\ell^n_{i+1}=\lambda_{i-1}+mb_{i+1}>nm$ and then
$$
\lambda_i+n(\ell^n_i-1)>(m-1)n\geq c_\Gamma\geq c_{\Lambda_{i-1}}.
$$
Then $\lambda_i+n(\ell^n_i-1)\in\Lambda_{i-1}$ and this contradicts the minimality of $\ell_i^n$.
\end{itemize}
Assume that we are in case b), that is $\ell^m_{i+1}+b_{i+1}-\ell_i^m=-n$. this implies that $\ell_i^m>n$ and we do the same argument as before to obtain a contradiction.

Assume that we are in case c), that is $\ell^m_{i+1}+b_{i+1}-\ell_i^m=n$. We have that
$\ell^m_{i+1}+b_{i+1}>n $. By Remark \ref{rk:separracionycolimites} we see that the separation $S(z_{i-1},z_i)$ is given by $S(z_{i-1},z_i)=b_{i+1}$ (recall that $k_{i}^n=i-1$) and that the separation $S(z_i, z_{k_i^m})$ is given by $S(z_i, z_{k_i^m})=\ell_{i+1}^m$. Noting that
$
z_i\in <z_{i-1},z_{k_i^m}>
$
and $z_{i-1}\ne z_{k_i^m}$, we conclude that
$$
n>S(z_{i-1},z_i)+S(z_i,z_{k_i^m})=b_{i+1}+\ell_{i+1}^m.
$$
This contradicts $b_{i+1}+\ell_{i+1}^m>n$. The proof that $\ell^m_{i+1}+b_{i+1}=\ell_i^m$ is ended. Moreover, since $\ell^m_{i+1}+b_{i+1}-\ell_i^m=0$, by equation \eqref{eq:colimites5}, we conclude that $\ell^n_{i+1}+a_{i+1}=a_i$, as desired.

The proof of statement (2) runs in a similar way to the above arguments.
\end{proof}

Next corollary will be useful in our computation of Saito bases:
\begin{corolario}
\label{cor:cotaseles}
Consider $2\leq j+1<q\leq s+1$. Then
$$
\ell_{j+1}^m-(\ell_{j+2}^m+\ell_{j+3}^m+\cdots+\ell_q^m)=b_q>0,
$$
under the assumption that $\tilde t_{j+1}=t_{j+1}^m$ and $\tilde t_\ell=t_\ell^n$, for $j+2\leq \ell\leq q-1$.
In a symmetric way, we have that
$$
\ell_{j+1}^n-(\ell_{j+2}^n+\ell_{j+3}^n+\cdots+\ell_q^n)=a_q>0,
$$
under the assumption that $\tilde t_{j+1}=t_{j+1}^n$ and $\tilde t_\ell=t_\ell^m$, for $j+2\leq \ell\leq q-1$.
\end{corolario}
\begin{proof} We prove the first assertion, the second one is similar. Let us consider the difference
$$
\ell^m_{j+1}-\ell^m_{j+2}.
$$
Since $\tilde t_{j+1}=t^m_{j+1}$, by Proposition \ref{prop:semi:limites}, we have that $k_{j+1}^n=j$. By Proposition \ref{prop:semi:li+ai}, we conclude that
$$
\ell^m_{j+1}-\ell^m_{j+2}=b_{j+2}.
$$
Now, let us study the difference $b_{j+2}-\ell_{j+3}^m$. Since $t_{j+2}=t_{j+2}^m$, we have that $k_{j+2}^m=j+1$. By Proposition  \ref{prop:semi:li+ai} we conclude that
$$
b_{j+2}-\ell_{j+3}^m=b_{j+3}.
$$
Following in this way, we conclude that
$$
\ell_{j+1}^m-(\ell_{j+2}^m+\ell_{j+3}^m+\cdots+\ell_q^m)=b_q>0,
$$
as desired.
\end{proof}

\section{Cuspidal Standard Bases}\label{sec:bases}
\subsection{Semigroup and semimodule of an analytic branch}
Let us consider the local ring $\mathcal O_{\mathbb C^2,\mathbf 0}$ of the analytic space $\mathbb C^2$ at the origin. Denote by $x,y$ the coordinates of $\mathbb C^2$, that we consider as elements $x,y\in \mathcal O_{\mathbb C^2,\mathbf 0}$. We recall that there is an identification
$
 \mathcal O_{\mathbb C^2,\mathbf 0}=\mathbb C\{x,y\}
$
between the local ring  $\mathcal O_{\mathbb C^2,\mathbf 0}$ and the ring of convergent power series in $x,y$ with complex coefficients. By definition, an {\em analytic plane branch $C$} at the origin of $\mathbb C^2$ is a principal prime ideal $C\subset \mathcal O_{\mathbb C^2,\mathbf 0}$. Any generator $f$ of $C$ is called an {\em equation of $C$}. It is known that there is a morphism
$$
\varphi: \mathcal O_{\mathbb C^2,\mathbf 0}\rightarrow \mathbb C\{t\}
$$
such that $C=\ker \varphi$. Such morphisms are given in terms of convergent series by
$$
\varphi(g(x,y))=g(a(t),b(t)),\quad a(t),b(t)\in \mathbb C\{t\}.
$$
We call them {\em parametrizations  of $C$}. The subset
$$
\Gamma_\varphi=\{\operatorname{ord}_t(\varphi(g));\; g\in \mathcal O_{\mathbb C^2,\mathbf 0}\}\subset \mathbb Z_{\geq 0}
$$
is a semigroup of $\mathbb Z_{\geq 0}$.  A parametrization $\varphi$ is called {\em primitive} if and only if the following equivalent properties hold:
\begin{enumerate}
	\item The semigroup $\Gamma_\varphi$ has a conductor, that is, there is $c_{\Gamma_\varphi}\in \mathbb Z$ minimal with the property that
	$n\in \Gamma_{\varphi}$ for any $n\geq c_{\Gamma_\varphi}$.
	\item  There is no series $\psi(t)\ne 0$ with $\operatorname{ord}_t(\psi(t))\geq 2$ such that
	$$
	\varphi(g(x,y))=g( a(t),b(t)),
	$$
	 where $ a(t)= \tilde a(\psi(t))$
and  $ b(t)= \tilde b(\psi(t))$. In particular, we have another parametrization $\tilde\varphi$ given by $\tilde\varphi(g(x,y))=g(\tilde a(t),\tilde b(t))$.
\end{enumerate}
There are always primitive parametrizations. If $\varphi$ and $\varphi'$ are primitive parametrizations of the plane branch $C$, we have that
$$
\operatorname{ord}_t(\varphi(g))= \operatorname{ord}_t(\varphi'(g)),\quad
g\in\mathcal O_{\mathbb C^2,\mathbf 0}.
$$
We denote $\nu_C(g)=\operatorname{ord}_t(\varphi(g))$. We also conclude that
$
\Gamma_{\varphi}=\Gamma_{\varphi'}
$ and we call
 this semigroup the {\em semigroup $\Gamma_C$ of the plane branch $C$}.

We say that $C$ is a {\em plane cusp} if the semigroup $\Gamma_C$ is generated by two integer numbers $2\leq n <m$ without common factor. In this paper we mainly deal with plane cusps.

Let $C$ be a plane cusp. We know that, after an adequate coordinate change, we can choose coordinates and a primitive parametrization
$$
x=t^n,\; y=b(t),\quad \operatorname{ord}_tb(t)=m,
$$
where $2\leq n<m$  and $n,m$ are without common factors. These coordinates are called {\em adapted coordinates} and the above parametrization is a {\em Puiseux parametrization associated to the adapted coordinates}.

\begin{remark} Zariski's Equisingularity Theory concerns all the plane branches, we are fixing our attention in the cusps, that are the branches with a single Puiseux pair $(n,m)$.
\end{remark}
Let us denote by $\Omega^1_{\mathbb C^2,\mathbf 0}$ the $\mathcal O_{\mathbb C^2,\mathbf 0}$-module of germs of holomorphic differential $1$-forms at the origin of $\mathbb C^2$. We know that $\Omega^1_{\mathbb C^2,\mathbf 0}$ is a free rank-two $\mathcal O_{\mathbb C^2,\mathbf 0}$-module generated by $dx$ and $dy$. In a similar way, we denote by $\Omega^1_{\mathbb C,0}$ the $\mathcal O_{\mathbb C,0}$-module of germs of holomorphic differential $1$-forms at the origin of $\mathbb C$. We also denote by $\Omega_{\mathbb C^2,\mathbf 0}^2$ the $\mathcal O_{\mathbb C^2,\mathbf{0}}$-module of germs of holomorphic differential 1-forms at the origin of $\mathbb C^2$, in this case it is a free module of rank 1 generated by $dx\wedge dy$.

 An element $\alpha\in \Omega^1_{\mathbb C,0}$ is written as $\alpha=\psi(t)dt$. The {\em order} of $\alpha$ is by definition the order of $\psi(t)$. More precisely, we write
 $$
 \operatorname{ord}_t(\alpha)=\operatorname{ord}_t(\psi(t)).
 $$

Given a primitive parametrization $\varphi$ of an analytic plane branch $C$, we have a ``pull-back'' application
 $$
 \varphi^*:\Omega^1_{\mathbb C^2,\mathbf 0}\rightarrow \Omega^1_{\mathbb C,\mathbf 0},\quad \omega\mapsto \varphi^*\omega,
 $$
defined by the properties that $\varphi^*(\omega+\omega')=\varphi^*(\omega)+\varphi^*(\omega')$,
$\varphi^*(h\omega)=\varphi(h)\varphi^*(\omega)$ and
$$
\varphi^*(dx)= \frac{\partial \varphi(x)}{\partial t}dt, \quad
\varphi^*(dy)= \frac{\partial \varphi(y)}{\partial t}dt.
$$
The set
$$
\Lambda_C=\{\operatorname{ord}_t(\varphi^*\omega)+1;\; \omega\in \Omega^1_{\mathbb C^2,\mathbf 0})\}\setminus\{\infty\}\subset \mathbb Z_{\geq 0}
$$
is the so-called {\em semimodule of differential values for $C$}. This set is independent of the chosen primitive parametrization and also of the analytic class of $C$, where two analytic plane branches are analytically equivalent when they correspond one to another via an automorphism of the local ring $\mathcal O_{\mathbb C^2,\mathbf 0}$.

We take the notation $\nu_C(\omega)=\operatorname{ord}_t(\varphi^*(\omega))+1$ and we call it the {\em differential order of $\omega$ with respect to $C$}.
\begin{remark} The set $\Lambda_C$ is a semimodule over the semigroup $\Gamma_C$. Namely, we have that
$$
\nu_C(h)+ \nu_C(\omega)= \nu_C(h\omega).
$$
Moreover, if $\nu_C(h)>0$, we have that $\nu_C(dh)= \nu_C(h)$, where the notation $dh$ stands for the differential
$$
dh=(\partial h/\partial x)dx+(\partial h/\partial y)dy.
$$
This means that $\Gamma_{C}\setminus \{0\}\subset \Lambda_C$ and $c_{\Lambda_C}\leq c_{\Gamma_C}$.
\end{remark}

Consequently, there is a basis $\mathcal B=(\lambda_{-1},\lambda_0,\lambda_1,\ldots,\lambda_s)$ for $\Lambda_C$. A set of differential forms
$$
\mathcal S=(\omega_{-1},\omega_0,\omega_1,\ldots,\omega_s)
$$
such that $\nu_C(\omega_i)=\lambda_i$, for $i=-1,0,1,\ldots,s$, will be called a {\em minimal standard basis of differential $1$-forms} for the analytic plane branch $C$.

\subsection{Minimal Standard Bases}
From now on, we consider a plane cusp $C$ with Puiseux pair $(n,m)$, we denote its semigroup by $\Gamma=\Gamma_C$ and its semimodule of differential values by $\Lambda=\Lambda_C$. We know that $\Gamma$ is generated by $n,m$, where $2\leq n<m$,  without common factor. We also fix the notation
$$
\mathcal B=(\lambda_{-1},\lambda_0,\ldots,\lambda_s)
$$
for the basis of $\Lambda$. We recall that $\Gamma\setminus\{0\}\subset \Lambda$ and we have that $\lambda_{-1}=n$, $\lambda_0=m$.

We also fix adapted coordinates $x,y$ and a Puiseux parametrization
$$
x=t^n,\quad y=b(t),\quad \operatorname{ord}_t{b(t)}=m.
$$
Note that $\nu_C(x)=\nu_C(dx)=n$ and $\nu_C(y)=\nu_C(dy)=m$.
\subsubsection{Divisorial Order}
Let $\pi_C:M\rightarrow (\mathbb C^2,\mathbf 0)$ be the minimal reduction of the singularities of the cusp  $C$. That is, the morphism $\pi_C$ is the minimal finite composition of blowing-ups centered at points in the successive strict transforms of $C$ in such a way that the strict transform $C'$ of $C$ by $\pi_C$ in $M$ has normal crossings with the total exceptional divisor $E=\pi^{-1}(\mathbf 0)$. Let us denote by $D_C$ the only irreducible component of $E$ such that $C'\cap D_C\ne \emptyset$.
\begin{lema} Let $C,\tilde C$ be two cusps with the same Puiseux pair $(n,m)$. The following statements are equivalent:
\begin{enumerate}
\item $\pi_C=\pi_{\tilde C}$.
\item $D_C=D_{\tilde C}$.
\item Given a local coordinate system $(x,y)$ of $(\mathbb C^2,\mathbf 0)$, we have that $(x,y)$ is adapted to $C$ if and only if it is adapted to $\tilde C$.
\end{enumerate}
\end{lema}
\begin{proof} We leave this proof to the reader.
\end{proof}
Denote $D=D_C$ and $\pi=\pi_C$. Let us define the {\em divisorial order $\nu_D$} as follows. Given an adapted coordinate system $(x,y)$ and a point $Q\in D$, there are local coordinates $(x',y^*)$ such at $Q$ such that $D=(x'=0)$ locally at $Q$ and there is a complex number $\lambda$ with $y^*=y'-\lambda$, where
$$
\begin{array}{ccc}
x&=&(x')^n (y')^b \\
y&=&(x')^m (y')^d
\end{array},
\quad nd-mb=1.
$$
Consider $h\in \mathcal O_{\mathbb C^2,\mathbf 0}$. Let us write the germ $(h\circ\pi)_{Q}$ of $h$ at $Q$ as
$$
(h\circ\pi)_{Q}=(x')^{\beta}h',
$$
where $x'$ does not divide $h'$. Then, we define $\nu_D(h)=\beta$.
In a similar way, given a $1$-form $\omega\in \Omega^1_{\mathbb C^2,\mathbf 0}$, let us write the germ $(\pi^*\omega)_{Q}$ as
$$
(\pi^*\omega)_{Q}=(x')^{\beta}\left\{a'\frac{dx'}{x'}+b'dy^*\right\},
$$
where $x'$ does not divide the pair $(a',b')$. Then, we define $\nu_D(\omega)=\beta$.

The above definition of divisorial order does not depend on the chosen adapted coordinates nor on the particular point $Q\in D$. We also call $\nu_D$ the {\em monomial order} in view of the next statement:
\begin{lema}
\label{lema:monomialorder}
Let $(x,y)$ be an adapted coordinate system for $C$ and consider $h\in\mathcal O_{\mathbb C^2,\mathbf 0}=\mathbb C\{x,y\}$ that we write as
$$
h=\sum_{i,j\geq 0}h_{ij}x^iy^j.
$$
Then we have that
$
\nu_D(h)=\min\{ni+mj;\; h_{ij}\ne 0\}= \min\{\nu_C(x^ix^j);\; h_{ij}\ne 0\}.
$
In a similar way, take $\omega\in \Omega^1_{\mathbb C^2,\mathbf 0}$ that we write as
$$
\omega=\sum_{i,j\geq 0}x^iy^j\left\{\lambda_{ij}\frac{dx}{x}+\mu_{ij}\frac{dy}{x}\right\}
$$
where $(\lambda_{ij},\mu_{ij})\in \mathbb C$. Then, we have that
$
\nu_D(h)= \min\{ni+mj;\; (\lambda_{ij},\mu_{ij})\ne (0,0)\}.
$
\end{lema}
\begin{proof}Left to the reader, see \cite{Can-C-SS-2023}.
\end{proof}
\begin{remark} Take $\omega=nydx-mxdy=xy\left(ndx/x-mdy/y\right)$, note that
$$
\nu_D(\omega)= n+m<\nu_C(\omega).
$$
Let us also note that if $\omega=(adx+bdy)$, then $\nu_D(\omega)=\min\{\nu_D(xa),\nu_D(yb)\}$.
\end{remark}

\begin{remark}  For the case of functions as well as for differential $1$-forms, we have that $\nu_C(h)\geq \nu_D(h)$ and $\nu_C(\omega)\geq \nu_D(\omega)$.
\end{remark}
\subsubsection{Initial parts}
\label{Initialparts}
Here, we fix a coordinate system $(x,y)$ adapted to the cusp $C$. Given a function $h=\sum_{i,j}h_{ij}x^iy^j\in \mathcal O_{\mathbb C^2,\mathbf 0}=
\mathbb C \{x,y \}$ with $\nu_D(h)=p$, we define the initial part $\operatorname{In}_{n,m}^{x,y}(h)$ by
$$
\operatorname{In}_{n,m}^{x,y}(h)= \sum_{ni+mj=p}h_{ij}x^iy^j.
$$
If there is no confusion, we write $ \operatorname{In}(h)=
\operatorname{In}_{n,m}^{x,y}(h)
$. In a similar way, given a $1$-differential form $\omega\in \Omega^1_{\mathbb C^2,\mathbf 0}$, with $p=\nu_D(\omega)$, that we write as
$$
\omega=\sum_{ni+mj\geq p}x^iy^j\left(\lambda_{ij}\frac{dx}{x}+\mu_{ij}\frac{dy}{y}\right),
$$
we define the initial part $\operatorname{In}(\omega)$ by
$$
\operatorname{In}(\omega)=\sum_{ni+mj= p}x^iy^j\left(\lambda_{ij}\frac{dx}{x}+\mu_{ij}\frac{dy}{y}\right).
$$
\begin{remark}
 \label{rk:partieinicialvalorbajo}
 Assume that $\nu_D(\omega)=p< nm$. Then there are unique $\alpha,\beta\in \mathbb Z_{\geq 0}$, with $(\alpha,\beta)\ne (0,0)$,  such that $n\alpha+m\beta=p$ and hence the initial part is a single differential monomial
$$
\operatorname{In}(\omega)= x^\alpha y^\beta\left(\lambda\frac{dx}{x}+\mu\frac{dy}{y}\right).
$$
We can do a similar argument for the case of a function $h\in\mathcal O_{\mathbb C^2,\mathbf 0}$.
\end{remark}

Let $\omega=adx+bdy$ be a differential $1$-form in $\Omega^1_{\mathbb C^2,\mathbf 0}$ such that $\nu_D(\omega)<nm$ and such that $\nu_C(\omega)>\nu_D(\omega)$. Then, the initial part $\operatorname{In}(\omega)$ is given by
$$
\operatorname{In}(\omega)=\mu x^\alpha y^\beta\left(n\frac{dx}{x}-m\frac{dy}{y}\right)
$$
where $\mu\ne 0$, and $\alpha,\beta\geq 1$ are such that $\nu_D(\omega)=n\alpha+m\beta<nm$. Any differential $1$-form as above will be called a {\em resonant differential $1$-form}.

Given two differential $1$-forms $\omega$ and $\eta$, we say that $\eta$ is {\em reachable} from $\omega$ if there is a monomial $\mu x^\alpha y^\beta$ such that
$$
\operatorname{In}(\eta)= \mu x^\alpha y^\beta\operatorname{In}(\omega).
$$
\subsubsection{Semimodule Versus Minimal Standard Bases}
Let us consider a minimal standard basis $\mathcal S=(\omega_{-1},\omega_0,\omega_1,\ldots,\omega_s)$ of the cusp $C$. Recall that $\nu_C(\omega_{i})=\lambda_i$, for $i=-1,0,1,\ldots,s$.

\begin{lema}
 \label{lema:firstvalues}
 We have that $\lambda_{-1}=n$ and $\lambda_0=m$. More precisely, the initial parts of $\omega_{-1}$ and $\omega_0$ are respectively given by
$\operatorname{In(\omega_{-1})}=\lambda dx$, with $\lambda\ne 0$, and $\operatorname{In}(\omega_{0})=\mu dy$, with $\mu\ne 0$.
\end{lema}
\begin{proof} We have that $\nu_C(adx)=\nu_C(a)+\nu_C(dx)$,  $\nu_C(bdy)=\nu_C(b)+\nu_C(dy)$ and
$$
\nu_C(adx+bdy)\geq \min\{\nu_C(adx),\nu_C(bdy)\}.
$$
Since $\nu_C(dx)=n$ and $\nu_C(dy)=m$, with $n<m$, we conclude that $n=\min\Lambda=\lambda_{-1}$. We also have that $\nu_C(\omega_{-1})=n$. Let us write
$$
\omega_{-1}=\lambda dx+\eta;\quad \eta=x\eta_1+y\eta_2+hdy.
$$
We have that $\nu_D(\eta)>n$ and hence $\nu_C(\eta)\geq\nu_D(\eta)>n$. The only possibility to have that $\nu_C(\omega_{-1})=n$ is that $\lambda\ne 0$ and, in this case, we see that $\operatorname{In}(\omega_{-1})=\lambda dx$.

Let us show that $\lambda_0=m$ and that $\operatorname{In}(\omega_0)=\mu dy$. Let $k\geq 1$ be the integer number defined by the property that $kn<m<(k+1)n$. Take a differential $1$-form $\omega$ that we write as
$$
\omega=(c_1+c_2x+c_3x^2+\cdots+c_{k}x^{k-1})dx+\eta,\quad \eta=x^k\eta_1+y\eta_2+hdy
$$
We have that $m\leq \nu_D(\eta)\leq\nu_C(\eta)$. If $(c_1,c_2,\ldots,c_k)\ne (0,0,\ldots, 0)$, taking the first $j$ such that $c_j\ne 0$ we conclude that
$$
\nu_C(\omega)=jn\in \lambda_{-1}+\Gamma.
$$
Thus, the next differential value $\lambda_0$ is given by a differential $1$-form $\eta$ written as $\eta=x^k\eta_1+y\eta_2+hdy$. Let us decompose
$$
\eta=\mu dy+\tilde\eta;\quad \tilde\eta=x^k\eta_1+y\eta_2+(xh_1+yh_2)dy.
$$
We have that $m<\nu_D(\tilde\eta)<\nu_C(\tilde\eta)$. Thus,  If $\mu\ne 0$, we get
$
m=\nu_D(\eta)=\nu_C(\eta).
$
The desired result follows.
\end{proof}
\begin{remark} Recall the decomposition chain
$\{
\Lambda_{i}
\}_{i=-1}^s$ of $\Lambda$. In view of Lemma \ref{lema:firstvalues}, we see  that $\Lambda_{-1}=n+\Gamma,$ $\Lambda_0=(n+\Gamma)\cup (m+\Gamma)=\Gamma\setminus\{0\}$.
\end{remark}

Assume that $s\geq 1$ and let us describe the initial part $\operatorname{In}(\omega_1)$ of the element $\omega_1$ in the minimal standard basis $\mathcal S$. Take $\omega$ such that $\nu_C(\omega)\notin \Lambda_0$; the value $\lambda_1$ is the minimum of the differential values $\nu_C(\omega)$ for such $1$-forms $\omega$. The first remark is that $\nu_C(\omega)<nm$, since $c_\Gamma=(n-1)(m-1)$ and hence $\nu_D(\omega)<nm$; the second remark is that $\omega$ is resonant. Namely, in view of Remark  \ref{rk:partieinicialvalorbajo}, we can write
$$
\operatorname{In}(\omega)=x^\alpha y^\beta\left\{\lambda\frac{dx}{x}+\mu\frac{dy}{y}\right\},
\quad
n\alpha+m\beta=\nu_D(\omega)<nm.
$$
If $\lambda n-\mu m\ne 0$, we have that $\nu_C(\omega)=\nu_D(\omega)=n\alpha+m\beta\in\Lambda_0$. Then, the $1$-form $\omega$ is necessarily resonant, that is
$$
\operatorname{In}(\omega)=\mu x^\alpha y^\beta\left\{n\frac{dx}{x}-m\frac{dy}{y}\right\}.
$$
Assume now that $\nu_C(\omega)=\lambda_1$ (this property is satisfied by $\omega_1$), let us show that $\alpha=\beta=1$ (note that we necessarily have that $\alpha\geq 1$ and $\beta\geq 1$). Let us reason by contradiction, assuming that  $n\alpha+m\beta>n+m$. We start with the differential $1$-form
$$
\eta=nydx-mxdy.
$$
We know that $\nu_D(\eta)=n+m<\nu_C(\eta)$. If $\nu_C(\eta)=a_1n+b_1m\in \Lambda_0$ and $\nu_C(\eta)<nm$, there is $\mu_1\ne 0$, such that
$$
\nu_C(\eta+\mu_1d(x^{a_1}y^{b_1}))>\nu_C(\eta).
$$
Put $\eta^1=\eta+\mu_1d(x^{a_1}y^{b_1})$ and re-start the procedure with $\eta^1$. In this way, we obtain a differential $1$-form $\tilde \eta$ with the following properties:
\begin{enumerate}
\item Either $\nu_C(\tilde \eta)\notin \Lambda_0$ or $\nu_C(\tilde \eta)>nm$.
\item $\operatorname{In}(\tilde\eta)=\eta=nydx-mxdy$ and hence $\nu_D(\tilde\eta)=n+m$.
\end{enumerate}
Now, we compare $\tilde\eta$ with $\omega$ as follows. We know that $\lambda_1\leq \nu_C(\tilde\eta)$. Moreover, there is a constant $\mu\ne 0$ such that
$$
\nu_D(\omega-\mu x^{\alpha-1}y^{\beta-1}\tilde\eta)>\nu_D(\omega)=\alpha n+\beta m.
$$
Put $\omega^1=\omega-\mu x^{\alpha-1}y^{\beta-1}\tilde\eta$. We have that $\nu_C(\omega^1)=\lambda_1$, since $\nu_C(x^{\alpha-1}y^{\beta-1}\tilde\eta)>\lambda_1$. We re-start with $\omega^1$. Repeating the procedure, we can get $\tilde\omega$ such that
$$
\nu_D(\tilde\omega)>nm,\quad \nu_C(\tilde\omega)=\lambda_1.
$$
Since $\lambda_1=\nu_C(\tilde\omega)\geq \nu_D(\tilde\omega)>nm$; this should imply that $\lambda_1\in \Lambda_0$, contradiction.

As a consequence, we have that $\operatorname{In}(\omega_1)=\mu(nxdy-mydx)$ as desired. The above arguments are generalized in \cite{Can-C-SS-2023} (Proposition B.1) to obtain the following statement:

\begin{theorem}
\label{teo:standard:delorme}
For each $1\leq i\leq s$ we have the following statements:
\begin{enumerate}
\item $\lambda_i=\sup\{\nu_C(\omega); \; \omega\in \Omega^1_{\mathbb C^2,\mathbf 0},\text{ with }\nu_D(\omega)=t_i\}$.
\item If $\nu_C(\omega)=\lambda_i$, then $\nu_{D}(\omega)=t_i$.
\item For each $1$-form $\omega$ with $\nu_C(\omega)\notin \Lambda_{i-1}$, there is a unique pair $a,b\geq 0$ such that $\nu_{D}(\omega)=\nu_{D}(x^ay^b\omega_i)$. Moreover, we have that $\nu_C(\omega)\geq \lambda_i+na+mb$.
\item Let $k=\lambda_i+na+mb$, then $k\notin \Lambda_{i-1}$ if and only if for all $\omega$ such that $\nu_C(\omega)=k$ we have that $\nu_{D}(\omega)\leq \nu_{D}(x^ay^b\omega_i)$.
\item We have that $\lambda_i>u_i$.
\end{enumerate}
\end{theorem}
Let us note that the critical values $t_i$ of $\Lambda$ correspond exactly to the divisorial values of the elements $\omega_i$ of any minimal standard basis. Let us also note that the semimodule of differential values $\Lambda$ is an increasing cuspidal semimodule.
\begin{corolario}
 \label{cor:partiesinicialesenlabase}
 For each $1\leq i\leq s$, the $1$-forms $\omega_i$ are resonant.  In particular, taking an adapted coordinate system $(x,y)$, the initial parts can be written as
$$
\operatorname{In}(\omega_i)=\mu_i x^{e_i} y^{f_i}\left(n\frac{dx}{x}-m\frac{dy}{y}\right),
\quad ne_i+ mf_i=\nu_D(\omega_i)=t_i<nm.
$$
\end{corolario}
\begin{proof} Applying Lemma \ref{lema:lambdasytes}, since  $\Lambda$ is increasing, we have that
$$
\nu_D(\omega_i)=t_i<\lambda_i=\nu_{C}(\omega_i)<nm.
$$
The statement follows from these inequalities.
\end{proof}

\subsection{Generalized Delorme's Decomposition}

In this subsection, we state and prove a decomposition result for 1-forms which generalizes Delorme's decomposition \cite{delorme} and Theorem 8.5 in \cite{Can-C-SS-2023}:

Along this subsection, we fix a cusp $C$, with semimodule $\Lambda$ and basis
$\mathcal B=(\lambda_{-1},\ldots,\lambda_s)$.
We consider a minimal standard basis $\mathcal S=(\omega_{-1},\omega_0,\ldots,\omega_s)$ of $C$.
We also fix an element $*\in \{n,m\}$ (that is $*$ is either equal to $n$ or to $m$). Let us recall that we denote by
$$
k_i,\quad 0\leq i\leq s
$$
the bounds corresponding to the axes $u_{i+1}$, as  introduced in Section \ref{sec:semi}. In the same way, we denote by $k_i^*$ the bounds corresponding to the axes $u^*_{i+1}$. That is, we have
\begin{equation}
k_i=\left\{
\begin{array}{ccc}
k_i^n,&\text{ if }& u_{i+1}=u_{i+1}^n,\\
k_i^m,&\text{ if }& u_{i+1}=u_{i+1}^m;
\end{array}
\right.
\quad
k_i^*=\left\{\begin{array}{ccc}
k_i^n,&\text{ if }& u_{i+1}^*=u_{i+1}^n,\\
k_i^m,&\text{ if }& u_{i+1}^*=u_{i+1}^m.
\end{array}\right.
\end{equation}

\begin{theorem}
\label{teo:standard:decomposition3}
Consider indices $0\leq j\leq i\leq s$ and let us give a $1$-form
$\omega$ such that $\nu_D(\omega)=t_{i+1}^*$ and $\nu_C(\omega)>u_{i+1}^*$. Then, there is a decomposition of the $1$-form $\omega$ given by
\begin{equation}
\textstyle
\omega=\sum_{\ell=-1}^j f_\ell^{ij}\omega_\ell,
\end{equation}
 such that the following properties hold, where $v_{ij}^*=\nu_C(f_j^{ij}\omega_j)$:
\begin{enumerate}
\item $v_{ij}^*=\min\{\nu_C(f_\ell^{ij}\omega_\ell); -1\leq \ell<j\}$.
\item $v_{ij}^*=\lambda_j+t_{i+1}^*-t_j$, in particular, if $j=i$ we have that
$v_{ii}^*=\lambda_i+t_{i+1}^*-t_i=u_{i+1}^*$.
\item If $j<i$, we have that $\nu_{C}(f_\ell^{ij}\omega_\ell)=v_{ij}^*$, for $\ell=k_j$ and
$\nu_{C}(f_\ell^{ij}\omega_\ell)>v_{ij}^*$, for any $\ell\ne k_j$ and $-1\leq\ell<j$.
\item If $j=i$,  we have that $\nu_{C}(f_\ell^{ii}\omega_\ell)=v_{ii}^*$, for $\ell=k^*_j$ and
$\nu_{C}(f_\ell^{ii}\omega_\ell)>v_{ii}^*$, for any $\ell\ne k^*_j$ and $-1\leq\ell<j$.
\end{enumerate}
\end{theorem}
\begin{remark}
 \label{rk:descompositionofbasiselements}
 Let $\mathcal E=(\omega_{-1},\omega_0,\omega_1,\ldots,\omega_{s},\omega_{s+1})$ be an extended standard basis for $C$, that is, $(\omega_{-1},\omega_0,\omega_1,\ldots,\omega_{s})$ is a minimal standar basis of $C$ and $\omega_{s+1}$ is a 1-form with $\nu_D(\omega_{s+1})=t_{s+1}$ and $C$ is invariant for $\omega_{s+1}$ (see section~\ref{sec:saito}). Notice that Theorem \ref{teo:standard:decomposition3} can be used to write $\omega_{i+1}$ in terms of $\omega_{-1},\omega_0,\ldots,\omega_{i}$ for $0\leq i\leq s$. Indeed, let us choose $*\in\{n,m\}$ be such that $u^*_{i+1}=u_{i+1}$, and hence $t^*_{i+1}=t_{i+1}$. We have that $\nu_D(\omega_{i+1})=t_{i+1}$ and
$$
\nu_C(\omega_{i+1})=\left\{
\begin{array}{ccc}
\lambda_{i+1}>u_{i+1},&\text{ if }& i\leq s-1,\\
\infty>u_{s+1}, &\text{ if }& i=s.\\
\end{array}
\right.
$$
Now, by a direct application of Theorem \ref{teo:standard:decomposition3}, if we fix $j$ with $0\leq j\leq i$ we have an expression
\begin{equation}
\label{eq:descomposiciondeomegaimasuno}
\omega_{i+1}=f^{ij}_j\omega_{j}+f^{ij}_{j-1}\omega_{j-1}+\cdots f^{ij}_0\omega_0+f^{ij}_{-1}\omega_{-1},
\end{equation}
such that $
\lambda_j+t_{i+1}-t_j=
\nu_C(f^{ij}_j\omega_{j})=\nu_C(f^{ij}_{k_j}\omega_{k_j})<\nu_C(f^{ij}_{\ell}\omega_{\ell})
$, for any $\ell\ne k_j$, with $-1\leq \ell\leq j-1$.
\end{remark}

Before starting the proof of Theorem \ref{teo:standard:decomposition3}, we state the following Lemma, whose proof follows closely the proof of (\cite{Can-C-SS-2023}, Lemma C.1).

\begin{lema}\label{lem:standard:subida1} Consider $0\leq i\leq s$.  Given a 1-form $\eta$ with $\nu_C(\eta)> u_{i+1}^*$ and $\nu_D(\eta)> t^*_{i+1}$,  we have that:
	\begin{enumerate}
		\item If $\nu_D(\eta)<nm$, there is a 1-form $\alpha$ such that:
			\begin{enumerate}
			\item $\nu_D(\eta-\alpha)>\nu_D(\eta)$.
			\item There is a decomposition  $\alpha=\sum_{\ell=-1}^ig_\ell\omega_\ell$, where $\nu_C(g_\ell\omega_\ell)>u_{i+1}^*$ and $\nu_D(g_\ell\omega_\ell)>t_{i+1}^*$,  for any $-1\leq \ell \leq i$.
			\end{enumerate}
		\item If $\nu_D(\eta)\geq nm$, there is a decomposition $\eta=\sum_{\ell=-1}^ih_\ell\omega_\ell$ where each term satisfies that $\nu_C(h_\ell\omega_\ell) > u_{i+1}^*$.
	\end{enumerate}
\end{lema}
\begin{proof} Let us prove first statement (2).  Since $\{\omega_{-1},\omega_0\}$ is a basis of $\Omega_{(\mathbb{C}^2,0)}^1$, we can write in a unique way
\begin{equation}
\label{eq:descomposicionconnm}
\eta=g_{-1}\omega_{-1}+g_0\omega_0.
\end{equation}
Moreover, since $\operatorname{In}(\omega_{-1})=\lambda dx$ and $\operatorname{In}(\omega_{-1})=\mu dy$, we have that
$$
\nu_D(\eta)=\min\{\nu_D(g_{-1}\omega_{-1}),\, \nu_D(g_{0}\omega_{0})\}
$$
Noting that $\nu_D(\omega)\geq nm$, we have that $\nu_D(g_{-1}\omega_{-1})\geq nm$ and $\nu_D(g_{0}\omega_{0})\geq nm$. By Lemma~\ref{lem:semi:u-tilde}, we have that $\tilde u_{i+1}<\tilde u_1$, besides $u_{i+1}^*\leq \tilde u_{i+1}$, hence
$$
u_{i+1}^*\leq \tilde u_{i+1}<\tilde u_1=nm.
$$
We conclude that $\nu_C(g_\ell\omega_\ell)\geq \nu_D(g_\ell\omega_\ell)\geq nm>u_{i+1}^*$, for $\ell=-1,0$. Then the decomposition in equation \eqref{eq:descomposicionconnm} satisfies the required properties.

Let us prove now statement (1).  By Remark~\ref{rk:partieinicialvalorbajo}, the initial part of $\eta$ has a single term:
	$$
	\operatorname{In}(\eta)=x^a y^b \left( \mu_{-1}\frac{dx}{x}+\mu_0\frac{dy}{y}\right); \qquad \nu_D(\eta)=n a+m b.
	$$
There are two possibilities: either $\eta$ is resonant or not. If $\eta$ is not resonant, we have that $\nu_C(\eta)=\nu_D(\eta)=na+mb$. Let us consider
\begin{equation}
\label{eq:descomposicionnoresonante}
\alpha=\operatorname{In}(\eta)=g_{-1}\operatorname{In}(\omega_{-1})+g_0\operatorname{In}(\omega_{0}),\quad g_{-1}=\mu_{-1} x^{a-1}y^b,\; g_0=\mu_0 x^ay^{b-1}.
\end{equation}
We have that $\nu_D(\eta-\alpha)>\nu_D(\eta)$.
Moreover we also have that $\nu_D(\alpha)=\nu_D(\eta)>u^*_{i+1}$. Since
$$
\nu_D(\alpha)=\min\{\nu_D(g_{-1}\omega_{-1}),\nu_D(g_0\omega_0)\},
$$
we conclude that $\nu_C(g_\ell\omega_\ell)\geq \nu_D(g_\ell\omega_\ell)\geq \nu_D(\alpha)>u_{i+1}^* $, for $i=-1,0$. Moreover, in view of Corollary \ref{cor: ejesyvalorescriticos}  we have that $u_{i+1}^*>t_{i+1}^*$, hence we also get that $\nu_D(g_\ell\omega_\ell)>t_{i+1}^*$, for $\ell=-1,0$. Thus, the expression in equation \eqref{eq:descomposicionnoresonante} satisfies the desired properties.

Now, let us assume that $\eta$ is resonant. Up to multiply $\eta$ by a non-null scalar, we have that
$$
\operatorname{In}(\eta)=x^ay^b\left(n\frac{dx}{x}-m\frac{dy}{y}\right),\quad \nu_D(\eta)=na+mb>t_{i+1}^*.
$$
 Let us define the index
		$
		k:=\max\{\ell\leq i: \eta\text{ is reachable by }\omega_\ell\}
		$.
Since $\eta$ is resonant, then $k\geq 1$. By definition of $k$, there exists  a monomial
$\mu x^c y^d$ such that $\nu_D(\mu x^cy^d\omega_k)=\nu_D(\eta)$ and
$$
\nu_D(\eta')>\nu_D(\eta),\quad \text{ where }
\eta'=\eta-\mu x^cy^d\omega_k.
$$
The desired decomposition will be given just by the expression
 $
 \alpha=\mu x^cy^d\omega_k
 $.
 Since $\nu_D(\alpha)=\nu_D(\eta)>t_{i+1}^*$, we only need to verify that $\nu_C(x^cy^d\omega_k)>u_{i+1}^*$.
		
		First, let us  assume that $k=i$ and hence $\alpha=\mu x^cy^d\omega_i$. Write
 $$
 \nu_D(\alpha)=nc+md+t_i = \nu_D(\eta)>t^*_{i+1}.
 $$
Recalling that $t^*_{i+1}=t_i+u^*_{i+1}-\lambda_i$, we obtain that $nc+md+\lambda_i>u^*_{i+1}$. Hence,  we conclude by noting that
 $$
 \nu_C(x^cy^d\omega_i)=cn+dm+\lambda_i>u_{i+1}^*.
 $$
		
Now, let us  consider the case when $1\leq k\leq i-1$. Assume by contradiction that $\nu_C(x^cy^d\omega_k)\leq u^*_{i+1}<\nu_C(\eta)$. Taking into account that $\eta'=\eta-\mu x^cy^d\omega_k$, we see the following:
		\begin{align*}
			\nu_D(\eta') & >\nu_D(x^cy^d\omega_k)=nc+md+t_k; \\
			\nu_C(\eta') & =\nu_C(x^cy^d\omega_k)=nc+md+\lambda_k.
		\end{align*}
By statement (4) in Theorem~\ref{teo:standard:delorme}, we have that $nc+md+\lambda_k\in \Lambda_{k-1}$. In view of  Lemma \ref{lema:propiedadescrecimiento}, this implies that either $c\geq \ell_{k+1}^n$ or $d\geq \ell_{k+1}^m$. There are four possibilities:
		\begin{align*}
			u_{k+1}&=\lambda_k+n\ell_{k+1}^n\text{ and } c\geq \ell_{k+1}^n;  & u_{k+1}&=\lambda_k+n\ell_{k+1}^n\text{ and } d\geq \ell_{k+1}^m; \\
			u_{k+1}&=\lambda_k+m\ell_{k+1}^m\text{ and } c\geq \ell_{k+1}^n;  & u_{k+1}&=\lambda_k+m\ell_{k+1}^m\text{ and } d\geq \ell_{k+1}^m.
		\end{align*}
		The cases from the first line behave in a similar way as those in the second one, therefore, we will only show what happens in the first two cases.

{\em Case $u_{k+1}=u_{k+1}^n=\lambda_k+n\ell_{k+1}^n\text{ and } c\geq \ell_{k+1}^n$.} In this case we have that $\eta$ is reachable from $x^{\ell_{k+1}^n}\omega_{k}$. If we show that $x^{\ell_{k+1}^n}\omega_{k}$ is reachable from $\omega_{k+1}$, we contradict the maximality of $k$, as desired. In view of Corollary \ref{cor:partiesinicialesenlabase}, noting that both $\omega_{k+1}$ and $\omega_k$ are resonant, it is enough to show that
$$
\nu_D(x^{\ell_{k+1}^n}\omega_{k})=\nu_D(\omega_{k+1}).
$$
We have that $\nu_D(\omega_{k+1})=t_{k+1}$ and $\nu_D(x^{\ell_{k+1}^n}\omega_{k})=t_k+n\ell_{k+1}^n$. Let us see that $t_{k+1}=t_k+n\ell_{k+1}^n$ in our case.  In a general way, we have that $t_{k+1}^n=t_k+n\ell_{k+1}^n$; moreover, the fact that $u_{k+1}=u_{k+1}^n$ implies also that
$t_{k+1}=t_{k+1}^n$ and we are done.

{\em Case $u_{k+1}=u_{k+1}^n=\lambda_k+n\ell_{k+1}^n\text{ and } d \geq \ell_{k+1}^m$. }
By Lemma~\ref{lem:semi:u-tilde}, we see that
		$$
		nc+md+\lambda_k\geq \lambda_k + m\ell_{k+1}^m = \tilde{u}_{k+1}\geq \tilde{u}_i>u^*_{i+1}.
		$$
This ends the proof.
\end{proof}

\begin{proof}[Proof of Theorem \ref{teo:standard:decomposition3}] 
Let us take $\omega$ being such that $\nu_D(\omega)=t^*_{i+1}$ and $\nu_C(\omega)>u^*_{i+1}$ as in the statement.
We will consider three cases:
$$
\text{a) }\, i=0;\quad \text{b) }\, i>0,\, j=i;\quad \text{c) }\, i>0,\,0\leq j<i.
$$

\noindent{\em Case a): $i=0$.} Since $\{\omega_{-1},\omega_0\}$ is a basis of $\Omega_{\mathbb{C}^2,\mathbf 0}^1$, the 1-form $\omega$ can be written as
\begin{equation}
\label{eq:descomposicioniigualcero}
\omega=f^{00}_{-1}\omega_{-1} +f^{00}_0\omega_0.
\end{equation}
 Looking at the computations in Example \ref{ej:semi:t_1}, we see that
  $t^*_1=u^*_1\leq nm$ and $k^*_0=-1$.
  Therefore, we need to prove that $\nu_C(f_{-1}^{00}\omega_{-1})=\nu_C(f_0^{00}\omega_0)=u_1^*$.

Recall that, up to constant, we have that $\mathrm{In}(\omega_{-1})=dx$ and $\mathrm{In}(\omega_0)=dy$, so one of the following cases occurs:
\begin{enumerate}
 \item[(i)] $\nu_{D}(f^{00}_{0}\omega_{0})=t_1^*$ and $\nu_{D}(f^{00}_{-1}\omega_{-1})\geq t_1^*$.
 \item[(ii)] $\nu_{D}(f^{00}_{-1}\omega_{-1})=t_1^*$ and $\nu_{D}(f^{00}_{0}\omega_{0})\geq t_1^*$.
 \end{enumerate}
Assume that we are in case (i). Since $\nu_{D}(f_0^{00})+\nu_D(\omega_0)\leq nm$, we have that
$$
\nu_D(f^{00}_0)<nm.
$$
This implies that $\nu_D(f^{00}_0)=\nu_C(f^{00}_0)$ (left to the reader).  Therefore, we can write
 $$
 \nu_C(f^{00}_0\omega_0)=\nu_C(f^{00}_0) +\nu_C(\omega_0)=\nu_D(f^{00}_0) +\nu_D(\omega_0)=
 \nu_D(f^{00}_0\omega_0)=t^*_1=u^*_1.
 $$
 Moreover, since $\nu_D(f^{00}_{-1}\omega_{-1})\geq t_1^*=u^*_1$, we have that
 $$
 \nu_C(f^{00}_{-1}\omega_{-1})\geq \nu_D(f^{00}_{-1}\omega_{-1})\geq t_1^*=u^*_1.
 $$
  Noting that $\nu_C(f^{00}_{-1}\omega_{-1} +f^{00}_0\omega_0)>u^*_1$,  that $\nu_C(f^{00}_0\omega_0)=u^*_1$ and that $\nu_C(f^{00}_0\omega_0)\geq u^*_1$, we conclude that $\nu_C(f^{00}_{-1}\omega_{-1})=\nu_C(f^{00}_{0}\omega_{0})=u^*_1$.

 We do a similar argument in the case that $\nu_C(f^{00}_{-1}\omega_{-1})=t_1^*$.

 \bigskip

\noindent{\em Case b):  $i>0$ and $j=i$ .} We do the proof in the case $*=n$, the case $*=m$ runs in a similar way. Note that:
$$
\nu_D(\omega)=t^n_{i+1}<nm,\quad \nu_D(\omega)=t^n_{i+1}<u^n_{i+1}<\nu_C(\omega),
$$
in view of Corollary \ref{cor:cotastes} and Corollary \ref{cor: ejesyvalorescriticos}.
We deduce that the $1$-form $\omega$ is resonant. Since $\omega_i$ is also resonant and we have that
$$
\nu_D(\omega)=\nu_D(x^{\ell_{i+1}^n}\omega_i),\quad  (\text{recall that }\; t_{i+1}^n=t_i+n\ell_{i+1}^n),
$$
we deduce that there is a non-null scalar $\mu\ne 0$ such that
$$
\operatorname{In}(\omega)=\mu\operatorname{In}(x^{\ell_{i+1}^n}\omega_i)=
\mu x^{\ell_{i+1}^n}\operatorname{In}(\omega_i).
$$
Thus, the 1-form $\eta_1=\omega-\mu x^{\ell_{i+1}^n}\omega_i$ satisfies the following two properties:
$$
\nu_D(\eta_1)>t^n_{i+1}, \quad \nu_C(\eta_1)=\nu_C( x^{\ell_{i+1}^n}\omega_i)=u^n_{i+1},\quad (\text{recall that } u_{i+1}^n=\lambda_i+n\ell_{i+1}^n).
$$
The second one comes from the fact that $\nu_C(\omega)>u_{i+1}^n$. Take the bound $k=k_i^n$ and the co-limit $b=b_{i+1}$. We recall that
$$
u_{i+1}^n=\lambda_i+n\ell_{i+1}^n=\lambda_k+mb.
$$
Hence, the $1$-form $y^b\omega_k$ satisfies that $\nu_C(y^b\omega_k)=u_{i+1}^n$.
On the other hand, the divisorial value $\nu_D(y^b\omega_k)$ is given by
$$
\nu_D(y^b\omega_k)=mb+t_k.
$$
Let us show that $\nu_D(y^b\omega_k)>t_{i+1}^n=\nu_D(\omega)$. We have
\begin{eqnarray*}
t_{i+1}^n<mb+t_k\Leftrightarrow
t_i+n\ell_{i+1}^n<t_k+mb \Leftrightarrow\\
t_i-t_k< mb-n\ell_{i+1}^n=mb-n\ell_{i+1}^n+u_{i+1}^n-u_{i+1}^n\Leftrightarrow\\
t_i-t_k< (u_{i+1}^n-n\ell_{i+1}^n)-(u_{i+1}^n-mb)=\lambda_i-\lambda_k.
\end{eqnarray*}
We conclude, since $t_i-t_k<\lambda_i-\lambda_k$ in view of Lemma
\ref{lema:lambdasytes}.

Take $\eta_2=\eta_1-\mu_2 y^b\omega_k$ such that $\nu_C(\eta_2)>u_{i+1}^n$. Note that $\nu_D(\eta_2)>t_{i+1}^n$. Applying Lemma
\ref{lem:standard:subida1}, we get a decomposition
$$
\eta_2=\omega-\mu x^{\ell_{i+1}^n}\omega_i -\mu_2 y^b\omega_k=
\sum_{\ell=-1}^i h_\ell\omega_\ell,\quad \nu_C(h_\ell\omega_\ell)>u_{i+1}^n,\; \nu_D(h_\ell\omega_\ell)>t_{i+1}^n,
$$
having the desired properties.

{\em Case c): $i>0$, $0\leq j<i$. }
Let us reason by inverse induction on $j$, recalling that we are done when $j=i$. By induction hypothesis,  we can decompose $\omega$ as:
\begin{equation}\label{eq:omega2}
\textstyle
\omega=\sum_{\ell=-1}^{j+1} f_\ell^{ij+1}\omega_\ell,
\end{equation}
where
$
v_{ij+1}^*=\nu_C(f_{j+1}^{ij+1}\omega_{j+1})=\min\{\nu_C(f_\ell^{ij+1}\omega_\ell); -1\leq \ell<j+1\}.
$
Notice that in the case b), we have proven the case where $j+1=i$.  In view of Remak
\ref{rk:descompositionofbasiselements}, we can apply case b) to $\omega_{j+1}$ to obtain a decomposition:
\begin{equation}\label{eq:omega-j2}
\textstyle
\omega_{j+1}=\sum_{\ell=-1}^{j} f_\ell^{jj}\omega_\ell,
\end{equation}
where
$
u_{j+1}=\nu_C(f_j^{jj}\omega_j)=\min\{\nu_C(f_\ell^{jj}\omega_\ell);\; \ell<j\},
$
and the minimum is only reached at the bound $k=k_{j+1}$.  If we substitute the expression of $\omega_{j+1}$ given in \eqref{eq:omega-j2} into the expression of $\omega$ given in \eqref{eq:omega2}, we obtain
\begin{equation}
\label{eq:descomposicioncasoc}
\textstyle
\omega=\sum_{\ell=-1}^{j} (f_\ell^{ij+1}+f_{j+1}^{ij+1}f_\ell^{jj})\omega_\ell.
\end{equation}
Let us show that equation \eqref{eq:descomposicioncasoc} gives the desired decomposition. In order to do this, we only have to show that
\begin{enumerate}
\item[i)] $
\nu_C((f_j^{ij+1}+f_j^{jj}f_{j+1}^{ij+1})\omega_j)=\nu_C((f_k^{ij+1}+f_k^{jj}f_{j+1}^{ij+1})\omega_k)=v^*_{ij}
$.
\item[ii)]  $\nu_C((f_\ell^{ij+1}+f_\ell^{jj}f_{j+1}^{ij+1})\omega_\ell) > v^*_{ij}$ for $\ell \neq j,k$.
\end{enumerate}
Recall that $v^*_{ij+1}=\lambda_{j+1}+t^*_{i+1}-t_{j+1}$ and $v^*_{ij}=\lambda_j+t^*_{i+1}-t_j$. Hence, by Lemma \ref{lema:lambdasytes}, we have that  $v^*_{ij}<v^*_{ij+1}$. Moreover, by the properties of the decomposition given in equation \eqref{eq:omega2}, we get that:
\begin{align}
\nu_C(f^{ij+1}_{j+1})&=v^*_{ij+1}-\lambda_{j+1}; \label{eq:nuC1} \\
\nu_C(f_\ell^{ij+1}\omega_\ell) & \geq v^*_{ij+1}>v^*_{ij}, \text{ for } \ell < j+1. \label{eq:nuC2}
\end{align}
Using the expression given in  \eqref{eq:nuC1} and the properties of the decomposition given in \eqref{eq:omega-j2}, it follows that:
\begin{align*}
\nu_C(f_{j+1}^{ij+1}f_\ell^{jj}\omega_\ell)=&\nu_C(f_{j+1}^{ij+1})+\nu_C(f_\ell^{jj}\omega_\ell) = \\
 =&v^*_{ij+1}-\lambda_{j+1}+\nu_C(f_\ell^{jj}\omega_\ell) \geq  \\
 \geq & v_{ij+1}^*-\lambda_{j+1}+u_{j+1},
\end{align*}
where the last inequality is an equality just for $\ell=j,k$.
Now, taking into account that $u_{j+1}=\lambda_j+t_{j+1}-t_j$ and that $v_{ij+1}^*=t_{i+1}^*+\lambda_{j+1}-t_{j+1}$, we obtain that
$$
v_{ij+1}^*-\lambda_{j+1}+u_{j+1} =\lambda_j+t^*_{i+1}-t_j=v^*_{ij}.
$$
Finally, since $\nu_C(f_\ell^{ij+1}\omega_\ell) >v^*_{ij}$  for $\ell < j+1$, by expression \eqref{eq:nuC2}, we get that
$$
\nu_C((f_\ell^{ij+1}+f_\ell^{jj}f_{j+1}^{ij+1})\omega_\ell) \geq  v^*_{ij},
$$
where, again, we have an equality just for $\ell = j,k$.
\end{proof}

\begin{corolario}\label{cor:standard:partes inciales} Consider $1\leq i\leq s$ and let $\omega$ be a $1$-form such that $\nu_D(\omega)=t_{i+1}^*$ and $\nu_C(\omega)>u_{i+1}^*$. For any decomposition
$$
\textstyle
\omega=\sum_{\ell=-1}^j f_\ell^{ij}\omega_\ell,\qquad 1\leq j\leq i.
$$
satisfying the stated properties in Theorem \ref{teo:standard:decomposition3}, we have that $\mathrm{In}(\omega)=\mathrm{In}(f_j^{ij}\omega_j)$.
\end{corolario}
\begin{proof}
We only need to show that $\nu_D(f_j^{ij}\omega_j)<\nu_D(f_\ell^{ij}\omega_\ell)$ for $\ell<j$. We know that $\nu_C(f_j^{ij}\omega_j)\leq \nu_C(f_\ell^{ij}\omega_\ell)$ for $\ell<j$. Besides, $\nu_D(\omega)=t_{i+1}^*<nm$, because $i\geq 1$. Therefore, we have that $nm>\nu_D(f_j^{ij})$, consequently, the monomial order and the differential value coincide, $\nu_D(f_j^{ij})=\nu_C(f_j^{ij})$. Furthermore:
$$
\nu_C(f_\ell^{ij}\omega_\ell)=\nu_C(f_\ell^{ij})+\lambda_\ell\geq \nu_C(f_j^{ij}\omega_j)=\nu_C(f_j^{ij})+\lambda_j.
$$
By Lemma~\ref{lema:lambdasytes}, we have that  $\lambda_j-\lambda_\ell>t_j-t_\ell$, thus
$$
\nu_C(f_\ell^{ij})+t_\ell> \nu_C(f_j^{ij}\omega_j)=\nu_C(f_j^{ij})+t_j.
$$
If $\nu_C(f_\ell^{ij})>nm$, then its monomial order is at least $nm$. Hence, we have that $$\nu_D(f_\ell^{ij}\omega_\ell)>t_{i+1}^*=\nu_D(f_j^{ij}\omega_j).$$Indeed, if $\nu_C( f_\ell^{ij})\leq nm$, we get that $\nu_C(f_\ell^{ij})=\nu_D(f_j^{ij})$. With this, we conclude that
$$
\nu_D(f_\ell^{ij}\omega_\ell)=\nu_D(f_\ell^{ij})+t_\ell>\nu_D(f_j^{ij})+t_j=\nu_D(f_j^{ij}\omega_j).
$$
\end{proof}
\begin{prop}
\label{cor:standard:partes inciales2} Let $\omega$ be a $1$-form such that $\nu_D(\omega)=t_{1}^*$ and $\nu_C(\omega)>u_{1}^*$. Let us write (in a unique way)
$$
\omega=f_{-1}\omega_{-1}+f_0\omega_0.
$$
Then, we have that
\begin{enumerate}
\item If $t_1^*=t_1$, we have that $\operatorname{In}(\omega)=\mu(mydx-nxdy)$, where $\mu\ne 0$.
\item If $t_1^*=\tilde t_1$, we have that $\operatorname{In}(\omega)=\mu(\operatorname{In}(df))$, where $f=0$ is a reduced equation of the cusp $C$.
\end{enumerate}
In particular, we have that $\operatorname{In}(\omega)=\operatorname{In}(
f_{-1}\omega_{-1})+\operatorname{In}(f_{0}\omega_{0})$.
\end{prop}
\begin{proof} If $t_1^*=t_1$, since $t_1=n+m$, we see that $\operatorname{In}(\omega)$ can be written as
$$
\operatorname{In}(\omega)=\mu_{-1}ydx-\mu_0xdy.
$$
Moreover, we have that $t_1=u_1=n+m$ and hence $\nu_C(\omega)>\nu_D(\omega)$. Hence $\omega$ is resonant and we are done.

If $t_1^*=\tilde t_1=nm$, we also have that $\tilde u_1=nm$. The initial part $\operatorname{In}(\omega)$ has the form
$$
\operatorname{In}(\omega)= \mu_{-1}x^{m-1}dx+\mu_0y^{n-1}dy.
$$
If this initial part is not a multiple of $\operatorname{In}(f)$, we get that $\nu_C(\omega)=nm$, contradiction.
\end{proof}

\section{Standard Systems and Saito Bases}\label{sec:saito}
\subsection{Standard Systems}
Consider a  minimal standard basis
 $\mathcal S=(\omega_{-1},\omega_0,\omega_1,\ldots,\omega_s)$ of $C$. In \cite{Can-C-SS-2023} we have found differential $1$-forms $\omega_{s+1}$ with the following properties:
 \begin{enumerate}
 \item $\nu_D(\omega_{s+1})=t_{s+1}$.
 \item The cusp $C$ is invariant for $\omega_{s+1}$, that is $\nu_C(\omega_{s+1})=\infty$.
 \end{enumerate}
 We call {\em extended standard basis} of $C$ to any sequence
$$
\mathcal E=(\omega_{-1},\omega_0,\omega_1,\ldots,\omega_s,\omega_{s+1}),
$$
such that $\mathcal S=(\omega_{-1},\omega_0,\omega_1,\ldots,\omega_s)$ is a minimal standard basis and $\omega_{s+1}$ satisfies to the above properties.

In next definition we present systems of differential $1$-forms, where the axes $\tilde u_{s+1}$, instead of $u_{s+1}$, are essential in their construction.
\begin{defn}\label{def:special-standard-system}
A {\em  standard system $(\mathcal E,\mathcal F)$} for a cusp  $C$ is the data of an extended standard basis
$
\mathcal E=(\omega_{-1},\omega_0,\omega_1,\ldots,\omega_s,\omega_{s+1})
$
and a family $\mathcal F=(\tilde\omega_1,\tilde\omega_2,\ldots,\tilde\omega_{s},\tilde\omega_{s+1})$ of $1$-forms
satisfying that
$$
\nu_D(\tilde\omega_{j})=\tilde t_j,\quad \nu_C(\tilde\omega_j)=\infty,\quad 1\leq j  \leq s+1.
$$
We say that a  standard system $(\mathcal E,\mathcal F)$ for $C$ is a {\em special standard system} if there are expressions
$\tilde\omega_{j}=h_j\omega_{s+1}+f_j\tilde\omega_{s+1}$, where $h_j,f_j\in\mathcal O_{\mathbb C^2,\mathbf 0}$ for any  $1\leq j\leq s$.
\end{defn}

\subsection{Saito Bases}
Let $C$ be a cusp. Let us denote by $\Omega^1_{\mathbb C^2,\mathbf 0}[C]$ the $\mathcal O_{\mathbb C^2,\mathbf 0}$-submodule of $\Omega^1_{\mathbb C^2,\mathbf 0}$ given by the $1$-forms $\omega$ such that $C$ is invariant for $\omega$, that is, $\nu_C(\omega)=\infty$.

It is known that $\Omega^1_{\mathbb C^2,\mathbf 0}[C]$ is isomorphic to the $\mathcal O_{\mathbb C^2,\mathbf 0}$-module $\Omega^1_{\mathbb C^2,\mathbf 0}[\log C]$ of logarithmic meromorphic $1$-forms having poles along $C$. It is also known that these modules are free
$\mathcal O_{\mathbb C^2,\mathbf 0}$-module of rank two (see \cite{saito}). A basis of $\Omega^1_{\mathbb C^2,\mathbf 0}[C]$ will be called a {\em Saito basis for $C$}. The main result in this paper is the following one
\begin{theorem} \label{th:Saitobases}
Let $C$ be a cusp and let $\mathcal B=(\lambda_{-1},\lambda_{0},\lambda_1,\ldots,\lambda_s)$ be the basis of the semimodule $\Lambda$ of differential values for $C$. Denote by $t_{s+1}$ and $\tilde t_{s+1}$ the last critical values of $\Lambda$. Then,  there are two $1$-forms $\omega_{s+1},\tilde\omega_{s+1}$ having $C$ as an invariant curve and such that $\nu_D(\omega_{s+1})=t_{s+1}$ and $\nu_D(\tilde\omega_{s+1})=\tilde t_{s+1}$. Moreover, for any pair of $1$-forms as above, the set $\{\omega_{s+1},\tilde \omega_{s+1}\}$ is a Saito basis for $C$.
\end{theorem}

We prove Theorem \ref{th:Saitobases} in several steps:
\begin{enumerate}
\item We prove Theorem \ref{th:Saitobases} in the case $s=0$.
\item We show the existence of $\omega_{s+1},\tilde\omega_{s+1}$ having $C$ as invariant curve and such that $\nu_D(\omega_{s+1})=t_{s+1}$ and $\nu_D(\tilde\omega_{s+1})=\tilde t_{s+1}$.
\item  We show that $\mathcal F\cup \{\omega_{s+1}\}$ generates the $\mathcal O_{\mathbb C^2,\mathbf 0}$-module $\Omega^1_{\mathbb C^2,\mathbf 0}[C]$, for any  standard system $(\mathcal E,\mathcal F)$ that includes $\omega_{s+1}$ and $\tilde \omega_{s+1}$.
\item  We show that any pair of $1$-forms $\omega_{s+1},\tilde\omega_{s+1}$ having $C$ as invariant curve and such that $\nu_D(\omega_{s+1})=t_{s+1}$ and $\nu_D(\tilde\omega_{s+1})=\tilde t_{s+1}$ are included in at least one special  standard system $(\mathcal E,\mathcal F)$.
\item We conclude as follows. We start with $\{\omega_{s+1},\tilde \omega_{s+1}\}$ and we consider a special  standard system $(\mathcal E,\mathcal F)$ containing them. By statement (3), any $1$-form $\omega$ in the Saito module $\Omega^1_{\mathbb C^2,\mathbf 0}[C]$  is a combination
    $$
    \textstyle
    \omega=h\omega_{s+1}+\sum_{\ell=-1}^{s+1}f_\ell\tilde\omega_{\ell}.
    $$
    Since $(\mathcal E,\mathcal F)$ is a special standard system, each 1-form $\omega_\ell$ is a combination of $\omega_{s+1},\tilde \omega_{s+1}$, for any $\ell=-1,0,,1,\ldots,s$. In this way, we find a writing $\omega=f\omega_{s+1}+g\tilde\omega_{s+1}$, as desired.
\end{enumerate}

In next subsections we prove Theorem \ref{th:Saitobases}, following the above steps.
\subsection{The quasi-homogeneous case}
\label{The quasi-homogeneous case}
The statement of Theorem \ref{th:Saitobases} when $s=0$ is well know, see for instance \cite{saito}. Let us show it, for the sake of completeness. From Zariski's introduction of Zariski invariant, we know that the cusp $C$ is analytically equivalent to
the curve $f=0$, where $f=y^n-x^m$. We can take
$$
\omega_1= nxdy-mydx,\quad \tilde\omega_1=df=-mx^{m-1}dx+ny^{n-1}dy.
$$
Let us use Saito's criterion \cite{saito}
 which states that two $1$-forms $\omega,\tilde\omega$ in $\Omega_{\mathbb C^2,\mathbf 0}^1[C]$ give a Saito basis if and only if
 $$
 \omega\wedge\tilde\omega'=ufdx\wedge dy,
 $$
 where $u$ is a unit. Then $\omega_1$ and $\tilde\omega_1$ provide a Saito basis. Take now
 $\omega,\tilde\omega$ in $\Omega_{\mathbb C^2,\mathbf 0}^1[C]$  being such that
 $$
 \nu_D(\omega)=t_1=n+m,\quad \nu_D(\tilde \omega)=\tilde t_1=nm.
 $$
 Write
 $$
 \omega=A\omega_1+B\tilde\omega_1,\quad \tilde\omega=\tilde A\omega_1+\tilde B\tilde\omega_1.
 $$
 Since $n+m<nm$, we see that $A$ is a unit. It is also obvious that $\tilde A$ is not a unit. If we show that $\tilde B$ is a unit, the determinant $A\tilde B-B\tilde A$ is a unit and hence $\omega,\tilde\omega$ is a Saito basis. Note that
 $$
 \nu_D(\tilde A\omega_1)\ne nm.
 $$
 Indeed, if $\nu_D(\tilde A\omega_1)=\nu_D(\tilde A)+n+m=nm$, we conclude that
 $$
 \nu_D(\tilde A)=nm-n-m=c_\Gamma-1\in \Gamma.
 $$
 This is a contradiction. Then, we have that
 $$
 \nu_D(\tilde\omega)=nm=\nu_D(\tilde B\tilde\omega_1)=\nu_D(\tilde B)+nm.
 $$
 This implies that $\tilde B$ is a unit and we are done.

\subsection{Existence of $1$-forms with the last critical values}

In next proposition we show the existence of $\omega_{s+1},\tilde\omega_{s+1}$ with respective divisorial values $t_{s+1}$ and $\tilde t_{s+1}$ and such that $C$ is invariant for both $1$-forms. The proof follows the one in (\cite{Can-C-SS-2023}, Proposition 8.3).

\begin{prop} Let $C$ be a cusp and $\Lambda$ its semimodule of differential values, with basis $\mathcal B=(\lambda_{-1},\lambda_0,\lambda_1,\ldots,\lambda_s)$. Assume that $s\geq 1$. There are two $1$-forms $\omega_{s+1}$ and $\tilde\omega_{s+1}$ having $C$ as an invariant curve, such that $\nu_D(\omega_{s+1})=t_{s+1}$ and
 $\nu_D(\tilde \omega_{s+1})=\tilde t_{s+1}$.
\end{prop}
\begin{proof} Let us select a minimal standard basis
$
\mathcal S=(\omega_{-1},\omega_0,\omega_1,\ldots,\omega_s)
$
of the cusp $C$. As we have already done, let us denote by $*$ a chosen element $*\in\{n,m\}$. We have to find $\omega_{s+1}^*\in \Omega^1_{\mathbb C^2,\mathbf 0}[C]$ such that $\nu_D(\omega^*_{s+1})=t^*_{s+1}$.

We do the detailed proof for $*=n$. The case $*=m$ runs in a similar way. Then, we have to find $\omega^n_{s+1}\in \Omega^1_{\mathbb C^2,\mathbf 0}[C]$ such that $\nu_D(\omega^n_{s+1})=t^n_{s+1}$.

Let us recall that $u^n_{s+1}=\lambda_s+n\ell=\lambda_k+mb$, where we denote $\ell=\ell^n_{s+1}$, $k=k_{s}^n$ and $b=b_{s+1}$. Recall that $k<s$. Consider the $1$-forms
$$
\eta_0=x^\ell\omega_s,\quad  \eta_1= y^b\omega_k.
$$
Note that $\nu_D(\eta_0)=t^n_{s+1}=t_s+n\ell$ and $\nu_D(\eta_1)>t^n_{s+1}$. Indeed, we have
\begin{eqnarray*}
\nu_D(\eta_1)=bm+t_k>\nu_D(\eta_0)&=&t_s+n\ell\Leftrightarrow \\ t_s-t_k<bm-n\ell&=&(u^n_{s+1}-\lambda_k)-(u^n_{s+1}-\lambda_s)=\lambda_s-\lambda_k
\end{eqnarray*}
and we are done by Lemma \ref{lema:lambdasytes}. Moreover, the differential orders coincide
$$
\nu_C(\eta_0)=\nu_C(\eta_1)=u^n_{s+1}.
$$
Thus, there is a constant $\mu\ne 0$ such that if we take $\theta_1=\eta_0-\mu\eta_1$, we get that
$$
\nu_D(\theta_1)=t_{s+1}^n,\quad \nu_C(\theta_1)>\nu_C(\eta_0)=\nu_C(\eta_1)=u^n_{s+1}.
$$
We consider three cases:
\begin{enumerate}
\item[a)] $\nu_C(\theta_1)=\infty$. Then we are done by taking $\omega^n_{s+1}=\theta_1$.
\item[b)] $\nu_C(\theta_1)\geq nm$.
\item[c)] $\nu_C(\theta_1)<nm$.
\end{enumerate}
Assume that we are in case b) and let $\varphi$ be a primitive parametrization of $C$.
We have that $\varphi^*(\theta_1)=\psi(t)dt$, with $\operatorname{ord}_t(\psi(t))\geq nm-1> c_\Gamma$. In view of the classical theory of equisingularity \cite{wall}, there is a function $h(x,y)$ such that $\varphi^*(dh)=\psi(t)dt$. If we take $\omega^n_{s+1}=\theta_1-dh$, we have that $\nu_C(\omega^n_{s+1})=\infty$. In order to finish, we have to see that $\nu_D(dh)>t^n_{s+1}$. Since $t^n_{s+1}<\tilde t_1=nm$ (see Lemma \ref{lem:semi:u-tilde}), if we see that $\nu_D(dh)\geq nm$, we are done. If $\nu_D(dh)<nm$, we obtain that $\nu_C(dh)=\nu_D(dh)$, in contradiction with the fact that $\nu_C(dh)>nm$.

Assume now that we are in case c). Write $\nu_C(\theta_1)=\lambda_i+\alpha n+\beta m >u_{s+1}^n$, for a certain index $-1\leq i\leq s$. Consider the $1$-form $\eta_2$ given by
$$
\eta_2= x^\alpha y^\beta\omega_i,\quad \nu_D(\eta_2)= t_i+n\alpha+m\beta.
$$
Let us see that $\nu_D(\eta_2)>t_{s+1}^n=t_{s}+\ell n=\nu_D(\theta_1)$. Assume first that $i=s$, we know that $u^{n}_{s+1}=\lambda_s+n\ell<\nu_C(\eta_2)=\lambda_s+\alpha n+\beta m$, hence $n\alpha+m\beta>\ell n$ as desired. Assume now that $i<s$. We have
\begin{eqnarray*}
\nu_C(\eta_2)=\lambda_i+\alpha n+\beta m>u^{n}_{s+1}=\lambda_s+\ell n\Rightarrow\\
\Rightarrow
\alpha n+\beta m- \ell n>\lambda_s-\lambda_i>t_s-t_i\Rightarrow \\
\Rightarrow t_i+n\alpha+m\beta>t_s+\ell n.
\end{eqnarray*}
On the other hand, we have that $\nu_C(\eta_2)=\nu_C(\theta_1)$. Hence, there is a constant $\mu\ne 0$ such that if we take $\theta_2=\theta_1-\nu\eta_2$, we obtain that
$$
\nu_D(\theta_2)=\nu_D(\theta_1)=t^n_{s+1},\quad \nu_C(\theta_2)>\nu_C(\theta_1).
$$
We re-start the procedure with $\theta_2$, since the differential value is strictly increasing, in a finite number of steps we arrive to case b) or to case a) and we are done.
\end{proof}

\subsection{Generators of Saito Module}
Let us consider a  standard system $(\mathcal E,\mathcal F)$ of $C$, given by
$$
\mathcal E=(\omega_{-1},\omega_0,\omega_1,\ldots,\omega_{s},\omega_{s+1}),\quad
\mathcal F=(\tilde\omega_1,\tilde\omega_2,\ldots,\tilde\omega_s,\tilde\omega_{s+1}).
$$
In next Proposition \ref{prop:generadoresSaito} we describe a generator system of the Saito module $\Omega^1_{\mathbb C^2,\mathbf 0}[C]$.

Our arguments run by first considering the initial forms and finally by applying Artin's Approximation Theorem. Moreover, we work in an ordered way in terms of the divisorial values of the forms. In order to do this, we just need the concept of ``partial standard system''.

Consider an index $0\leq j\leq s$.
A {\em $j$-partial standard system associated to the extended standard basis $\mathcal E$} is a pair $(\mathcal E, \mathcal F^{j})$, where $\mathcal F^{j}$ is a list
$$
\mathcal F^{j}=(\tilde\omega_{j+1},\tilde\omega_{j+2},\ldots,\tilde\omega_{s+1}),
$$
such that $\nu_D(\tilde\omega_\ell)=\tilde t_\ell$ and $\omega_\ell\in \Omega^1_{\mathcal O_{\mathbb C^2,\mathbf 0}}[C]$, for $j+1\leq \ell\leq s+1$.

We start by a lemma concerning the structure of critical values:
\begin{lema}
\label{lema:valorescriticosespeciales}
Let $\Lambda$ be an increasing cuspidal semimodule of length $s\geq 1$.  Assume that the basis $
\mathcal B=(\lambda_{-1},\lambda_0,\lambda_1,\ldots,\lambda_s)
$,
satisfies that $\lambda_{-1}=n$ and $\lambda_0=m$. Consider the set
$$
T=\{t_{s+1},\tilde t_2,\tilde t_3,\ldots,\tilde t_{s+1}\},
$$
where $t_j,\tilde t_j$ are the critical values of $\Lambda$ corresponding to the index $j$.
Then, there are two nonnegative integer numbers $p,q\in \mathbb Z_\geq 0$ such that $$\{pn+n+m,qm +n+m\}\subset T.$$
Moreover, we have that $p<m-2$ and $q<n-2$.
\end{lema}
\begin{proof} We know that one of the following mutually excluding properties holds:
\begin{enumerate}
\item[(I)] $\tilde t_2=t_1+n\ell_{2}^n=n+m+n\ell_{2}^n$
\item[(II)] $\tilde t_2=t_1+m\ell_{2}^m=n+m+m\ell_{2}^m$
\end{enumerate}
Let us do the proof in the case (I), the case (II) has a similar proof. We can write
$
\tilde t_2=n+m+pn \in T
$, where $p=\ell_{2}^n$; thus, it is enough to find an element of $T$ of the form $n+m+qm$.
Assume first that $s=1$. Then $t_{s+1}=t_2=t_1+m\ell_{2}^m=n+m+m\ell_{2}^m$. Taking $q=\ell_2^m$, we have that $t_{s+1}=n+m+qm\in T$ and we are done.

Assume now that $s>1$. There are two cases:
\begin{enumerate}
    \item[a)] For any $2\leq i\leq s$, we have that $t_{i+1}-t_i=m\ell_{i+1}^m$.
    \item[b)] There is an index (that we take to be the minimum one) with $2\leq i\leq s$ such that $t_{i+1}-t_i=n\ell_{i+1}^n$.
\end{enumerate}
Assume we are in case a). Recall that $t_2=t_1+m\ell_{2}^m$, since $\tilde t_2=t_1+n\ell_{2}^n$. By a telescopic computation, we see that $t_{s+1}\in T$ may be written as
$$
\textstyle
t_{s+1}=t_1+\left(\sum_{\ell=2}^{s+1}\ell_{\ell}^m\right)m= n+m+qm,
$$
and we are done.

Assume we are in case b). For any $2\leq j\leq i$, we have that $t_j=t_{j-1}+m\ell_j^m$. By a telescopic computation, we obtain that
$
t_i=t_1+q_im
$.
The element $\tilde t_{i+1}\in T$ is given by $\tilde t_{i+1}=t_i+m\ell_{i+1}^m$ and hence we have that
$$
\tilde t_{i+1}=t_i+(m+\ell_{i+1}^m)m,
$$
as desired. This ends the proof.
\end{proof}
\begin{remark}
\label{rk:dosunoformasespeciales}
As a consequence of Lemma \ref{lema:valorescriticosespeciales}, we have the following property.  Assume that $\Lambda$ is the semimodule of differential values of a cusp $C$ and $(\mathcal E,\mathcal F)$ is a standard system, where
$$
\mathcal E=(\omega_{-1},\omega_0,\omega_1,\ldots,\omega_s,\omega_{s+1}),\quad
\mathcal F= (\tilde\omega_{1},\tilde \omega_2,\ldots,\tilde\omega_s,\tilde\omega_{s+1}).
$$
Consider the set
$
\mathcal T=\{\omega_{s+1},\tilde \omega_2,\ldots,\tilde\omega_s,\tilde\omega_{s+1}\}
$.
Assuming that $(x,y)$ is a system of adapted coordinates with respect to the cusp $C$, there are two $1$-forms $\eta_1,\eta_2\in\mathcal T$ such that
$$
\operatorname{In}(\eta_1)=\mu_1x^p(mydx-nxdy),\quad \operatorname{In}(\eta_2)=\mu_2y^q(mydx-nxdy),
$$
where $\mu_1\ne0\ne\mu_2$ and  $p,q\in \mathbb Z_{\geq 0}$.
\end{remark}
Next lemma is the key argument for finding our generator system of Saito's module. It will be also important in order to find the Saito bases we are looking for.

\begin{lema}
\label{lema:formasinicialestildes}
 Let us consider a standard system $(\mathcal E,\mathcal F)$ and a $1$-form
 $\omega\in \Omega^1_{\mathcal O_{\mathbb C^2,\mathbf 0}}[C]$.  Assume that $(x,y)$ is a system of adapted coordinates with respect to $C$. Then, the initial form $\operatorname{In}(\omega)$ is a linear combination, with quasi-homogeneous coefficients, of the initial forms
 $$
 \operatorname{In}(\tilde \omega_1), \ldots,\operatorname{In}(\tilde \omega_{s+1}), \operatorname{In}(\omega_{s+1}).
 $$
 \end{lema}
 \begin{proof} The initial form $W=\operatorname{In}(\omega)$ has the invariant curve $C_1$ given by $y^n=\mu x^m$, for a certain $\mu\ne 0$ (we leave this property to the reader). Let us invoke the result of Theorem \ref{th:Saitobases} for the case of length zero established in subsection \ref{The quasi-homogeneous case}. In this case we consider the two $1$-forms
$$
W_1=nxdy-mydx,\quad \tilde W_1=ny^{n-1}dy-\mu m x^{m-1}dx,
$$
that give a Saito basis $\{W_1,\tilde W_1\}$ of $C_1$. This gives a decomposition
$$
W=HW_1+\tilde G_1 \tilde W_1,
$$
where we can take $H, \tilde G_1$ to be quasi-homogeneous with respect to the weights $n,m$. By statement (2) of Proposition \ref{cor:standard:partes inciales2} and up to multiply $\tilde\omega_1 $ by a constant, we have that
$$
\operatorname{In}(\tilde\omega_1)=\tilde W_1.
$$
Now, we are going to show the existence of a decomposition
\begin{equation}
\label{eq:descompartesiniciales}
HW_1= G_{s+1}W_{s+1}+\sum_{\ell=2}^{s+1}\tilde G_\ell \tilde W_\ell,\quad\text{ where }\; \tilde W_\ell=\operatorname{In}(\tilde\omega_\ell), \; W_{s+1}=\operatorname{In}(\omega_{s+1}),
\end{equation}
with all the coefficients $G_{s+1}$ and $\tilde G_\ell$ being quasi-homogeneous.

Let $\delta=\nu_D(HW_1)$. Since  $H$ is a quasi-homogeneous polynomial, we can write
$$
HW_1= \sum_{_{\alpha n+\beta m=\delta}}W_{\alpha\beta},\quad W_{\alpha\beta}=\mu x^{\alpha}y^{\beta}\left(n\frac{dy}{y}-m\frac{dx}{x}\right)
,\quad  \alpha,\beta\in \mathbb Z_{\geq 1}.
$$
Now, it is enough to show that each of the $1$-forms $W_{\alpha\beta}$ is reachable by one of the $1$-forms in the set
$$
\mathcal T=\{\omega_{s+1},\tilde \omega_2,\ldots,\tilde \omega_{s+1}\}.
$$
We consider two cases:
\begin{enumerate}
\item[a)] There is a differential monomial $W_{\alpha\beta}\ne 0$ such that $\alpha\geq m$ or $\beta\geq n$.
\item[b)] For any $W_{\alpha\beta}\ne 0$ we have that $\alpha<m$ and $\beta<n$.
\end{enumerate}
Assume we are in case a). By a straightforward verification, we see that all the monomials $W_{\alpha\beta}\ne 0$ satisfy the condition that either  $\alpha\geq m$ or $\beta\geq n$.
In view of Lemma \ref{lema:valorescriticosespeciales} and Remark \ref{rk:dosunoformasespeciales}, we see that each $W_{\alpha\beta}\ne 0$ is reachable by an element of $\mathcal T$.

Assume now that we are in case b). Then, there is only one monomial $W_{\alpha\beta}\ne 0$ and hence, we have
$$
HW_1=\mu x^{\alpha-1} y^{\beta-1}(mydx-nxdy), \quad  1\leq \alpha<m,\, 1\leq \beta<n.
$$
Moreover, we have that $\tilde G_1\tilde W_1=0$. Indeed, we know that
$$
\tilde G_1\tilde W_1=\tilde G_1(ny^{n-1}dy-\mu m x^{m-1}dx)
$$
and, if this expression is nonzero, it contributes to monomials corresponding to case a), contradiction. We conclude that
$$
\operatorname{In}(\omega)=W=HW_1= \mu x^{\alpha-1} y^{\beta-1}(mydx-nxdy)=  \mu' x^{\alpha-1} y^{\beta-1}W_1.
$$
Note that $\omega$ is then reachable by $\omega_1$. Let $q$ be the maximum index $1\leq q\leq s+1$ such that $\omega$ is reachable by $\omega_q$. If $q=s+1$, we are done. Assume that $1\leq q\leq s$. Write
$$
\eta=\omega-\mu''x^ay^b\omega_q,\quad  \nu_D(\eta)>\nu_D(\omega).
$$
We have that $\nu_C(\eta)=\nu_C(x^ay^b\omega_q)$. We can invoke property (4) in Theorem \ref{teo:standard:delorme} to obtain that $\nu_C(\eta)\in \Lambda_{q-1}$, that is
$$
\lambda_q+na+mb\in \Lambda_{q-1}.
$$
By Lemma \ref{lema:propiedadescrecimiento}, we have that either $a\geq \ell_{q+1}^n$ or $b\geq \ell_{q+1}^m$. Assume that $a\geq \ell_{q+1}^n$. If $u_{q+1}=u_{q+1}^n$, then $\omega$ is reachable by $\omega_{q+1}$, contradiction the maximality of $q$, if $u_{q+1}=u_{q+1}^m$, we obtain that $\omega$ is reachable by $\tilde\omega_{q+1}$ and we are done. Same arguments for the case that $b\geq \ell_{q+1}^m$. This ends the proof.
\end{proof}
\begin{remark}
\label{rk:formasiniciales orden bajo}
Let $(\mathcal E, \mathcal F^{j})$, with $
\mathcal F^{j}=(\tilde\omega_{j+1},\tilde\omega_{j+2},\ldots,\tilde\omega_{s+1}),
$ be a $j$-partial standard system, with $j\geq 1$ and take a $1$-form $\omega\in \Omega_{\mathbb C^2,\mathbf 0}^1[C]$ such that $\nu_D(\omega)< \tilde t_j$. By the same arguments as in preceding lemma, noting that $\tilde t_j<\tilde t_{j-1}<\cdots<\tilde t_1$, we see that there is a combination
\begin{equation}
\label{eq:partialcombination}
\operatorname{In}(\omega)= G_{s+1}W_{s+1}+\sum_{\ell=j+1}^{s+1}\tilde G_\ell \tilde W_\ell,\quad\text{ where }\; \tilde W_\ell=\operatorname{In}(\tilde\omega_\ell), \; W_{s+1}=\operatorname{In}(\omega_{s+1}),
\end{equation}
all the coefficients being quasi-homogeneous of the corresponding degree.
\end{remark}

\begin{prop}
 \label{prop:generadoresSaito}
 The set $\mathcal T=\{\omega_{s+1},\tilde \omega_1,\tilde \omega_2,\ldots,\tilde \omega_{s+1}\}$ is a generator system of the Saito $\mathcal O_{\mathbb C^2,\mathbf 0}$-module $\Omega^1_{\mathbb C^2,\mathbf 0}[C]$.
\end{prop}
\begin{proof} Take $\omega\in \Omega^1_{\mathbb C^2,\mathbf 0}[C]$,  we know the existence of a decomposition
\begin{equation*}
\label{eq:descompartesiniciales2}
\operatorname{In}(\omega)= G_{s+1}W_{s+1}+\sum_{\ell=1}^{s+1}\tilde G_\ell \tilde W_\ell,\quad\text{ where }\; \tilde W_\ell=\operatorname{In}(\tilde\omega_\ell), \; W_{s+1}=\operatorname{In}(\omega_{s+1}),
\end{equation*}
with all the coefficients $G_{s+1}$ and $\tilde G_\ell$ being quasi-homogeneous.
We re-start the procedure of Lemma \ref{lema:formasinicialestildes}  with
$$
\omega'=\omega-\left(G_{s+1}\omega_{s+1}+\sum_{\ell=1}^{s+1}\tilde G_\ell \tilde \omega_\ell\right ).
$$
In this way, we obtain a formal expression
$
\textstyle
\omega= \hat g_{s+1}\omega_{s+1}+\sum_{\ell=1}^{s+1}\widehat{\tilde g}_\ell \tilde \omega_\ell
$.
By a direct application of Artin's Approximation Theorem \cite{artin}, we obtain the desired convergent expression
$$
\textstyle
\omega= g_{s+1}\omega_{s+1}+\sum_{\ell=1}^{s+1}\tilde{g}_\ell \tilde \omega_\ell.
$$
\end{proof}

\subsection{Existence of Special Standard Systems}

This subsection is devoted to provide a proof of the following result
\begin{prop}\label{prop:existencia-special-system}
 Assume that the length $s$ of the semimodule $\Lambda$ of differential values of the cusp $C$ is $s\geq 1$. Take two $1$-forms $\omega_{s+1}$ and $\tilde\omega_{s+1}$ in $\Omega^1_{\mathbb C^2,\mathbf 0}[C]$ such that $\nu_D(\omega_{s+1})=t_{s+1}$ and $\nu_D(\tilde\omega_{s+1})=\tilde t_{s+1}$. Then, there is a special standard system  $(\mathcal E,\mathcal F)$ for $C$ containing $\omega_{s+1},\tilde\omega_{s+1}$ in the sense that
$$
\mathcal E=(\omega_{-1},\omega_0,\omega_1,\ldots,\omega_{s},\omega_{s+1}),\quad
\mathcal F=(\tilde\omega_1,\tilde\omega_2,\ldots,\tilde\omega_s,\tilde\omega_{s+1}).
$$
\end{prop}
The proof of the above proposition follows directly from next result
\begin{prop}
 \label{prop:existenceofspecial2}
 Assume that the length $s$ of the semimodule $\Lambda$ of differential values of the cusp $C$ is $s\geq 1$. Take two $1$-forms $\omega_{s+1}$ and $\tilde\omega_{s+1}$ in $\Omega^1_{\mathbb C^2,\mathbf 0}[C]$ such that $\nu_D(\omega_{s+1})=t_{s+1}$ and $\nu_D(\tilde\omega_{s+1})=\tilde t_{s+1}$. For any index $1\leq j\leq s $ there are functions $f_j,\tilde f_j$ such that
$$
\nu_D(\tilde\omega_j)=\tilde t_j,
$$
where $\tilde\omega_j=f_j\omega_{s+1}+\tilde f_j\tilde\omega_{s+1}$.
\end{prop}
Along the whole proof, we consider an extended standard basis
$$\mathcal E=(\omega_{-1},\omega_0,\omega_1,\ldots,\omega_s,\omega_{s+1})$$ ending at $\omega_{s+1}$.
The proof of Proposition \ref{prop:existenceofspecial2} is quite long. In order to make clear the arguments, we do it in two steps:
\begin{itemize}
\item {\em Step 1: case $j=s$.} That is, we find $\tilde\omega_{s}\in \Omega_{\mathbb C^2,\mathbf 0}^1[C]$ such that $\nu_D(\tilde\omega_s)=\tilde t_s$.
\item Step 2: The general case.
\end{itemize}
\subsubsection{Existence of Special Standard Systems. First Case}
This subsection is devoted to the proof of Proposition \ref{prop:existenceofspecial2} when $j=s$. That is, we are going to prove that there is a linear combination
$$
\tilde\omega_s=\tilde f_s\tilde\omega_{s+1}+ f_s\omega_{s+1}
$$
such that $\nu_D(\tilde\omega_{s})=\tilde t_s$.

There are two possible cases: $t_{s+1}=t_s+n\ell^n_{s+1}$ and $t_{s+1}=t_s+m\ell^m_{s+1}$.
Both cases run in a similar way. We assume from now on that $t_{s+1}=t_s+n\ell^n_{s+1}$ and hence we have $\tilde t_{s+1}=t_s+m\ell^m_{s+1}$.
Let us write Delorme's decompositions of $\tilde \omega_{s+1}$ and $\omega_{s+1}$ as follows
\begin{eqnarray}
\label{eq:delormetildeomegaesemasuno}
\textstyle
\tilde\omega_{s+1}&=&\tilde\mu_1y^{\ell^m_{s+1}}\omega_s+
\tilde\mu_2x^{a_{s+1}}\omega_{k^m_s}+\tilde \eta,\quad \tilde\eta=\sum_{\ell=-1}^s\tilde h_\ell\omega_\ell,\\
\label{eq:delormeomegaesemasuno}
\omega_{s+1}&=&\mu_1x^{\ell^n_{s+1}}\omega_s+\mu_2y^{b_{s+1}}\omega_{k^n_s}+\eta,\quad \eta=\sum_{\ell=-1}^s h_\ell\omega_\ell
\end{eqnarray}
where we have the following properties:
\begin{enumerate}
\item $\operatorname{In}(\omega_{s+1})=\mu_1\operatorname{In}(x^{\ell^n_{s+1}}\omega_s)$. Recall that $t_{s+1}=t_s+n\ell_{s+1}^n$.
\item $\operatorname{In}(\tilde\omega_{s+1})=
    \tilde\mu_1\operatorname{In}(y^{\ell^m_{s+1}}\omega_s)$. Recall that $\tilde t_{s+1}=t_s+m\ell_{s+1}^m$.
\item $\nu_C(\mu_1x^{\ell^n_{s+1}}\omega_s+\mu_2y^{b_{s+1}}\omega_{k^n_s})>
\nu_C(\mu_1x^{\ell^n_{s+1}}\omega_s)=
\nu_C(\mu_2y^{b_{s+1}}\omega_{k^n_s})=u_{s+1}^n=u_{s+1}
$. Recall that $u_{s+1}^n=\lambda_s+n\ell^n_{s+1}=\lambda_{k^n_s}+mb_{s+1}$.
\item $\nu_C(\tilde\mu_1y^{\ell^m_{s+1}}\omega_s+\tilde\mu_2x^{a_{s+1}}\omega_{k^m_s})>
\nu_C(\mu_1y^{\ell^m_{s+1}}\omega_s)=
\nu_C(\mu_2x^{a_{s+1}}\omega_{k^m_s})=u_{s+1}^m=\tilde u_{s+1}
$. Recall that $u_{s+1}^m=\lambda_s+m\ell^m_{s+1}=\lambda_{k^m_s}+na_{s+1}$.
\item For any $-1\leq \ell\leq s$, we have that $\nu_C(h_\ell\omega_\ell)>u_{s+1}^n$ and
$\nu_C(\tilde h_\ell\omega_\ell)>u_{s+1}^m$.
\end{enumerate}

Let us consider the $1$-form $\theta_\mathbf 0\in \Omega^1_{\mathbb C^2,\mathbf 0}[C]$ defined by
$$
\theta_\mathbf 0=  \mu_1 x^{\ell_{s+1}^n}\tilde\omega_{s+1}- \tilde\mu_1 y^{\ell_{s+1}^m}\omega_{s+1}=\xi+\zeta_\mathbf 0,
$$
where
$
\xi= \tilde\mu_3x^{\ell_{s+1}^n+a_{s+1}}\omega_{k^m_s}- \mu_3 y^{\ell_{s+1}^m+b_{s+1}}\omega_{k^n_s}$, with $\tilde\mu_3=\mu_1\tilde\mu_2$, $\mu_3=\tilde \mu_1\mu_2
$ and such that
$
\textstyle
\zeta_\mathbf 0=\sum_{\ell=-1}^sg^\mathbf 0_\ell\omega_\ell
$.
In a more general way, given a pair of functions $\tilde f, f\in \mathcal O_{\mathbb C^2,\mathbf 0}$, we write
$$
\theta_{\tilde f,f}= \theta_\mathbf 0+\tilde f\tilde\omega_{s+1}+f\omega_{s+1}=\xi+\zeta_{\tilde f,f}\in  \Omega^1_{\mathbb C^2,\mathbf 0}[C],
$$
where $\zeta_{\tilde f,f}=\zeta_\mathbf 0+\tilde f\tilde\omega_{s+1}+f\omega_{s+1}$.
We also write
$
\zeta_{\tilde f,f}=\sum_{\ell=-1}^{s}g^{\tilde f,f}_\ell\omega_\ell
$.
Let us note that $\theta_\mathbf 0=\theta_{0,0}$, $\zeta_\mathbf 0=\zeta_{0,0}$ and  $g_\ell^{\mathbf 0}=g_\ell^{0,0}$, for $-1\leq \ell\leq s$.

In order to prove the desired result, we are going to show the existence of a pair $\tilde f, f$ such that  $\nu_D(\theta_{\tilde f,f})=\tilde t_s$.

We have two options: $u_s=u^n_s$ and $u_s=u^m_s$. Both cases run in a similar way. So, we fix the case that $u_s=u^n_s$. Hence, we have $t_s=t_s^n$, $\tilde u_s=u^m_s$ and  $\tilde t_{s}=t^m_s$. By Proposition \ref{prop:semi:limites}, we know  that $k^n_s=s-1$ and $k^m_s=k^m_{s-1}$.
\begin{lema}
	\label{lema:step1specialsystems0}
$\nu_D(\xi)=\tilde t_s$.
\end{lema}
\begin{proof}
By Proposition \ref{prop:semi:li+ai}, the co-limits $a_{s+1}$ and $b_{s+1}$ satisfy that
$
b_{s+1}+\ell_{s+1}^m=\ell_s^m$ and $a_{s+1}+\ell_{s+1}^n=a_s.
$
Hence,we have that
$$
\xi= -\mu_3 y^{\ell_{s}^m}\omega_{k^n_s}+\tilde\mu_3x^{a_s}\omega_{k^m_s}=
-\mu_3 y^{\ell_{s}^m}\omega_{s-1}+\tilde\mu_3x^{a_s}\omega_{k^m_{s-1}}.
$$
Let us show that $\nu_D(\xi)=\tilde t_s$.
Note that $\nu_D(y^{\ell_{s}^m}\omega_{s-1})=m\ell^m_s+t_{s-1}=t^m_s=\tilde t_s$. Thus, it is enough to show that $\nu_D(x^{a_s}\omega_{k^m_{s-1}})>\tilde t_s=t^m_s$. We have
$
\nu_D(x^{a_s}\omega_{k^m_{s-1}})=na_s+t_{k^m_{s-1}}
$.
Since  $u^m_s=m\ell_s^m+\lambda_{s-1}=na_s+\lambda_{k^m_{s-1}}$, then
\begin{eqnarray*}
na_s-m\ell^m_s&=&\lambda_{s-1}-\lambda_{k^m_{s-1}}>t_{s-1}-t_{k^m_{s-1}}\Rightarrow \\
&\Rightarrow& n a_s+t_{k^m_{s-1}}>\tilde t_s=t_{s-1}+m\ell^m_s.
\end{eqnarray*}
See Lemma \ref{lema:lambdasytes}. We conclude that $\nu_D(\xi)=\tilde t_s$.
\end{proof}

The problem is reduced to finding $\tilde f,f$ such that $\nu_D(\zeta_{f,\tilde f})>\tilde t_s$. Let us do it.

We say that a pair of functions $\tilde f,f$ is a {\em good pair} if and only if we have that $\nu_C(g^{\tilde f,f}_\ell\omega_\ell)>\tilde u_s$, for any $\ell=-1,0,\ldots,s$.

We end the proof as a direct consequence of the following lemmas:
\begin{lema}
\label{lema:step1specialsystems1}
The pair $\tilde f=0, f=0$ is a good pair.
\end{lema}
\begin{lema}
	\label{lema:step1specialsystems2}
	If $\tilde f,f$ is a good pair, then
	 $\nu_D(g_\ell^{\tilde f,f}\omega_\ell)>\tilde t_s$, for $-1\leq \ell\leq s-1$.
		and $\nu_D(g_s^{\tilde f,f}\omega_s)\ne \tilde t_s$.
\end{lema}
\begin{corolario}
	\label{cor:step1specialsystems2}
Assume that $\tilde f,f$ is a good pair. Then, we have that either
 $\nu_D(\theta_{\tilde f,f})=\tilde t_s$
 or
 $
 \nu_D(\theta_{\tilde f,f})=\nu_D(g^{\tilde f,f}_s\omega_s)<\tilde t_s
 $.
\end{corolario}
\begin{lema}
	\label{lema:step1specialsystems3}
		If $\tilde f,f$ is a good pair and $\nu_D(\theta_{\tilde f,f})<\tilde t_s$, then there is another good pair $\tilde f_1,f_1$ such that
	 $\nu_D(g_s^{\tilde f_1,f_1}\omega_s)>\nu_D(g_s^{\tilde f,f}\omega_s)$.
\end{lema}
Indeed, by Lemma \ref{lema:step1specialsystems1}, there is at least one good pair, by Lemma \ref{lema:step1specialsystems2} and Lemma \ref{lema:step1specialsystems0} we obtain Corollary \ref{cor:step1specialsystems2}. Now, we apply repeatedly Lemma \ref{lema:step1specialsystems3} to get that
$\nu_D(g^{\tilde f,f}_s\omega_s)\geq\tilde t_s$, hence, in view of Lemmas \ref{lema:step1specialsystems0} and \ref{lema:step1specialsystems1}, we get that  $\nu_D(g^{\tilde f,f}_s\omega_s)>\tilde t_s$ and $\nu_D(\theta_{\tilde f,f})=\tilde t_s$ as desired.

The rest of this subsection is devoted to proving the above three Lemmas \ref{lema:step1specialsystems1}, \ref{lema:step1specialsystems2} and \ref{lema:step1specialsystems3}.

\begin{proof}[Proof of Lemma \ref{lema:step1specialsystems1}]
We have to prove that $$
\nu_C(g_\ell^\mathbf 0\omega_\ell)> \tilde u_s, \text{ for any } \ell=-1,0,\ldots,s.
$$
Note that $\zeta_\mathbf 0=\mu_1 x^{\ell^n_{s+1}}\tilde \eta-\tilde \mu_1 y^{\ell^m_{s+1}} \eta$. Then, we have that
$
g_\ell^\mathbf 0= \mu_1 x^{\ell^n_{s+1}}\tilde h_\ell-\tilde \mu_1 y^{\ell^m_{s+1}} h_\ell$, for any  $\ell=-1,0,\ldots,s$. Now, it is enough to show that
$$
\nu_C(x^{\ell^n_{s+1}}\tilde h_\ell\omega_\ell)>\tilde u_s\; \text{ and }\;
\nu_C(y^{\ell^m_{s+1}} h_\ell\omega_\ell)>\tilde u_s.
$$
We have that
\begin{eqnarray*}
\nu_C(x^{\ell_{s+1}^n}\tilde h_\ell\omega_\ell)&=&n\ell^{n}_{s+1}+\nu_C(\tilde h_\ell\omega_\ell)
>n\ell^{n}_{s+1}+u_{s+1}^m=n\ell^{n}_{s+1}+\tilde u_{s+1}\\
&=&n\ell_{s+1}^n+na_{s+1}+\lambda_{k_s^m}= n(\ell_{s+1}^n+a_{s+1})+\lambda_{k_{s-1}^m}\\
&=& n a_s+\lambda_{k_{s-1}^m}=u_s^m=\tilde u_s\\
\nu_C(y^{\ell_{s+1}^m}h_\ell\omega_\ell)&=&m\ell^{m}_{s+1}+\nu_C(h_\ell\omega_\ell)
>m\ell^{m}_{s+1}+u_{s+1}^n=m\ell^{m}_{s+1}+u_{s+1}\\
&=&m\ell_{s+1}^m+mb_{s+1}+\lambda_{k_s^n}= m(\ell_{s+1}^m+b_{s+1})+\lambda_{s-1}\\
&=& m\ell_s^m+\lambda_{s-1}=u_s^m=\tilde u_s.
\end{eqnarray*}
This ends the proof of Lemma \ref{lema:step1specialsystems1}.
\end{proof}

\begin{proof}[Proof of Lemma \ref{lema:step1specialsystems2}]
Along the proof of this lemma, we just write $g^{\tilde f,f}_\ell=g_\ell$, in order to simplify the notation.

Let us first show that $\nu_D(g_\ell \omega_\ell)>\tilde t_s$,
for any $-1\leq \ell\leq s-1$. Recall that $\nu_C(g_\ell\omega_\ell)>\tilde u_s$ and write
$$
\nu_C(g_\ell\omega_\ell)=\nu_C(g_\ell)+\lambda_\ell>\tilde u_s=u_s^m=\lambda_{s-1}+m\ell^m_s.
$$
Noting that $\lambda_{s-1}-\lambda_\ell\geq t_{s-1}-t_\ell$, in view of Lemma \ref{lema:lambdasytes}, we have that
$$
\nu_C(g_\ell)+\lambda_{s-1}>\lambda_{s-1}+t_{s-1}-t_\ell+m\ell^m_s
$$
and thus we have
$
\nu_C(g_\ell)+t_\ell>t_{s-1}+m\ell^m_s=t^m_s=\tilde t_s
$.

Recall that $\nu_C(g_\ell)=\nu_D(g_\ell)$ when $\nu_D(g_\ell)<nm$. Noting that $\tilde t_s\leq nm$,  we conclude that
$$
\nu_D(g_\ell\omega_\ell)=\nu_D(g_\ell)+t_\ell> \tilde t_s,
$$
as desired.

Let us show that $\nu_D(g_s\omega_s)\ne \tilde t_s$. Assume by contradiction that $\nu_D(g_s\omega_s)=\tilde t_s$. Recalling that $t_s=t_s^n$, $\tilde t_s=t_s^m$, $t_s^n=t_{s-1}+n\ell^n_s$ and $t_s^m=t_{s-1}+m\ell^m_s$, we have
\begin{eqnarray*}
\nu_D(g_s\omega_s)=\tilde t_s
&\Rightarrow&\nu_D(g_s)+t_s=\tilde t_s \Rightarrow\nu_D(g_s)+t_s^n= t_s^m\Rightarrow \\
&\Rightarrow&\nu_D(g_s)+t_{s-1}+n\ell^n_s= t_{s-1}+m\ell^m_s\Rightarrow\\
&\Rightarrow& m\ell_s^m=\nu_D(g_s)+n\ell^n_s.
\end{eqnarray*}
This implies that $m\ell_s^m\in \Gamma$ is written in two different ways as a combination of $n,m$ with nonnegative integer coefficients. This is not possible, since $m\ell_s^m<nm$, in view of Remark \ref{rk:cotasparalimites}. The proof of Lemma \ref{lema:step1specialsystems2} is ended.
\end{proof}

\begin{proof}[Proof of Lemma \ref{lema:step1specialsystems3}]
 Assume that $\tilde f,f$ is a good pair with $\nu_D(\theta_{\tilde f,f})<\tilde t_s$. Let us find another good pair $\tilde f_1,f_1$ such that
$\nu_D(g_s^{\tilde f_1,f_1}\omega_s)>\nu_D(g_s^{\tilde f,f}\omega_s)$.

Since $\nu_D(\xi)=\tilde t_s$, $\theta_{\tilde f,f}=\xi+\zeta_{\tilde f,f}$ and $\nu_D(\theta_{\tilde f,f})<\tilde t_s$, we know that
$
\operatorname{In}(\theta_{\tilde f,f})=\operatorname{In}(\zeta_{\tilde f,f})$.
In particular $\nu_D(\zeta_{\tilde f,f})=\nu_D(\theta_{\tilde f,f})$. Applying Lemma \ref{lema:step1specialsystems2}, we get that
$$
\operatorname{In}(\theta_{\tilde f,f})=\operatorname{In}(\zeta_{\tilde f,f})=\operatorname{In}(g^{\tilde f,f}_s\omega_s)= \operatorname{In}(g^{\tilde f,f}_s)\operatorname{In}(\omega_s).
$$
Noting that $\nu_D(\theta_{\tilde f,f})<\tilde t_s\leq nm$ and $\nu_C(\theta_{\tilde f,f})=\infty$, we have that $\theta_{\tilde f,f}$ is a resonant $1$-form and, by the results in subsection~\ref{Initialparts}, there is a monomial $\mu x^ay^b$ such that
$$
\operatorname{In}(\theta_{\tilde f,f})= \mu x^ay^b\left(m\frac{dx}{x}-n\frac{dy}{y}\right),\quad a,b\geq 1,\;  na+mb=\nu_D(\theta_{\tilde f,f}).
$$
We conclude that there are $0\leq a'< a$, $0\leq b'<b$ and $\mu'\mu''=\mu$ such that
$$
\operatorname{In}(g^{\tilde f,f}_s)=\mu'x^{a'}y^{b'},\quad
\operatorname{In}(\omega_s)=\mu''x^{a-a'}y^{b-b'}\left(m\frac{dx}{x}-n\frac{dy}{y}\right).
$$
Let us consider the decomposition
$$
\theta_{\tilde f,f}= \mu'x^{a'}y^{b'}\omega_s+\eta',\quad \nu_D(\eta')>\nu_D(x^{a'}y^{b'}\omega_s).
$$
Noting that $\nu_C(\theta_{\tilde f,f})=\infty$, we have that $\nu_C(\eta')=\nu_C( \mu'x^{a'}y^{b'}\omega_s)=na'+mb'+\lambda_s$.
Let us apply Theorem \ref{teo:standard:delorme}, statement (4), to the integer number
$k=\lambda_s+na'+mb'$. Since there is $\eta'$ such that $\nu_C(\eta')=k$ and $\nu_D(\eta')> \nu_D(x^{a'}y^{b'}\omega_s)$, we conclude that $k\in\Lambda_{s-1}$. By Lemma \ref{lema:propiedadescrecimiento}, we know that one of the following properties holds:
$$
a'\geq \ell^n_{s+1}\quad \text{ or }\quad b'\geq \ell^m_{s+1}.
$$
Let us show that $\theta_{\tilde f,f}$ is reachable from $\omega_{s+1}$ or from $\tilde\omega_{s+1}$. Assume that $a'\geq \ell^n_{s+1}$, then we have that
\begin{eqnarray*}
\nu_D(\theta_{\tilde f,f})&=&an+bm=(a-a')n+(b-b')m+a'n+b'm\\
&=&t_s+n\ell^n_{s+1}+(a'-\ell^n_{s+1})n+b'm=t_{s+1}^n+(a'-\ell^n_{s+1})n+b'm.
\end{eqnarray*}
Noting that $t_{s+1}=t_{s+1}^n$, we have that $\theta_{\tilde f,f}$ and $x^{a'-\ell^n_{s+1}}y^{b'}\omega_{s+1}$ have the same initial parts (up to a constant) and thus $\theta_{\tilde f,f}$ is reachable from $\omega_{s+1}$. In the same way, if we assume that
$b'\geq \ell^m_{s+1}$, we have
\begin{eqnarray*}
\nu_D(\theta_{\tilde f,f})&=&an+bm=(a-a')n+(b-b')m+a'n+b'm\\
&=&t_s+m\ell^m_{s+1}+a'n+(b'-\ell^m_{s+1})m=t_{s+1}^m+a'n+(b'-\ell^m_{s+1})m\\
&=& \tilde t_{s+1}+a'n+(b'-\ell^m_{s+1})m.
\end{eqnarray*}
We conclude as above that $\theta_{\tilde f,f}$ is reachable from $\tilde\omega_{s+1}$.

Assume now that $a'\geq \ell^n_{s+1}$ and hence $\theta_{\tilde f,f}$ is reachable from $\omega_{s+1}$. Thus, there is a constant $\mu_3\ne 0$ such that
$$
\nu_D(\theta_{\tilde f,f}-\mu_3x^{a'-\ell^n_{s+1}}y^{b'}\omega_{s+1})>\nu_D(\theta_{\tilde f,f}).
$$
Let us put $\tilde f_1=\tilde f$ and $f_1=f-\mu_3x^{a'-\ell^n_{s+1}}y^{b'}$. Note that
$$
\theta_{\tilde f_1,f_1}= \theta_{\tilde f,f}-\mu_3x^{a'-\ell^n_{s+1}}y^{b'}\omega_{s+1}
$$
and hence $\nu_D(\theta_{\tilde f_1,f_1})>\nu_D(\theta_{\tilde f,f})$.

Let us verify that $\tilde f_1,f_1$ is a good pair. Let us write
$$
\textstyle
x^{a'-\ell^n_{s+1}}y^{b'}\omega_{s+1}=\sum_{\ell=-1}^s g'_\ell\omega_\ell,
$$
coming from the decomposition of $\omega_{s+1}$ in equation \eqref{eq:delormeomegaesemasuno}.
Noting that
$$
\zeta_{\tilde f_1,f}=\zeta_{\tilde f,f}-\mu_3x^{a'-\ell^n_{s+1}}y^{b'}\omega_{s+1},
$$
we see that $\tilde f_1,f_1$ is a good pair if  $\nu_C(g'_\ell\omega_\ell)>\tilde u_s$, for $\ell=-1,0,\ldots,s$. Let us show that this is true.
Since the terms $g_\ell'\omega_\ell$, for $-1\leq \ell\leq s$, come from the decomposition of $\omega_{s+1}$ times a monomial, we can apply Remark \ref{rk:descompositionofbasiselements} to see that
$$
\nu_C(g_s'\omega_s)\leq \nu_C(g_\ell'\omega_\ell),\quad \text{ for }-1\leq \ell\leq s.
$$
Hence, it is enough to show that $\nu_C(g_s'\omega_s)>\tilde u_s$.
Notice that
$$
\operatorname{In}(\zeta_{\tilde f,f})=\operatorname{In}(g_s^{\tilde f,f}\omega_s)=\mu_3\operatorname{In}(x^{a'-\ell^n_{s+1}}y^{b'}\omega_{s+1})=\mu_3\operatorname{In}(g_s'\omega_s),
$$
where the last equality comes from Corollary \ref{cor:standard:partes inciales}. Thus, we have $$\nu_D(g_s^{\tilde f,f}\omega_s)=\nu_D(g_s'\omega_s)<\tilde{t}_s\leq nm.$$ Therefore, $\nu_D(g_s^{\tilde f,f})=\nu_D(g_s')< nm$.  This implies that
$$
\nu_D(g_s^{\tilde f,f})=\nu_C(g_s^{\tilde f,f})=\nu_C(g_s')=\nu_D(g_s').
$$
Since $\tilde{f},f$ is a good pair, we conclude that $\nu_C(g_s')=\nu_C(g_s^{\tilde f,f})>\tilde u_s$. If $b'\geq \ell_{s+1}^m$, then $\theta_{\tilde{f},f}$ is reachable by $\tilde{\omega}_{s+1}$ and we proceed in a similar way.  This ends the proof of Lemma  \ref{lema:step1specialsystems3}.
\end{proof}

\subsubsection{Existence of Special Standard Systems. Induction Step}
This subsection is devoted to the proof of Proposition \ref{prop:existenceofspecial2} when $1\leq j<s$, assuming that the result is true for $j+1,j+2,\ldots,s$. That is, we are going to prove that there is a linear combination
$$
\tilde\omega_j=\tilde f_j\tilde\omega_{s+1}+ f_j\omega_{s+1}
$$
such that $\nu_D(\tilde\omega_{j})=\tilde t_j$, under the assumption that for any $j+1\leq \ell\leq s$ there is a linear combination
$
\tilde\omega_\ell=\tilde f_\ell\tilde\omega_{s+1}+ f_\ell\omega_{s+1},
$
such that $\nu_D(\tilde\omega_{\ell})=\tilde t_\ell$.

The proof is very similar to the case $j=s$. Recall that $\nu_D(\tilde \omega_{j+1})=\tilde t_{j+1}$. There are two options, either $\tilde t_{j+1}=t^n_{j+1}$ or $\tilde t_{j+1}=t^m_{j+1}$. In both cases, the proof runs in a similar way. So, we fix from now on the option $\tilde t_{j+1}=t^m_{j+1}$.

Let us define the number $q\in\{j+2,\ldots,s+1\}$ as follows
$$
q=
\left\{
\begin{array}{cc}
	{s+1}, \text{ if } \tilde t_\ell=t^n_\ell, \text{ for }\ell=j+2,j+3,\ldots,s+1,\\
	 \min\{\ell;\; \tilde t_\ell=t^m_\ell,\; j+2\leq \ell\leq s+1\},\text{ otherwise }.
\end{array}
\right.
$$
and define the $1$-form $\widehat\omega_q$ as follows:
$$
\widehat\omega_q=
\left\{
\begin{array}{cl}
\omega_{s+1},& \text{ if } \tilde t_\ell=t^n_\ell, \text{ for }\ell=j+2,j+3,\ldots,s+1,\\
\tilde\omega_q,&\text{ otherwise} .
\end{array}
\right.
$$
Let us note that $\nu_D(\widehat\omega_q)=t^m_q$ in both cases.

Now,
we proceed as follows:
\begin{enumerate}
	\item First, we find a  linear combination $\theta_{\mathbf 0}$ of $\tilde\omega_{j+1}$ and $\widehat\omega_q$ such that $\nu_D(\theta_{\mathbf 0})\leq\tilde t_j$. Note that $\theta_{\mathbf 0}$ should be a linear combination of $\tilde\omega_{s+1}$ and $\omega_{s+1}$, in view of the induction  hypothesis.
	\item Next, we find a $1$-form $\tilde\omega_j-\theta_{\mathbf 0}$ that is a linear combination of
	$$
	\tilde\omega_{j+1},\tilde\omega_{j+2},\ldots,\tilde\omega_{s+1},\omega_{s+1},
	$$
	in such a way that $\nu_D(\tilde\omega_j)=\tilde t_j$.
\end{enumerate}

Consider  Delorme's decompositions of $\tilde \omega_{j+1}$ and $\widehat \omega_q$ as introduced in Theorem \ref{teo:standard:decomposition3}, that we write  as follows
\begin{eqnarray}
	\label{eq:delormetildeomegaesemasuno1}
	\textstyle
	\tilde\omega_{j+1}&=&\tilde\mu_1y^{\ell^m_{j+1}}\omega_j+
	\tilde\mu_2x^{a_{j+1}}\omega_{k^m_j}+\tilde \eta,\quad \tilde\eta=\sum_{\ell=-1}^j\tilde h_\ell\omega_\ell,
	\\
	\label{eq:delormeomegaesemasuno1}
	\widehat \omega_q&=&M\omega_j+N\omega_{k^n_j}+\eta,\quad \eta=\sum_{\ell=-1}^j h_\ell\omega_\ell,
\end{eqnarray}
where $M,N$  are monomials in such a way that we have the following properties:
\begin{enumerate}
\item
$\operatorname{In}(\tilde\omega_{j+1})=
	\tilde\mu_1\operatorname{In}(y^{\ell^m_{j+1}}\omega_j)=\tilde \mu_1y^{\ell^m_{j+1}}\operatorname{In}(\omega_j)$. Recall that $\tilde t_{j+1}=t_j+m\ell_{j+1}^m$.
 \item
 $\nu_C(\tilde\mu_1y^{\ell^m_{j+1}}\omega_j+\tilde\mu_2x^{a_{j+1}}\omega_{k^m_j})>
	\nu_C(y^{\ell^m_{j+1}}\omega_j)=
	\nu_C(x^{a_{j+1}}\omega_{k^m_j})=u_{j+1}^m=\tilde u_{j+1}
	$. Recall that $u_{j+1}^m=\lambda_j+m\ell^m_{j+1}=\lambda_{k^m_j}+na_{j+1}$.
\item $\nu_C(\tilde h_\ell\omega_\ell)>\tilde u_{j+1}=u^m_{j+1}$, for $\ell=-1,0,1,\ldots,j$.
\item $\operatorname{In}(\widehat\omega_q)=\operatorname{In}(M\omega_j)=M\operatorname{In}(\omega_j)$.
\item
 $\nu_C(M\omega_j+N\omega_{k^n_j})>
	\nu_C(M\omega_j)=
	\nu_C(N\omega_{k^n_j})=\lambda_j+t_q^m-t_j=v^m_{q-1,j}
	$.
\item $\nu_C(h_\ell\omega_\ell)>\lambda_j+t_q^m-t_j=v^m_{q-1,j}$, for $\ell=-1,0,1,\ldots,j$.
\end{enumerate}
Let us compute the monomials $M$ and $N$. We have that
$$
t^m_q=\nu_D(\widehat\omega_q)=\nu_D(M)+\nu_D(\omega_j)\Rightarrow \nu_D(M)=t^m_q-t_j.
$$
By a telescopic argument, we obtain
\begin{eqnarray*}
t^m_q-t_j&=&t^m_q -t_{j+1}+(t_{j+1}-t_j)\\
&=& t^m_q -t_{j+1}+n\ell_{j+1}^n\\
&=&t^m_q-t_{j+2}+(t_{j+2} -t_{j+1})+n\ell_{j+1}^n\\
&=&t^m_q-t_{j+2}+m\ell_{j+2}^m+n\ell_{j+1}^n\\
&=&t^m_q-t_{j+3}+(t_{j+3} -t_{j+2})+m\ell_{j+2}^m+n\ell_{j+1}^n\\
&=&t^m_q-t_{j+3}+m(\ell_{j+3}^m+\ell_{j+2}^m)+n\ell_{j+1}^n\\
&\cdots&\cdots\\
&=&t_q^m-t_{q-1}+m(\ell^m_{q-1}+\cdots+\ell_{j+3}^m+\ell_{j+2}^m)+n\ell_{j+1}^n\\
&=&m(\ell^m_q+\ell^m_{q-1}+\cdots+\ell_{j+3}^m+\ell_{j+2}^m)+n\ell_{j+1}^n.
\end{eqnarray*}
This implies that $M=\mu_1x^ay^b$, where
$$
a=\ell^n_{j+1},\quad b= \ell^m_q+\ell^m_{q-1}+\cdots+\ell_{j+3}^m+\ell_{j+2}^m.
$$
Let us compute now the monomial $N$. We know that
$$
\nu_C(N\omega_{k_j^n})=\nu_D(N)+\lambda_{k_j^n}=\nu_C(M\omega_j)=\lambda_j+na+mb.
$$
Then, we have that
$$
\nu_D(N)=\lambda_j-\lambda_{k^n_j}+na+mb.
$$
Recalling that $u_{j+1}^n=\lambda_j+n\ell_{j+1}^n=\lambda_{k_j^n}+mb_{j+1}$, we obtain that
\begin{eqnarray*}
\nu_D(N)&=&\lambda_j-\lambda_{k^n_j}+na+mb=\\
&=&mb_{j+1}-n\ell_{j+1}^n+ na+mb=m(b_{j+1}+b).
\end{eqnarray*}
This implies that $N=\mu_2y^{b_{j+1}+b}$.

Let us note that
$b<\ell_{j+1}^m$,
in view of Corollary \ref{cor:cotaseles}. In a more precise way, we have that
$
\ell_{j+1}^m-b=b_q
$. Now, we consider the $1$-form $\theta_{\mathbf 0}$ given by
$$
\theta_{\mathbf 0}=\mu_1x^a \tilde\omega_{j+1}-\tilde\mu_1y^{b_q}\widehat\omega_q=
\mu_1x^{\ell_{j+1}^n}\tilde\omega_{j+1}-\tilde\mu_1y^{b_q}\widehat\omega_q.
$$
We write
$
\theta_{\mathbf 0}=\xi+\zeta_{\mathbf 0}
$,
where
\begin{eqnarray*}
\xi&=&\mu_1\tilde\mu_2x^{a+a_{j+1}}\omega_{k^m_j}-\tilde\mu_1\mu_2y^{b_q+b+b_{j+1}}\omega_{k^n_j}=\\
&=&\mu_1\tilde\mu_2x^{\ell_{j+1}^n+a_{j+1}}\omega_{k^m_j}-\tilde\mu_1\mu_2y^{\ell_{j+1}^m+b_{j+1}}\omega_{k^n_j}
\end{eqnarray*}
and $\zeta_\mathbf 0=\sum_{\ell=-1}^j g^{\mathbf 0}_\ell\omega_\ell=\sum_{\ell=-1}^j(\mu_1x^{\ell_{j+1}^n}\tilde h_\ell-\tilde\mu_1y^{b_q} h_\ell)\omega_\ell$.

In a more general way, given a list of functions $\widetilde{\mathbf f}, f$ in $\mathcal O_{\mathbb C^2,\mathbf 0}$, where
$$
\widetilde{\mathbf f}=(\tilde f_{j+1},\tilde f_{j+2},\ldots,\tilde f_{s+1}),
$$
 we write
$$
\textstyle
\theta_{\widetilde{\mathbf f},f}= \theta_\mathbf 0+\sum_{\ell=j+1}^{s+1}\tilde f_\ell\tilde\omega_{\ell}+f\omega_{s+1}=\xi+\zeta_{\widetilde{\mathbf f},f}\in  \Omega^1_{\mathbb C^2,\mathbf 0}[C],
$$
where $\zeta_{\widetilde{\mathbf f},f}=\zeta_\mathbf 0+\sum_{\ell=j+1}^{s+1}\tilde f_\ell\tilde\omega_{\ell}+f\omega_{s+1}$.
We also write
$
\zeta_{\widetilde{\mathbf f},f}=\sum_{\ell=-1}^{s}g^{\widetilde{\mathbf f},f}_\ell\omega_\ell
$.
Let us note that $\theta_\mathbf 0=\theta_{\mathbf 0,0}$, $\zeta_\mathbf 0=\zeta_{\mathbf 0,0}$ and  $g_\ell^{\mathbf 0}=g_\ell^{\mathbf 0,0}$, for $-1\leq \ell\leq s$.

In order to prove the desired result, we are going to show the existence of a list $\widetilde{\mathbf f}, f$ such that  $\nu_D(\theta_{\widetilde{\mathbf f},f})=\tilde t_j$.

We have two options: $u_j=u^n_j$ and $u_j=u^m_j$. Both cases run in a similar way. So, we fix the case that $u_j=u^n_j$. Hence, we have $t_j=t_j^n$, $\tilde u_j=u^m_j$ and  $\tilde t_{j}=t^m_j$. By Proposition \ref{prop:semi:limites}, we know  that $k^n_j=j-1$ and $k^m_j=k^m_{j-1}$.

\begin{lema}
	\label{lema:stepjspecialsystems0}
$\nu_D(\xi)=\tilde t_j$.
\end{lema}

Now, the problem is reduced to finding a list  $(\widetilde {\mathbf f},f)$ such that $\nu_D(\zeta_{\widetilde {\mathbf f},f})>\tilde t_j$. Let us do it.
We say that a list of functions $(\widetilde {\mathbf f},f)$ is a {\em good list} if and only if we have that $\nu_C(g^{\widetilde {\mathbf f},f}_\ell\omega_\ell)>\tilde u_j$, for any $\ell=-1,0,\ldots,j$.

We end the proof as a direct consequence of the following lemmas:
\begin{lema}
\label{lema:stepjspecialsystems1}
The list $(\widetilde {\mathbf f},f)=(\mathbf 0,0)$ is a good list.
\end{lema}
\begin{lema}
	\label{lema:stepjspecialsystems2}
	If $(\widetilde {\mathbf f},f)$ is a good list, then
	 $\nu_D(g_\ell^{\widetilde {\mathbf f},f}\omega_\ell)>\tilde t_j$, for $-1\leq \ell\leq j-1$ and $\nu_D(g_j^{\widetilde {\mathbf f},f}\omega_j)\ne \tilde t_j$.
\end{lema}
\begin{corolario}
	\label{cor:stepjspecialsystems2}
Assume that $(\widetilde {\mathbf f},f)$ is a good list. Then, either we have that
 $\nu_D(\theta_{\widetilde {\mathbf f},f})=\tilde t_j$
 or
 $
 \nu_D(\theta_{\widetilde {\mathbf f},f})=\nu_D(g^{\widetilde {\mathbf f},f}_j\omega_j)<\tilde t_j
 $.
\end{corolario}
\begin{lema}
	\label{lema:stepjspecialsystems3}
		If $(\widetilde {\mathbf f},f)$ is a good list and $\nu_D(\theta_{\widetilde {\mathbf f},f})<\tilde t_j$, then there is another good list $(\widetilde {\mathbf f}^1,f^1)$ such that
	 $\nu_D(g_j^{\widetilde {\mathbf f}^1,f^1}\omega_j)>\nu_D(g_j^{\widetilde {\mathbf f},f}\omega_j)$.
\end{lema}
Indeed, by Lemma \ref{lema:stepjspecialsystems1}, there is at least one good list, by Lemma \ref{lema:stepjspecialsystems2} and Lemma \ref{lema:stepjspecialsystems0} we obtain Corollary \ref{cor:stepjspecialsystems2}. Now, we apply repeatedly Lemma \ref{lema:stepjspecialsystems3} to get that
$\nu_D(g^{\widetilde {\mathbf f},f}_j\omega_j)\geq\tilde t_j$, hence, in view of Lemmas \ref{lema:stepjspecialsystems0} and \ref{lema:stepjspecialsystems1}, we get that  $\nu_D(g^{\widetilde {\mathbf f},f}_j\omega_j)>\tilde t_j$ and
$\nu_D(\theta_{\widetilde {\mathbf f},f})=\tilde t_j$ as desired.

The rest of this subsection is devoted to proving the above four Lemmas
\ref{lema:stepjspecialsystems0}, \ref{lema:stepjspecialsystems1}, \ref{lema:stepjspecialsystems2} and \ref{lema:stepjspecialsystems3}.

\begin{proof}[Proof of Lemma \ref{lema:stepjspecialsystems0}]
By Proposition \ref{prop:semi:li+ai}, the co-limits $a_{j+1}$ and $b_{j+1}$ satisfy that
$
b_{j+1}+\ell_{j+1}^m=\ell_j^m$ and $a_{j+1}+\ell_{j+1}^n=a_j.
$
Hence,we have that
\begin{eqnarray*}
\xi&=& \mu_1\tilde\mu_2x^{\ell_{j+1}^n+a_{j+1}}\omega_{k^m_j}
-\tilde\mu_1\mu_2y^{\ell_{j+1}^m+b_{j+1}}\omega_{k^n_j}\\
&=& \mu_1\tilde\mu_2x^{a_j}\omega_{k^m_{j-1}}
-\tilde\mu_1\mu_2y^{\ell_{j}^m}\omega_{j-1}.
\end{eqnarray*}
Note that $\nu_D(y^{\ell_{j}^m}\omega_{j-1})=m\ell^m_j+t_{j-1}=t^m_j=\tilde t_j$. Thus, it is enough to show that $\nu_D(x^{a_j}\omega_{k^m_{j-1}})>\tilde t_j=t^m_j$. We have
$
\nu_D(x^{a_j}\omega_{k^m_{j-1}})=na_j+t_{k^m_{j-1}}
$.
Since  $u^m_j=m\ell_j^m+\lambda_{j-1}=na_j+\lambda_{k^m_{j-1}}$, then
\begin{eqnarray*}
na_j-m\ell^m_j&=&\lambda_{j-1}-\lambda_{k^m_{j-1}}>t_{j-1}-t_{k^m_{j-1}}\Rightarrow \\
&\Rightarrow& n a_j+t_{k^m_{j-1}}>\tilde t_j=t_{j-1}+m\ell^m_j.
\end{eqnarray*}
See Lemma \ref{lema:lambdasytes}. We conclude that $\nu_D(\xi)=\tilde t_j$.
\end{proof}

\begin{proof}[Proof of Lemma \ref{lema:stepjspecialsystems1}]
 We have to prove that $$
\nu_C(g_\ell^\mathbf 0\omega_\ell)> \tilde u_j, \text{ for any } \ell=-1,0,\ldots,j.
$$
Note that $g_\ell^\mathbf 0\omega_\ell=(\mu_1x^{\ell_{j+1}^n}\tilde h_\ell-\tilde\mu_1y^{b_q} h_\ell)\omega_\ell$. Now, it is enough to show that
$$
\nu_C(x^{\ell^n_{j+1}}\tilde h_\ell\omega_\ell)>\tilde u_j\; \text{ and }\;
\nu_C(y^{b_q} h_\ell\omega_\ell)>\tilde u_j.
$$
We have
\begin{eqnarray*}
\nu_C(x^{\ell_{j+1}^n}\tilde h_\ell\omega_\ell)&=&n\ell^{n}_{j+1}+\nu_C(\tilde h_\ell\omega_\ell)
>n\ell^{n}_{j+1}+u_{j+1}^m=n\ell^{n}_{j+1}+\tilde u_{j+1}\\
&=&n\ell_{j+1}^n+na_{j+1}+\lambda_{k_j^m}= n(\ell_{j+1}^n+a_{j+1})+\lambda_{k_{j-1}^m}\\
&=& n a_j+\lambda_{k_{j-1}^m}=u_j^m=\tilde u_j.
\end{eqnarray*}
Let us consider now $\nu_C(y^{b_q} h_\ell\omega_\ell)$. We have that
$$
\nu_C(y^{b_q}h_\ell\omega_\ell)>mb_q+\lambda_j+t_q^m-t_j.
$$
Let us show that $mb_q+\lambda_j+t_q^m-t_j=\tilde u_j$. Recall that $\tilde u_j=u_j^m=\lambda_{j-1}+m\ell_j^m$. Thus, we have to prove that
$$
mb_q+\lambda_j+t_q^m-t_j-\lambda_{j-1}-m\ell_j^m=0.
$$
Note that $k_j^n=j-1$ and then $\lambda_j-\lambda_{j-1}=-n\ell_{j+1}^n+mb_{j+1}$. Then we have to verify that
$$
mb_q-n\ell_{j+1}^n+mb_{j+1}+t_q^m-t_j-m\ell_j^m=0.
$$
Recalling that $t_q^m-t_j=na+mb=
n\ell^n_{j+1}+mb$ and that $b_q=\ell_{j+1}^m-b$, we have to verify that
$$
m(\ell_{j+1}^m-b)-n\ell_{j+1}^n+n\ell^n_{j+1}+mb+mb_{j+1}-m\ell_j^m=0.
$$
We have to see that $b_{j+1}+\ell_{j+1}^m=\ell_j^m$, and this follows from Proposition \ref{prop:semi:li+ai}.
\end{proof}

\begin{proof}[Proof of Lemma \ref{lema:stepjspecialsystems2}]
 Along the proof of this lemma, we just write $g^{\widetilde{\mathbf f},f}_\ell=g_\ell$, in order to simplify the notation.

Let us first show that $\nu_D(g_\ell)>\tilde t_j$,
for any $-1\leq \ell\leq j-1$. Recall that $\nu_C(g_\ell\omega_\ell)>\tilde u_j$ and write
$$
\nu_C(g_\ell\omega_\ell)=\nu_C(g_\ell)+\lambda_\ell>\tilde u_j=u_j^m=\lambda_{j-1}+m\ell^m_j.
$$
Noting that $\lambda_{j-1}-\lambda_\ell\geq t_{j-1}-t_\ell$, in view of Lemma \ref{lema:lambdasytes}, we have that
$$
\nu_C(g_\ell)+\lambda_{j-1}>\lambda_{j-1}+t_{j-1}-t_\ell+m\ell^m_j
$$
and thus we have
$
\nu_C(g_\ell)+t_\ell>t_{j-1}+m\ell^m_j=t^m_j=\tilde t_j
$.

Recall that $\nu_C(g_\ell)=\nu_D(g_\ell)$ when $\nu_D(g_\ell)<nm$. Noting that $\tilde t_j\leq nm$,  we conclude that
$$
\nu_D(g_\ell\omega_\ell)=\nu_D(g_\ell)+t_\ell> \tilde t_j,
$$
as desired.

Let us show that $\nu_D(g_j\omega_j)\ne \tilde t_j$. Assume by contradiction that $\nu_D(g_j\omega_j)=\tilde t_j$. Recalling that $t_j=t_j^n$, $\tilde t_j=t_j^m$, $t_j^n=t_{j-1}+n\ell^n_j$ and $t_j^m=t_{j-1}+m\ell^m_j$, we have
\begin{eqnarray*}
\nu_D(g_j\omega_j)=\tilde t_j
&\Rightarrow&\nu_D(g_j)+t_j=\tilde t_j \Rightarrow\nu_D(g_j)+t_j^n= t_j^m\Rightarrow \\
&\Rightarrow&\nu_D(g_j)+t_{j-1}+n\ell^n_j= t_{j-1}+m\ell^m_j\Rightarrow\\
&\Rightarrow& m\ell_j^m=\nu_D(g_j)+n\ell^n_j.
\end{eqnarray*}
This implies that $m\ell_j^m\in \Gamma$ is written in two different ways as a combination of $n,m$ with nonnegative integer coefficients. This is not possible, since $m\ell_j^m<nm$, in view of Remark \ref{rk:cotasparalimites}. \end{proof}

\begin{proof}[Proof of Lemma \ref{lema:stepjspecialsystems3}]
 Assume that $\widetilde{\mathbf f},f$ is a good list with $\nu_D(\theta_{\widetilde{\mathbf f},f})<\tilde t_j$. Let us find another good list $\widetilde {\mathbf f^1},f^1$ such that
$$\nu_D(g_j^{\widetilde {\mathbf f^1},f^1}\omega_j)>\nu_D(g_j^{\widetilde{\mathbf f},f}\omega_j).$$

Let us note that $\nu_D(\theta_{\widetilde{\mathbf f},f})=\nu_D(g_j^{\widetilde{\mathbf f},f}\omega_j)<\tilde t_j$ and, more precisely, we have that
$$
W=\operatorname{In}(g_j\omega_j)= \operatorname{In}(\theta_{\widetilde{\mathbf f},f}).
$$
In view of Remark \ref{rk:formasiniciales orden bajo}, there is a decomposition
$$
W=G_{s+1}W_{s+1}+\sum_{\ell=j+1}^{s+1} \tilde G_\ell \tilde W_\ell,
$$
where the coefficients are quasi-homogeneous. Moreover, all the forms $W,W_{s+1},\tilde W_{\ell}$, for $j+1\leq \ell\leq s+1$ are resonant with divisorial value $<nm$. We conclude that all those forms are given by a monomial times the $1$-form
$$
m\frac{dx}{x}-n\frac{dy}{y}.
$$
Up to multiply some of the terms for an adequate scalar number, we can assume without loss of generality that all the coefficients $G_{s+1},\tilde G_{j+1},\tilde G_{j+2},\ldots,\tilde G_{s+1}$ are zero except exactly one of them. So, we have that
$$
W=G_{s+1}W_{s+1}\text{ or there is $\ell_0$ such that } W=\tilde G_{\ell_0}\tilde W_{\ell_0}.
$$
Let us write $S= W_{s+1}$ in the first case and $S=\tilde W_{\ell_0}$ in the second one. Then we have that $W=GS$, where $G=G_{s+1}$ in the first case and $G=\tilde G_{\ell_0}$ in the second one.

Now we define the list $(\widetilde{\mathbf f}^1,f^1)$ by
\begin{eqnarray*}
(\tilde f^1_{j+1},\tilde f^1_{j+1},\ldots,\tilde f^1_{s+1},f^1)=
(\widetilde{\mathbf f},f)-(\tilde G_{j+1},\tilde G_{j+2},\ldots,\tilde G_{s+1}, G_{s+1}).
\end{eqnarray*}
It is a straightforward verification that
$\nu_D(g^{\widetilde{\mathbf f}^1,f^1}_j\omega_j)>\nu_D(g^{\widetilde{\mathbf f},f}_j\omega_j)$.

We have just to verify that $(\widetilde{\mathbf f}^1,f^1)$ is a good list. We do it in the case that $S=W_{s+1}$, the other cases run in a similar way. Note that
$$
\zeta_{\widetilde{\mathbf f}^1,f^1}=\zeta_{\widetilde{\mathbf f},f}-G_{s+1}\omega_{s+1}=
\sum_{\ell=-1}^jg^{\widetilde{\mathbf f}^1,f^1}_\ell\omega_\ell.
$$
Let us perform a Delorme decomposition of $\omega_{s+1}$:
$
\omega_{s+1}=\sum_{\ell=-1}^jc_\ell\omega_\ell
$,
where we know that
\begin{enumerate}
\item $\operatorname{In}(\omega_{s+1})=\operatorname{In}(c_j\omega_j)$.
\item $\nu_C(c_j\omega_j)\leq \nu_C(c_\ell\omega_\ell)$, for $\ell=-1,0,1,\ldots,j$.
\end{enumerate}
Note that
$
g^{\widetilde{\mathbf f}^1,f^1}_\ell= g^{\widetilde{\mathbf f},f}_\ell-G_{s+1}c_\ell$, for $\ell=-1,0,1,\ldots,j
$.
Then, in order to show that we have a good list, it is enough to show that $\nu_C(G_{s+1}c_j\omega_j)>\tilde u_j$. Let us do it.

We know that $\nu_D(g_j\omega_j)=\nu_D(G_{s+1}c_j\omega_j)<nm$, since they share initial part. Noting that the divisorial values are under $nm$, we have that
$$
\nu_D(g_j)=\nu_C(g_j),\quad \nu_D(G_{s+1}c_j)=\nu_C(G_{s+1}c_j).
$$
We conclude that
$
\nu_C(G_{s+1}c_j\omega_j)=\nu_C(g_j\omega_j)>\tilde u_j
$,
as desired.
\end{proof}

\section{New Discrete Analytic Invariants}\label{sec:example}

Let $\pi:M\rightarrow ( {\mathbb C^2,0})$ be the minimal reduction of singularities of a cusp  with Puiseux pair $(n,m)$. We know that $\pi$ is the composition
$$
\pi=\pi_1\circ\pi_2\circ\cdots\pi_N
$$
of blowing-ups $\pi_j:M_{j}\rightarrow M_{j-1}$ centered at points $P_{j-1}\in M_{j-1}$,  for $j=1,2,\ldots,N$, where $P_0=0\in \mathbb C^2$. Hence $M_0=(\mathbb C^2,\mathbf 0)$ and $M_N=M$. We also know that each infinitely near point $P_j$ belongs to the divisor
$$
D_{j}=\pi_j^{-1}(P_{j-1}), \qquad j=1,2,\ldots,N-1.
$$
We also put $D=D_N$ the last divisor of $\pi$.
Let us denote by $\mathcal C_\pi$ the set of all cusps $C$ such that $\pi$ is the minimal reduction of singularities of $C$.

\begin{remark}
	\label{rk:cotasumemultbases}
	For any Saito basis $\omega,\omega'$ of a cusp $C\in \mathcal C_\pi$, we have that
	$$
	\nu_{D_j}(\omega)+\nu_{D_j}(\omega')\leq \nu_{D_j}(xyf),
	\quad j=1,2,\ldots,N,
	$$
	where $f=0$ is a reduced equation of the cusp $C$. Indeed, since $\omega,\omega'$ is a Saito basis, we have that
	$$
	\omega\wedge\omega'=ufdx\wedge dy=uxyf\left(\frac{dx}{x}\wedge\frac{dy}{y}\right),
	$$
	where $u$ is  a unit. The property follows from the fact that $$\nu_{D_j}(\omega)+\nu_{D_j}(\omega')\leq \nu_{D_j}(\omega\wedge \omega').$$
\end{remark}

Given a divisor $D_j$, for $j=1,2,\ldots,N$, and a cusp $C\in \mathcal C_\pi$, we define {\em the pair $(\mathfrak{s}_{D_j}(C),\widetilde{\mathfrak s}_{D_j}(C))$ of Saito multiplicities at $D_j$} by
\begin{eqnarray}
\mathfrak s_{D_j}(C)&=&\min\{\nu_{D_j}(\omega); \; \omega \text{ belongs to a Saito basis of } C\}.\\
\widetilde{\mathfrak s}_{D_j
}(C)&=&\max\{\nu_{D_j}(\omega);\; \omega \text{ belongs to a Saito basis of } C\}.
\end{eqnarray}
Note that $\mathfrak s_{D_j}(C)$ is equal to the minimal divisorial order of the elements of any Saito basis, whereas $\widetilde{\mathfrak s}_{D_j}(C)$ does not follow directly from a given Saito basis.

The pair of Saito multiplicities is an analytic invariant of the cusp $C$. In \cite{Genzmer2}, the author introduces an invariant directly related with  the first pair $(\mathfrak s_{D_1}(C),\widetilde{\mathfrak s}_{D_1}(C))$.

A natural question is to know if the pairs of Saito multiplicities may be deduced from the knowledge of the semimodule of differential values. The answer is positive for the last pair $(\mathfrak s_D(C),\widetilde{\mathfrak s}_D(C))$. On the other hand, we present here an example of two cusps in $\mathcal C_\pi$ having the same semimodule of differential values such that the first pairs of Saito multiplicities do not coincide.

\begin{theorem} Take $C\in \mathcal C_\pi$, then
	$
	(\mathfrak s_D(C),\widetilde{\mathfrak s}_D(C))=(t_{s+1},\tilde t_{s+1})
	$,
	where $t_{s+1}$ and $\tilde t_{s+1}$ are the last critical values of the semimodule of differential values of $C$.
\end{theorem}
\begin{proof}We know that $\omega_{s+1}$ and $\widetilde{\omega}_{s+1}$ is a Saito basis of $C$ and
	$$
	\nu_D(\omega_{s+1})=t_{s+1}< \tilde t_{s+1}=\nu_D(\tilde \omega_{s+1}).
	$$
This proves that $\mathfrak s_D(C)=t_{s+1}$ and $\tilde t_{s+1}\leq \widetilde{\mathfrak s}_D(C)$. Now, let $\omega,\omega'$ be another Saito basis, with $\nu_D(\omega)=t_{s+1}$ and  $\nu_D(\omega')\geq \nu_D(\tilde \omega_{s+1})=\tilde t_{s+1}$. Let us write
$$
\omega=h\omega_{s+1}+\tilde h\tilde\omega_{s+1},\quad
\omega'=g\omega_{s+1}+\tilde g\tilde\omega_{s+1},
$$
where $\delta=h\tilde g-g\tilde h$ is a unit in $\mathcal O_{\mathbb C^2,\mathbf 0}$. By taking in consideration the divisorial order $\nu_D$, we have that  $\nu_D(h)=0$ and $\nu_D(g)>0$; hence $h$ is a unit and $g$ is not a unit. Since $\delta $ is a unit, we have that $\tilde g$ is a unit. If $\nu_D(\omega')>\tilde t_{s+1}=\nu_D(\tilde \omega_{s+1})$, we necessarily have that
$$
\nu_D(g\omega_{s+1})=\nu_D(\tilde g\tilde \omega_{s+1})= \nu_D(\tilde \omega_{s+1})=\tilde t_{s+1}.
$$
Let us see that this is not possible. Assume that $t_{s+1}=t_{s}+n\ell_{s+1}^n$ and hence
$\tilde t_{s+1}=t_{s}+m\ell_{s+1}^m$ (the case $t_{s+1}=t_{s}+m\ell_{s+1}^m$ runs in a similar way). We have
$$
\nu_D(g)+t_{s+1}=\tilde t_{s+1}\Rightarrow \nu_D(g)+ n\ell_{s+1}^n=m\ell_{s+1}^m.
$$
Noting that $\nu_D(g)\in \Gamma$, we obtain two different ways of writing  $m\ell_{s+1}^m<nm$ as a linear combination of $n,m$ with non-negative integer coefficients. This is a contradiction.
\end{proof}

We are going now to present the example of two cusps $C_1$ and $C_2$ corresponding to the Puiseux pair $(7,36)$, such that the (common) semimodule of differential values has a basis $\mathcal B=(7,36,123)$ and such that the Saito pairs of multiplicities with respect to the first divisor $D_1$ are different for $C_1$ and $C_2$.

\begin{remark}
	Let us note that for any $1$-form $\omega$, we have that
	$$
	\nu_\mathbf 0(\omega)=\nu_{D_1}(\omega)-1,
	$$
	where $\nu_\mathbf 0$ means the minimum of the multiplicity of the coefficients of $\omega$.
\end{remark}

\bigskip
\noindent{\em First example:}
Consider the cusp $C_1$ invariant to the 1-form
$$
\omega=36x^3(7xdy-36ydx)-560y^3dy,
$$
with a parametrization $\phi_1(t)=(t^7,t^{36}+t^{116}+\tfrac{28}{9}t^{196}+h.o.t.)$. The basis of  semimodule of differential values of $C_1$ is $(7,36,123)$, with minimal standard basis
$$
\mathcal S=(\omega_{-1}=dx,\omega_0=dy,\omega_1=7xdy-36ydx).
$$
We have $u_2^n=\lambda_1+n\ell_2^n=\lambda_0+mb_2$, that is $123+7\ell_2^n=36+36b_2$, we obtain that
$$\ell_2^n=b_2=3,\quad u_2^n=144.
$$
Similarly, we found out that
$$
u_2^m=231=123+36\ell_2^m=7+7a_2,
\quad \ell_2^m=3,\quad a_2=32.
$$
Hence $u_2=u_2^n$ and $\tilde u_2=u_2^m$. Moreover, we have
 $$
 t_2=t_2^n=t_1+n\ell_2^n=43+7\cdot 3=64,\quad \tilde t_2=t_2^m=t_1+m\ell_2^m=43+36\cdot 3=151.
 $$
 We see that
 $\nu_D(\omega)=t_2=64$. Hence we can take $\omega_2=\omega$ to obtain an extended standard basis and as being one of the generators of a Saito basis of $C_1$. Notice that $\nu_{D_1}(\omega)=4$, since $\nu_\mathbf 0(\omega)=3$.
We can take $\tilde \omega_2$ to be  a 1-form with divisorial order $\nu_D(\tilde\omega_2)=\tilde t_2=151$ and $C_1$ being invariant by $\tilde \omega_2$. By Delorme's decomposition in  Theorem \ref{teo:standard:decomposition3}, we can  write $\tilde{\omega}_2$ as
$$
\tilde{\omega}_2=y^3\omega_1+\mu x^{32}dx+\eta_2; \qquad\eta_2=f_{-1}dx+f_0dy+f_1(7xdy-36ydx),
$$
for an appropriate constant $\mu$ and such that  $\nu_{C_1}(f_\ell\omega_\ell)>\tilde{u}_2=231$, for $\ell=-1,0,1$.

Let us compute $\nu_{D_1}(\tilde \omega_2)$. Assume that we have  $\nu_{D_1}(f_\ell\omega_\ell)>5$, for $\ell=-1,0,1$, then we obtain that
$
\nu_{D_1}(\tilde \omega_2)=5
$.
In view of Remark \ref{rk:cotasumemultbases}, we know that
$$
\mathfrak{s}_{D_1}(C_1)+\widetilde{\mathfrak{s}}_{D_1}(C_1)\leq \nu_{D_1}(xyf)=7+2=9,
$$
Thus, we have $({\mathfrak{s}}_{D_1}(C_1),\widetilde{\mathfrak{s}}_{D_1}(C_1))=(4,5)$ since the Saito basis $\omega,\tilde\omega_2$ gives the maximal pair $(4,5)$.

It remains to show that
$\nu_{D_1}(f_\ell\omega_\ell)>5$, for $\ell=-1,0,1$. We consider two situations; $\nu_{D}(f_\ell)\geq nm$ and $\nu_D(f_\ell)<nm$. In the first situation we have that
$$
\nu_\mathbf 0(f_\ell)\geq n=7.
$$
In the case that $\nu_D(f_\ell)<nm$ we have that
$$
\nu_D(f_\ell)=\nu_{C_1}(f_\ell)>231-\lambda_\ell.
$$
Moreover, looking at the monomials in the expression of $f_\ell$, we have that
$$
\nu_D(f_\ell)\leq \nu_\mathbf 0(f_\ell) m=36  \nu_\mathbf 0(f_\ell).
$$
Thus we have:
$$
\nu_{D_1}(f_\ell\omega_\ell)=
\left\{
\begin{array}{lcrr}
	\nu_0(f_{-1})+1\geq \frac{\nu_D(f_{-1})}{36}+1&>&\frac{231-\lambda_{-1}}{36}+1=\frac{260}{36}\geq 5;& \ell=-1.\\
	\nu_0(f_{0})+1\geq
	\frac{\nu_D(f_{0})}{36}+1&>&\frac{231-\lambda_{0}}{36}+1
	=\frac{231}{36}\geq 5;& \ell=0.\\
		\nu_0(f_{1})+2\geq
	\frac{\nu_D(f_{1})}{36}+2&>&\frac{231-\lambda_{1}}{36}+2=
	\frac{180}{36}= 5;& \ell=1.
\end{array}
\right.
$$

\bigskip
\noindent{\em Second example:}
Take the cusp $C_2$ with Puiseux pair $(7,36)$ invariant by the $1$-form
\begin{align*}
	\omega'=36x^3(7xdy-36ydx)-560y^3dy+y(7xdy-36ydx).
\end{align*}
and defined by a parametrization as follows
$$
\phi_2(t)=(t^7,t^{36}+t^{116}-\tfrac{4}{171}t^{131} + \tfrac{1}{1782}t^{146} -\tfrac{1}{72900}t^{161}+h.o.t.).
$$
The basis of the semimodule of differential values is $(7,36,123)$. We can take
$$
\mathcal S=(\omega_{-1}=dx,\omega_0=dy,\omega_1=7xdy-36ydx).
$$
as minimal standard basis for $C_2$ (thus, it is the same one as for $C_1$). We repeat the arguments as for $C_1$. Namely, we can take $\omega_2'=\omega'$ as one of the generators of a Saito basis of $C_2$, with $\nu_{D}(\omega')=t_2$. Again, we obtain a partial standard system $(\omega_{-1},\omega_0,\omega_1,\omega'_2=\omega',\tilde\omega'_2)$, where $\tilde{\omega}_2'$ can be written as
\begin{align*}
\tilde{\omega}_2'=y^3\omega_1+\mu' x^{32}dx+\eta_2'; \qquad\eta_2'=\sum_{\ell=-1}^1f_\ell'\omega_\ell,
\end{align*}
with $\mu'$ being and appropriate constant and $\nu_{C_2}(f_\ell'\omega_\ell)>231$. Thus, again we found that $\nu_{D_1}(f_\ell'\omega_\ell)> 5$. We have that $\nu_{D_1}(\tilde{\omega}_2')=5$.

Now, we have that
$$(\nu_{D_1}(\omega'),\nu_{D_1}(\tilde\omega'_2) )=(3,5).
$$
This implies that $\mathfrak{s}_{D_1}(C_2)=3<4=\mathfrak{s}_{D_1}(C_1)$. Hence the Saito pairs of multiplicities for $C_1$ and $C_2$ are different.

Moreover, the pair $(3,5)$ is not maximal yet: the 1-form $\eta=\tilde{\omega}_2'-y^2\omega_2'$ satisfies that
 $\{\eta,\omega_2'\}$ is a Saito basis and
$
\nu_{D_1}(\eta)=6
$.
Hence the Saito's pair of multiplicities for the first divisor and the cusp $C_2$ is equal to $({\mathfrak{s}}_{D_1}(C_2),\widetilde{\mathfrak{s}}_{D_1}(C_2))=(3,6)$.
\begin{remark}
	Given a plane curve $S\subset (\mathbb C^2,\mathbf 0)$ and a finite sequence of blowing-ups
	$$
	\sigma:M\rightarrow (\mathbb C^2,\mathbf 0)
	$$
	with last exceptional divisor $D$, we can define in a similar way the Saito's pair of multiplicities $(\mathfrak{s}_D(S),\widetilde{\mathfrak{s}}_D(S))$. In this way we have infinitely many analytic invariants of $S$. An interesting question should be to describe the set of these invariants as a subset of the moduli of plane curves given in \cite{hefez2}.
\end{remark}

\end{document}